\documentclass[notitlepage,leqno,10pt]{article}
\usepackage{latexsym}
\usepackage{amsmath}
\usepackage{amsfonts}
\usepackage{amssymb}
\usepackage{amscd}
\usepackage[all]{xy}
\usepackage{graphicx,amssymb,amsmath,epsfig}
\renewcommand{\a }{\alpha }
\renewcommand{\b }{\beta }
\renewcommand{\d}{\delta }
\newcommand{\D }{\Delta }

\newcommand{\e }{\varepsilon }

\renewcommand{\i }{\iota}
\newcommand{\G }{\Gamma }
\renewcommand{\l }{\lambda }
\renewcommand{\L }{\Lambda }
\newcommand{\m }{\mu }
\newcommand{\n }{\nabla }

\newcommand{\Sig }{\Sigma}

\renewcommand{\O }{\Omega }

\newcommand{\ov}{\overline}
\newcommand{\intbar}{\mathop{\int\makebox(-13.5,0){\rule[4pt]{.7em}{0.3pt}}%
\kern-6pt}\nolimits}

\newcommand{\be}{\begin{equation}}
\newcommand{\ee}{\end{equation}}
\newcommand{\bes}{\begin{equation*}}
\newcommand{\ees}{\end{equation*}}
\newcommand{\ba}{\begin{eqnarray}}
\newcommand{\ea}{\end{eqnarray}}
\newcommand{\bas}{\begin{eqnarray*}}
\newcommand{\eas}{\end{eqnarray*}}
\newenvironment{pf}{\noindent{\sc Proof}.\enspace}{\rule{2mm}{2mm}\medskip}
\newenvironment{pfn}{\noindent{\sc Proof}}{\rule{2mm}{2mm}\medskip}

\newcommand{\R}{\mathbb{R}}

\newcommand{\Z}{\mathbb{Z}}

\newcommand{\N}{\mathbb{N}}

\newcommand{\bde}{\begin{definition}\rm}
\newcommand{\ede}{\end{definition}\rm}
\newcommand{\bex}{\begin{example}\rm}
\newcommand{\eex}{\end{example}}

\def\sj#1{\hbox{Sym}^{*{#1}}}
\def\ob#1{\overline{B}_{#1}}

\author{ Mohameden  AHMEDOU,  \; \; \; Sadok KALLEL,   \; \;\;   Cheikh Birahim NDIAYE  }

\date{}

\title{\bf The resonant boundary $Q$-curvature problem  \\
and   boundary-weighted barycenters}

\begin{document}

\newtheorem{lem}{Lemma}[section]
\newtheorem{pro}[lem]{Proposition}
\newtheorem{thm}[lem]{Theorem}
\newtheorem{rem}[lem]{Remark}
\newtheorem{cor}[lem]{Corollary}
\newtheorem{df}[lem]{Definition}

\maketitle


\centerline{   \emph{In fond memory of Abbas Bahri }}

\begin{center}
{\bf Abstract}

\end{center}
Given a compact four-dimensional Riemannian manifold $(M, g)$ with boundary, we study the problem of existence of Riemannian metrics on $M$ conformal to $g$ with prescribed $Q$-curvature in the interior $\mathring{M}$  of $M$, and zero $T$-curvature and mean curvature on the boundary $\partial M$ of $M$. This geometric problem is equivalent to solving a fourth-order elliptic boundary value problem (BVP) involving the Paneitz operator with boundary conditions of  Chang-Qing and Neumann operators. The corresponding BVP has a variational formulation but  the corresponding variational problem, in the case under study, is not compact. To overcome such a difficulty we perform  a systematic study, \`a  la Bahri, of the so called   {\it critical points at infinity }, compute their Morse indices, determine their contribution to the difference of topology between the sublevel sets of  associated Euler-Lagrange   functional and hence extend the  full Morse Theory to this noncompact variational problem. To establish  Morse inequalities we were led to investigate from the topological viewpoint the space of  boundary-weighted barycenters of the underlying manifold, which  arise in the description of the topology of very negative sublevel sets of the related functional. As an application of our approach we derive various existence results and provide a Poincar\'e-Hopf type criterium for  the prescribed $Q$-curvature problem on compact  four dimensional Riemannian  manifolds with boundary.

\begin{center}

\bigskip

\centerline{\bf AMS subject classification:  53C21, 35C60, 58J60, 55R80.}

\end{center}

\tableofcontents

\section{Introduction and statement of the main results}

On a four dimensional Riemannian manifold $(M, g)$, Paneitz\cite{p1} discovered in 1983, a conformally covariant  fourth order  differential operator denoted by $P^4_g$ and called the Paneitz operator. Brenson  and Oersted\cite{bo} associated to  this operator a natural concept of curvature called $Q$-curvature and denoted by $Q_g$. Both the Paneitz operator and the $Q$-curvature are defined in terms of the Ricci tensor\;$Ric_{g}$\;and the scalar curvature $R_{g}$\;of the Riemannian manifold\;$(M , g)$\; as follows:
\begin{equation*}
P^4_g=\D_{g}^{2}+div_{g}\left((\frac{2}{3}R_{g}g-2Ric_{g})\n_g\right),\;\;\;\;\;\;\;Q_g=-\frac{1}{12}(\D_{g}R_{g}-R_{g}^{2}+3|Ric_{g}|^{2}),
\end{equation*}
where \;$\n_g$\;is the covariant derivative with respect to $g$.
\vspace{4pt}

\noindent
Likewise Chang-Qing\cite{cq1} have discovered an operator\;$P^3_g$\; which is associated to the boundary of a compact four-dimensional Riemannian manifold $(M, g)$ with boundary and a third-order curvature \;$T_g$\;naturally associated to \;$P^3_g$. They are defined as follows
\begin{equation*}
P^3_g=\frac{1}{2}\frac{\partial {\D_g}}{\partial n_g}+\D_{\hat g}\frac{\partial }{\partial n_g}-2H_g\D_{\hat g}+L_g(\nabla_{\hat g},\nabla_{\hat g})+\nabla_{\hat g}H_g.\nabla_{\hat g}+(F_g-\frac{R_g}{3})\frac{\partial }{\partial n_g}.
\end{equation*}
\begin{equation*}
T_g=-\frac{1}{12}\frac{\partial R_g}{\partial n_g}+\frac{1}{2}R_gH_g-<G_g,L_g>+3H_g^3-\frac{1}{3}tr_g(L_g^3)-\D_{\hat g}H_g,
\end{equation*}
where $\hat g$\; is the metric induced by \;$g$\;on\;\;$\partial M$,  $\frac{\partial }{\partial n_g}$ is the inward Neumann operator on $\partial M$ with respect to $g$,\;$L_g=(L_{g,ab})=-\frac{1}{2}\frac{\partial g_{ab}}{\partial n_g}$ is the second fundamental form of \;$\partial M$ with respect to $g$,\;\;$H_g=\frac{1}{3}tr_g(L_g)=\frac{1}{3}\hat g^{ab}L_{g,ab}$\;($\hat g^{ab}$\;are the entries of the inverse \;$\hat g^{-1}$\;of the metric\; $\hat g$)\;is the mean curvature of \;$\partial M$,\;$R^k_{g,ijl}$\;is the Riemann curvature tensor of \;$(M, g)$\;\; $R_{g,ijkl}=g_{mi}R^{m}_{g,jkl}$\;($g_{ij}$\;are the entries of the metric \;$g$),\;\;$F_g=R^{a}_{g,nan}$ (with $n$ denoting the index corresponding to the normal direction in local coordinates)\;\;\;and\;\;$<G_g,L_g>=\hat g^{ac}\hat g^{bd}R_{g,anbn}L_{g,cd}$. Moreover, the notation $L_g(\nabla_{\hat g},\nabla_{\hat g})$, means $L_g(\nabla_{\hat g},\nabla_{\hat g})(u)=\n_{\hat g}^a(L_{g,ab}\n_{\hat g}^bu)$. We point out that in all those notations above $i,j,k,l=1,\cdots 4$ and $a,b,c,d=1, \cdots, 3$, and Einstein summation convention is used for repeated indices.
\vspace{6pt}

\noindent
As the Laplace-Beltrami operator and the Neumann operator are conformally covariant, we have that\;$P^4_g$ is conformally covariant of bidegree $(0, 4)$ and $P^3_g$ of bidegree $(0, 3)$. Furthermore, as they govern the transformation laws of the Gauss curvature and the geodesic curvature on compact surfaces with boundary, the couple $(P^4_g, P^3_g)$ does the same for\;$(Q_g,T_g)$\;on a compact four-dimensional Riemannian manifold with boundary $(M, g)$. In fact, under a conformal change of metric\;$ g_u=e^{2u}g$, we have
\begin{equation}\label{eq:law}
\left\{
    \begin{split}
P^4_{g_u}=e^{-4u}P^4_g, \mbox{ in }  \mathring{M},\\
P^3_{g_u}=e^{-3u}P^3_{g}, \mbox{ on }  \partial M.
    \end{split}
  \right.
\qquad \mbox{and}\qquad
\left\{
\begin{split}
P^4_gu+2Q_g=2Q_{ g_u}e^{4u}\;\;\text{in }\;\; \mathring{M},\\
P^3_gu+T_g=T_{ g_u}e^{3u}\;\;\text{on}\;\;\partial M.
\end{split}
\right.
\end{equation}
Besides the above  analogy, we have also an extension of the Gauss-Bonnet identity,  known as the Gauss-Bonnet-Chern formula
\begin{equation}\label{eq:gbc}
\int_{M}(Q_{g}+\frac{|W_{g}|^{2}}{8})dV_{g}+\oint_{\partial M}(T_g+Z_g)dS_g=4\pi^{2}\chi(M),
\end{equation}
where\;$W_g$\;denote the Weyl tensor of\;$(M,g)$\;and\;$Z_g$ is given by the following formula
$$
Z_g=R_gH_g-3H_gRic_{g,nn}+\hat g^{ac}\hat g^{bd}R_{g,anbn}L_{g,cd}-\hat g^{ac}\hat g^{bd}R_{g,acbc}L_{g,cd}+6H_g^3-3H_g|L_g|^2+tr_g(L_g^3),
$$
with $tr_g$ denoting the trace with respect to the metric induced on $\partial M$ by $g$ (namely $\hat g$) and $\chi(M)$ the Euler  characteristic of $M$. Concerning the quantity $Z_g$, we have that it vanishes when the boundary is totally geodesic and $\oint_{\partial M}Z_gdV_g$ is always conformally invariant, see \cite{cq1}. Thus, setting
\begin{equation}\label{eq:invc}
\kappa_{(P^4,P^3)}:=\kappa_{(P^4,P^3)}[g]:=\int_{M}Q_gdV_g+\oint_{\partial M}T_gdS_{g},
\end{equation}
we have that thanks to\;$\eqref{eq:gbc}$, and to the fact that \;$|W_g|^2dV_g$ is pointwise conformally invariant, $\kappa_{(P^4,P^3)}$\;is a conformal invariant.  We remark that $4\pi^2$ is the total integral of the $(Q, T)$-curvature of the standard four-dimensional hemisphere.
\vspace{8pt}

\noindent
As  was addressed  in \cite{no1} by the first and third authors, a natural question is  the following: given a compact four-dimensional Riemannian manifold with boundary\;$(M, g)$\; and $K: M\longrightarrow \R_+$ smooth, under which conditions on $K$ does $M$ carries a Riemannian metric conformal to $g$\;for which the corresponding \;$Q$-curvature is $K$ and the corresponding \;$T$-curvature and mean curvature vanishes. Related questions have been studied for compact Riemannian surfaces with boundary regarding Gauss curvature and geodesic curvature, see for example \cite{bre1}, \cite{bre2}, \cite{changliu}, \cite{cherrier}, and \cite{ops}, and for compact Riemannian manifolds with boundary of dimension bigger that $2$ regarding scalar curvature and mean curvature, see  for example, \cite{amy}, \cite{bre3}, \cite{bs}, \cite{cherrier}, \cite{dmo}, \cite{es1}, \cite{es2}, \cite{HL1}, \cite{HL2} and \cite{marques} and the references therein. Thanks to \eqref{eq:law}, this problem is equivalent to finding a smooth solution to the following BVP:
\begin{equation}\label{eq:bvpg}
\left\{
\begin{split}
P^4_gu+2Q_g&=2Ke^{4u}\;\;&\text{in}\;\;\mathring{M},\\
P^3_gu+T_g&=0\;\;&\text{on}\;\;\partial M,\\
\frac{\partial u}{\partial n_g}-H_gu&=0\;\;&\text{on}\;\;\partial M,
\end{split}
\right.
\end{equation}
where\;$\frac{\partial}{\partial n_g}$\;is the inward normal derivative with respect to \;$g$.
 Since the problem is conformally invariant, it is not restrictive to assume\;$H_g=0$, since this can be always achieved  through a conformal transformation of the background metric. Thus, from now on, we will assume that we are working with a background metric $g$ satisfying $H_g=0$ and hence BVP \eqref{eq:bvpg} becomes the following one with Neumann homogeneous boundary condition:
\begin{equation}\label{eq:bvps}
\left\{
\begin{split}
P^4_gu+2Q_g&=2Ke^{4u}\;\;&\text{in}\;\; \mathring{M},\\
P^3_gu+T_g&=0\;\;&\text{on}\;\;\partial M,\\
\frac{\partial u}{\partial n_g}&=0\;\;&\text{on}\;\;\partial M.
\end{split}
\right.
\end{equation}
Defining \;$\mathcal{H}_{\frac{\partial}{\partial n}}$\;as
\begin{equation*}
\mathcal{H}_{\frac{\partial}{\partial n}}=\Big\{u\in W^{2, 2}(M):\;\;\;\frac{\partial u}{\partial n_g}=0\;\;\;\text{on}\;\;\partial M\Big\},
\end{equation*}
where $W^{2, 2}(M)$ denotes the space of functions on $M$ which are square integrable together with their first and second derivatives, and \;$P^{4,3}_g$\;as follows, for every \;$u,v\in\mathcal{H}_{\frac{\partial }{\partial n}}$
\begin{equation*}
\begin{split}
\left<P^{4,3}_gu,v\right>_{L^2(M)}=\int_{M}\left(\D_g u\D_gv+\frac{2}{3}R_g\nabla_g u\cdot\nabla_g v\right)dV_g-2\int_{M}Ric_g(\nabla_g u,\nabla_g v)dV_g\\-2\oint_{\partial M}L_g( \nabla_{\hat g} u, \nabla_{\hat g} v)dS_g.
\end{split}
\end{equation*}
It follows from standard  elliptic regularity theory that smooth solutions to \eqref{eq:bvps} can be found by looking at critical points of the geometric functional
\begin{equation*}
\begin{split}
II(u)=\left<P^{4,3}u,u\right>_{L^2(M)}+4\int_{M}Q_gudV_g+4\oint_{\partial M}T_gudS_g-
\kappa_{(P^4,P^3)}\ln \int_{M}e^{4u}dV_g,\;\;\;u\in \mathcal{H}_{\frac{\partial }{\partial n}}.
\end{split}
\end{equation*}
\vspace{4pt}

\noindent
 Similarly to the case of closed four-dimensional Riemannian manifolds, the spectral properties of $P^{4,3}_g$ and the sign of $\kappa_{(P^4, P^3)}$ are strongly related. In fact, Catino-Ndiaye\cite{cn} proved that if $\partial M$ is umbilic in $(M , g)$, the Yamabe invariant $Y(M, \partial M, [g]):=\inf\{\int_{M}R_{g_u}dV_{g_u}+\oint_{\partial M}H_udS_{g_u}, g_u=e^{2u}g, \int_{M}e^{4u}dV_g=1\}$  is positive and $\kappa_{(P^4, P^3)}+\frac{1}{6}Y(M, \partial M, [g])^2>0$, then $ker P^{4, 3}_g\simeq \R$ and  $P^{4, 3}_g$ is also nonnegative. They  also observed that $\kappa_{(P^4, P^3)}$  satisfies a  rigidity type result, namely that if $\partial M$ is umbilic in $(M , g)$ and $Y(M, \partial M, [g])\geq 0$, then $\kappa_{(P^4, P^3)}\leq 4\pi^2$ with the equality holding if and only if $(M, g)$ is conformally equivalent to the standard four-dimensional hemisphere $\mathbb{S}^4_+$. We point out that, the latter rigidity result has been noticed by Chen\cite{cs}.
\vspace{6pt}

\noindent
As in the case of closed four-dimensional Riemannian manifolds, here also, the analytic features of $II$ depend strongly on the values taken by  $\kappa_{(P^4, P^3)}$. Indeed depending on whether $\kappa_{(P^4, P^3)}$ is a multiple of $4\pi^2$ or not, the way of finding critical points of $II$ changes drastically. To the best of our knowledge the first existence result for problem $\eqref{eq:bvps}$ has been obtained by Chang and Qing, see \cite{cq2}, under the assumptions that $P^{4,3}_{g}$ is nonnegative, $ker P^{4,3}_{g}\simeq \R$ and $\kappa_{(P^4,P^3)}<4\pi^2$. An alternative proof using geometric flows method has been given in \cite{nd4}. As already mentioned, a first sufficient condition to ensure those hypotheses (in the umbilic case) was given by Catino-Ndiaye\cite{cn}. Later, Ndiaye\cite{nd2} developed a variant of the min-max scheme of Djadli-Malchiodi\cite{dm} to extend the result of Chang-Qing\cite{cq2}. Precisely, he showed that problem $\eqref{eq:bvps}$ is solvable provided that $ker P^{4, 3}_g \simeq \R$ and $\kappa_{(P^4, P^3)}$ is not a positive integer multiple of $4\pi^2$.
\vspace{4pt}

\noindent
As in the case of closed four-dimensional Riemannian manifolds, here also, the assumptions $ker P^{4, 3}_g\simeq \R$ and $\kappa_{(P^4,P^3)}\notin 4\pi^2\N^{*}$ will be referred to as {\em nonresonant} case. This terminology is motivated by the fact that in that situation the set of solutions to some appropriate perturbations of BVP $\eqref{eq:bvps}$ (including it) is compact, see \cite{nd2}.  Naturally, we call {\em resonant} case when $ker P^{4, 3}_g\simeq \R$ and $\kappa_{(P^4, P^3)}\in 4\pi^2\N^{*}$. We divide the {\em resonant} case in two subcases. Precisely, we call the situation where $\kappa_{(P^4,P^3)}=4\pi^2$ the {\em critical case} and when $\kappa_{(P^4,P^3)}\in 4\pi^2(\N^{*}\setminus \{1\})$ the {\em supercritical} one. With these terminologies, we have that the works of Chang-Qing\cite{cq2} and Ndiaye\cite{nd2} answer affirmatively the question raised  above in the {\em non resonant} case. Unlike the { \it non resonant case}, up to the knowledge of the authors, there are  no existence results in the {\em resonant} one...
\vspace{4pt}

\noindent
To give some motivations of the study of the $(Q, T)$-curvature, we discuss some geometric applications of it.  We have two results proven by Chen\cite{cs}, and Catino-Ndiaye\cite{cn}.  The first one follows from the works of Chen\cite{cs} and Catino-Ndiaye\cite{cn} and says that if  the Yamabe invariant $Y(M,\partial M,[g])$ and $\kappa_{(P^4,P^3)}$ are both positive and \;$(M, g)$\; has umbilic boundary, then $M$ carries a conformal metric with positive Ricci curvature. Hence\;$M$\;has a finite fundamental group. The second one due to Catino-Ndiaye\cite{cn} says that if  $(M, g)$    has umbilic boundary,    $Y(M,\partial M,[g])>0$, and $\kappa_{(P^4, P^3)}>\frac{1}{8}\int_M|W_g|^2dV_g$, then $M$ admits a Riemannian metric   $\bar g$   such that $(M, \bar g)$ has constant positive sectional curvature and $\partial M$ is totally geodesic in $(M, \bar g)$. We remark also that the pair Paneitz, Chang-Qing operator, and the $(Q, T)$-curvature appear in the study of log-determinant formulas, Gauss-Bonnet type formulas, and the compactification of some locally conformally flat four-dimensional manifolds, see  \cite{cqy1}, \cite{cqy2}, \cite{cq1}, \cite{cq2}.
\vspace{4pt}

\noindent
In this paper, we are interested in the {\em  resonant} case, namely when $\ker P^{4, 3}_g\simeq \R $ and $\kappa_{(P^4, P^3)}=4k\pi^2$ with $k\in \N^*$. Namely we first completely identify  {\it  the critical points at infinity } of $II$, compute their Morse indices  and determine  their topological  contribution to the difference of topology between the sublevel sets.  Next, we combine the variational contribution of the critical points at infinity of $II$  with classical tools of Morse theory and  precise knowledge of the topology of the boundary-weighted barycenters  $B_*^{\partial}(M)$ (whose relevance in the  problem under study was first discovered by Ndiaye\cite{nd2}) to prove strong Morse inequalities and provide various type of existence results.


\vspace{4pt}

\noindent
To state  our main existence results, we  need to fix  some notation and make some definitions. For every $(p, q)\in (\N^*)^2$, we define \;$\mathcal{F}_{p}^M: (\mathring{M} ^p)^* \longrightarrow \R$\; as follows
\begin{equation}\label{eq:limitfsint}
\mathcal{F}_{p}^M(a_1,\cdots, a_p):=\sum_{i=1}^p\left(H(a_i, a_i)+\sum_{j\neq i}G(a_i, a_j)+\frac{1}{2}\ln(K(a_i))\right),
\end{equation}
and $\mathcal{F}_{q}^{\partial M}: (\partial M ^q)^* \longrightarrow \R$\; as follows
\begin{equation}\label{eq:limitfsbound}
\mathcal{F}_{q}^{\partial M}(a_1,\cdots, a_q):=\sum_{i=1}^q\left(H(a_i, a_i)+\sum_{j\neq i}G(a_i, a_j)+\ln(K(a_i))\right),
\end{equation}
where $$(\mathring{M}^p)^*:=M^p\setminus F(\mathring{M}^p), \;\;((\partial M)^q)^*:=(\partial M)^q\setminus F((\partial M)^q)$$
with $F(\mathring{M}^p)$ and $F((\partial M)^q)$ denoting respectively the fat Diagonal of $(\mathring{M})^p$ and $( \partial M)^p$, $G$ is the Green's function of $(P^4_g(\cdot)+\frac{4}{k}Q_g, P^3_g(\cdot)+\frac{2}{k}T_g)$ under homogenous Neumann condition with respect to $g$ and  satisfying the normalization $\int_{M} Q_g(x) G(\cdot, x)dV_g(x)+\oint_{\partial M}T_g(x)G(\cdot, x)dS_g(x)=0$, and $H$ is its regular part, see Section \ref{eq:notpre} for more information.
\vspace{4pt}

\noindent
On the other hand, for $(p, q)\in \N^2$ such that $2p+q=k$, we define
\;$\mathcal{F}_{p, \;q}^M: (\mathring{M} ^p)^*\times (\partial M ^q)^* \longrightarrow \R$\; as follows
\begin{equation}\label{eq:limitfsintpq}
\mathcal{F}_{p, \;q}^M(a_1, \cdots, a_{p+q}):=\mathcal{F}_{p}^M(a_1, \cdots, a_p)+\frac{1}{2}\sum_{i=1}^p\sum_{j=p+1}^{p+q}G(a_i, a_j),
\end{equation}
and\;\; $\mathcal{F}_{p, \;q}^{\partial M}: (\mathring{M} ^p)^*\times(\partial M ^q)^* \longrightarrow \R$\; as follows
\begin{equation}\label{eq:limitfsboundpq}
\mathcal{F}_{p, \;q}^{\partial M}(a_1,\cdots, a_{p+q}):=\mathcal{F}^{\partial M}_q(a_{p+1}, \cdots, a_{p+q})+2\sum_{i=1}^p\sum_{j=p+1}^{p+q}G(a_i, a_j).
\end{equation}
Moreover, we set
\begin{equation}\label{eq:limitfsintsectint}
\mathcal{F}_{p, \;q}(A):=2\mathcal{F}_{p, q}^M(A)+\frac{1}{2}\mathcal{F}_{p, q}^{\partial M}(A)=2\mathcal{F}_p^M(A_p)+\frac{1}{2}\mathcal{F}_q^{\partial M}(A_q)+2\sum_{i=1}^p\sum_{j=p+1}^{p+q}G(a_i, a_j),
\end{equation}
with $A=(a_1, \cdots, a_{p+q})$, $A_p:=(a_1, \cdots, a_p)$, $A_q:=(a_{p+1}, \cdots, a_{p+q})$, and define
\begin{equation}\label{eq:critfk}
Crit(\mathcal{F}_{p, \;q}):=\{A\in ((\mathring{M})^p)^*\times ((\partial M)^p)^*, \;\;A\;\;\;\text{critical point of} \;\;\mathcal{F}_{p, _q}\}.
\end{equation}
Furthermore, for $A:=(a_1, \cdots, a_{p+q})\in (\mathring{M}^p)^*\times ((\partial M)^p)^*$, we set
\begin{equation}\label{eq:partiallimitint}
\mathcal{F}^A_i(x):=e^{4(H(a_i, x)+\sum_{j=1, j\neq i}^pG(a_j, x)+\frac{1}{4}\ln(K(x))+\frac{1}{2}\sum_{j=p+1}^{p+q}G(a_j, x))},\;\;\;i=1, \cdots, p,
\end{equation}
and
\begin{equation}\label{eq:partiallimitbound}
\mathcal{F}^A_i(x):=e^{4(\frac{1}{2}H(a_i, x)+\frac{1}{2}\sum_{j=p+1,j\neq i}^{p+q}G(a_j, x)+\frac{1}{4}\ln(K(x))+\sum_{j=1}^pG(a_j, x))},\;\;i=p+1, \cdots, p+q.
\end{equation}
Moreover, we set
\begin{equation}\label{eq:defindexa}
\mathcal{L}_K(A):=
\begin{cases}
\sum_{i=1}^p (\mathcal{F}^{A}_i)^{\frac{1}{4}}(a_i)L_g((\mathcal{F}^{A}_i)^{\frac{1}{4}})(a_i),\;\;&\text{if}\,\;q=0,\\\\
\sum_{i=p+1}^{p+q}(\mathcal{F}_i^A)^{\frac{1}{4}}(a_i)\frac{\partial \ln  K}{\partial n_g}(a_i),&\text{if}\;\;q\neq 0,
\end{cases}
\end{equation}
and
\begin{equation}\label{eq:minf}
i_{\infty}(A):=5p+4q-1-Morse(\mathcal{F}_{p, q}, A),
\end{equation}
where $Morse(\mathcal{F}_{p, q}, A)$ denotes the Morse index of $\mathcal{F}_{p, \;q}$ at  $A$, . We set also
\begin{equation}\label{eq:critsett1}
\begin{split}
\mathcal{F}_{\infty}^{p, q}:=\{A\in Crit(\mathcal{F}_{p, q}):\;\;\;\mathcal{L}_K(A)<0\},
\end{split}
\end{equation}
and
\begin{equation}\label{eq:critsett}
\mathcal{F}_{\infty}:=\bigcup_{(p, q)\in \N^2: 2p+q=k}\;\mathcal{F}_{\infty}^{p, q}.
\end{equation}
Furthermore, we define
\begin{equation}\label{eq:mi1}
m_i^{p, q}:=card\{A\in \mathcal{F}_{\infty}^{p, q}: i_{\infty}(A)=i\}, \, \, i=p+q-1, \cdots, 5p+4q-1,
\end{equation}
and\begin{equation}\label{eq:mi}
m_i^{k}:=\sum_{(p, q)\in \N^2:\\\;2p+q=k, p+q-1\leq i\leq 5p+4q-1}m_i^{p, q}, \, \,  i=0, \cdots, 4k-1.
\end{equation}
For $ l \in \N^*$, we use the notation $B_{l}^{\partial}(M)$ to denote the following set of  barycentric type:
\begin{equation}\label{eq:baryp}
\begin{split}
B_{l}^{\partial}(M):=\bigcup_{(p. q)\in \N^2, \;0<2p+q\leq l}B_{q}^p(M, \partial M),
\end{split}
\end{equation}
where
\begin{equation}\label{eq:barypq}
\begin{split}
B_{q}^p(M, \partial M):= \big \{\sum_{i=1}^{p+q}\alpha_i\d_{a_i}, \;\;a_i\in \mathring{M}\;\; i=1, \cdots, p, \;\;a_i\in \partial M\; \;i=p+1, \cdots, p+q,\;\\ \alpha_i\geq 0\;\;i=1, \cdots, p+q, \;\;\text{and}\;\;\sum_{i=1}^{p+q}\alpha_i=k \big \}.
\end{split}
\end{equation}
Furthermore, we define
\begin{equation}\label{eq:cpm}
c^{k-1}_n=dim \;H_n(B_{k-1}^{\partial}(M)), \;\;n=1, \cdots 4k-5,
\end{equation}
where $H_n(B_{k-1}^{\partial}(M))$ denotes the $n$-th homology group of $B_{k-1}^{\partial}(M)$ with $\Z_2$ coefficients and $dim\;H_n(B_{k-1}^{\partial}(M))$ is its dimension. For $(p, q)\in \N^2$ such that $2p+q=k$, we say
\begin{equation}\label{eq:nondeg1}
\begin{split}
&(ND)_{p, q} \;\;\;\;\;\text{holds if} \;\;\mathcal{F}_{p, q}\;\;\text{is a Morse function and for every }\;\,A\in Crit(\mathcal{F}_{p, q}), \;\;\mathcal{L}_K(A)\neq 0.
\end{split}
\end{equation}
Finally. we say
\begin{equation}\label{eq:nondeg}
(ND)\;\;\;\text{holds if }\; \;\;(ND)_{p, q}\;\;\;\text{holds for every}\;\;(p, q)\in \N^2:\; 2p+q=k.
\end{equation}
\vspace{4pt}

\noindent
Now, we are ready to state our main results, starting from those which follow from our  strong Morse type inequalities. To do so, we start with the {\em critical} case, namely when $k=1$.
\begin{thm}\label{eq:morsepoincare1}
Let $(M, g)$ be a compact four-dimensional Riemannian manifold with boundary  such that $Ker P^{4, 3}_g\simeq \R$, $H_g=0$ and $\kappa_{(P^4, P^3)}=4\pi^2$. Assuming that $K$ is a smooth positive function on $M$ such that $(ND)$ holds and the following system of non negative integers $(n_i)$
\begin{equation}\label{eq:mp1}
\begin{cases}
m^{1}_0 \, = \, 1+n_0,\\
m^{1}_i \, = \, n_i \, + \, n_{i-1}, \;&i=1, \cdots, 3,\\
0 \, = \, n_4\\
n_i\geq 0,\;\;& i=0, \cdots, 3
\end{cases}
\end{equation}
has no solutions, then  $K$ is the $Q$-curvature of a Riemannian metric conformal to $g$ with zero $T$-curvature and vanishing mean curvature on $\partial M$.
\end{thm}
\vspace{4pt}

\noindent
As a corollary of the above theorem, we derive the following Hopf-Poincar\'e  index type criterion for existence.
\begin{cor}\label{eq:existence1}
Let $(M, g)$ be a compact four-dimensional Riemannian manifold with boundary  such that $Ker P^{4, 3}_g\simeq \R$, $H_g=0$, and $\kappa_{(P^4, P^3)}=4\pi^2$. Assuming that $K$ is a smooth positive function on $M$ such that $(ND)$ holds and
\begin{equation}\label{eq:ep1}
\sum_{A\in \mathcal{F}_{\infty}} (-1)^{i_{\infty}(A)}\neq 1,
\end{equation}
then  $K$ is the $Q$-curvature of a Riemannian metric conformal to $g$ with zero $T$-curvature and vanishing mean curvature on $\partial M$.
\end{cor}
\vspace{4pt}

\noindent
As a by product of our analysis, we extend the above Hopf-Poincar\'e  index criterium to include the case where the total sum equals 1 but a partial one is not, provided that there is a jump in the Morse indices of the function $\mathcal{F}_{0, 1}$. Namely we have the following criterium:
\begin{thm}\label{t:C}
Let $(M, g)$ be a compact four-dimensional Riemannian manifold with boundary  such that $Ker P^{4, 3}_g\simeq \R$, $H_g=0$, and $\kappa_P=4\pi^2$. Assuming that $K$ is a smooth positive function on $M$ such that $(ND)$ holds and there exists a positive integer \;$1 \leq l \leq 3$\; and $A^l\in \mathcal{F}_{\infty}$ such that
\begin{equation*}\label{eq:ep1c}
\begin{split}
\sum_{A\in \mathcal{F}_{\infty},\;  i_{\infty}(A) \leq l -1 } &(-1)^{i_{\infty}(A)}\neq 1\\&\text{and}\\
 \forall A \in \mathcal{F}_{\infty},\;\; &\quad i_{\infty}(A) \neq l,
\end{split}
\end{equation*}
then $K$ is the $Q$-curvature of a Riemannian metric conformal to $g$ with zero $T$-curvature and vanishing mean curvature on $\partial M$.
\end{thm}
\vspace{4pt}

\noindent
Next, we state our main result in the {\em supercritical} case, namely when $k\geq 2$. It read as follows:
\begin{thm}\label{eq:morsepoincare2}
Let $(M, g)$ be a compact four-dimensional Riemannian manifold with boundary  such that $Ker P^{4, 3}_g\simeq \R$, $H_g=0$, and $\kappa_P=4k\pi^2$ and $k\in \N$ with $k\geq 2$. Assuming that $K$ is a smooth positive function on $M$ such that $(ND)$ holds and the following system of non negative integers $(n_i)$
\begin{equation}\label{eq:mp3}
\begin{cases}
m^k_0 \, = \, n_0,\\
m^k_1 \, = \, n_0 \, + \, n_1,\\
m_i^k \, = \, c^{k-1}_{i-1}\, + \, n_i \, + \, n_{i-1}, \;&i=2,\cdots, 4k-4,\\
m_{i}^k \, = \, n_{i} \, + \, n_{i-1},\;& i=4k-3,\cdots, 4k-1,\\
0 \, = \, n_{4k-1},\\
n_i\geq 0,\;\;& i=0, \cdots, 4k-1,
\end{cases}
\end{equation}
has no solutions, then  $K$ is the $Q$-curvature of a Riemannian metric conformal to $g$ with zero $T$-curvature and vanishing mean curvature on $\partial M$.
\end{thm}
\vspace{4pt}
\begin{rem}
The  dimension of the homology groups of the boundary-weighted barycenters $B_l^{\partial}(M)$, namely
 $$c^l_n:= \dim H_n(B_l^{\partial}(M))$$ are computed in Theorem \ref{main} of Section 5 of this paper.
\end{rem}
\noindent
As a corollary of Theorem \ref{eq:morsepoincare2}, we derive the following Hopf-Poincar\'e  index type criterion for existence.
\begin{cor}\label{eq:existence2}
Let $(M, g)$ be a compact four-dimensional Riemannian manifold with boundary  such that $Ker P^{4, 3}_g\simeq \R$, $H_g=0$, and $\kappa_P=4k\pi^2$ and $k\in \N$ with $k\geq 2$. Assume that $K$ is a smooth positive function on $M$ such that $(ND)$ holds and
\begin{equation}\label{eq:ep2}
\sum_{A\in \mathcal{F}_{\infty}} (-1)^{i_{\infty}(A)}\neq 1-\chi(B_{k-1}^{\partial}(M)),
\end{equation}
where $ \chi(B_{k-1}^{\partial}(M))$ stands for the Euler Characteristic of the space of boundary-barycenters of order $k-1$.
 Then $K$ is the $Q$-curvature of a Riemannian metric conformal to $g$ with zero $T$-curvature and vanishing mean curvature on $\partial M$.
\end{cor}
\vspace{6pt}
\begin{rem}
$\chi(B_{k-1}^{\partial}(M))$ the Euler Characteristic of the space of boundary weighted-barycenters of order $k-1$  of $M$ is  computed in Theorem \ref{main1}
  of  Section 5 of this paper.\\
\end{rem}
Just as in the one mass case, we generalize the above criterium to the case where there is a jump in the indices of the functions $\mathcal{F}_{p, q}$ for all $(p, q)\in \N^2$ such that  $2p+q=k$. Namely we prove the following:
\begin{thm}\label{eq:Cm}
Let $(M, g)$ be a compact four-dimensional Riemannian manifold with boundary such that $Ker P^{4, 3}_g\simeq \R$, $H_g=0$, and $\kappa_P=4k\pi^2$ and $k\in \N$ with $k\geq 2$. Assuming that $K$ is a smooth positive function on $M$ satisfying the non degeneracy condition $(ND)$ and there exists a positive integer \;$1 \leq l \leq 4k-1$\;  and \;$A^l\in \mathcal{F}_{\infty}$ with \;$i_{\infty}(A^l) \leq l - 1$ such that
\begin{equation*}\label{eq:ep2c}
\begin{split}
\sum_{A\in \mathcal{F}_{\infty}, \; i_{\infty}(A) \leq l - 1} (-1)^{i_{\infty}(A)}&\neq 1-\chi(B_{k-1}^{\partial}(M)),\\
&\text{and}\\
\forall A  \in \mathcal{F}_{\infty},\;\;  &\qquad i_{ \infty}(A)  \neq l,
\end{split}
\end{equation*}
 then $K$ is the $Q$-curvature of a Riemannian metric conformal to $g$ with zero $T$-curvature and vanishing mean curvature on $\partial M$.
\end{thm}
\vspace{8pt}
Taking advantage from the precise knowledge of the  location of {\em critical points at infinity}, we put  condition on the function $K$ to insure that some subcritical approximations of the critical case do not blow up and hence leading to the following existence result:

\begin{thm}\label{th:thmodd1}
Let $(M, g)$ be a compact four-dimensional Riemannian manifold with boundary  such that $Ker P^{4, 3}_g\simeq \R$, $H_g=0$, $\kappa_{(P^4, P^3)}=4\pi^2$ and $\bar k=0$. Assuming that $K$ is a smooth positive function on $M$ such that for every maximum point $a$ of $\mathcal{F}_{0, 1}$, we have $\frac{\partial K}{\partial n_g}(a)>0$, then $K$ is the $Q$-curvature of a Riemannian metric conformal to $g$ with zero $T$-curvature, vanishing mean curvature on $\partial M$, and conformal factor minimizing the Paneitz functional $II$.
 \end{thm}
 In the special case of the half-sphere, the above statement reads as follows:
 \begin{cor}\label{c:thmodd1}
  Let $(\mathbb{S}^4_+, g_{\mathbb{S}^4_+})$ be the standard four-dimensional hemisphere, $K:\mathbb{S}^4_+\longrightarrow \R_+$ a smooth positive function and $\hat K:=K_{|\mathbb{S}^3}:   \mathbb{S}^3 \longrightarrow \R_+$ its restriction on $\mathbb{S}^3$. Assuming that $\frac{\partial K}{\partial n_{g_{\mathbb{S}^4_+}}}(a)>0$ for every  maximum point $a$ of $\hat K$, then
$K$ is the $Q$-curvature of a Riemannian metric conformal to $g_{\mathbb{S}^4_+}$ with zero $T$-curvature, vanishing mean curvature on $\mathbb{S}^3$, and conformal factor minimizing the Paneitz functional $II$.
 \end{cor}
\vspace{4pt}

\noindent
Next, we present a new type of existence results based on the use of spectral information to rule out the blow up of some supercritical approximations... It read as follows:
\begin{thm}\label{eq:algtopg}
Let $(M, g)$ be a compact four-dimensional Riemannian manifold with boundary  such hat $Ker P^{4, 3}_g\simeq \R$, $H_g=0$, $\kappa_{(P^4, P^3)}=4k\pi^2$, and $k\in \N$ with  $k\geq 2$. Assuming that $K$ is a smooth positive function on $M$ satisfying the non degeneracy  $(ND)$  such that at every local minimum of $A$  of $\mathcal{F}_{0, k}$ we have that $\mathcal{L}_K(A)  < 0,$ then
$K$ is the $Q$-curvature of a Riemannian metric conformal to $g$ with zero $T$-curvature and vanishing mean curvature on $\partial M$.
\end{thm}

\begin{rem}
The above theorem has a counterpart in the closed case which improves Corollary 1.4 of \cite{nd6} by dropping the assumption on critical points of index one. However the proof, which is rather analytic, differers drastically  from the algebraic topological argument of \cite{nd6}.
\end{rem}


\vspace{4pt}

\noindent
The remainder of this paper is organized as follows. In section 2 we fix the notation used in the paper, give asymptotic of the Green's function of the  pair of operators \;$(P_{g}^4(\cdot)+\frac{4}{k}Q_g, P^3_g(\cdot)+\frac{2}{k}T_g)$ under  homogeneous Neumann boundary condition and recall the local description of blowing up solutions of some appropriate perturbations of  BVP \eqref{eq:bvps}. In section 3 we perform a refined blow up analysis around blow up points and derive a deformation Lemma  which provides an accurate description of the lost of compactness, while Section 4 is devoted to a Morse type reduction near  the {\it end} of noncompact orbits of an appropriate pseudogradient. We take advantage of such a reduction to identify the {\it critical points at Infinity} of this noncompact variational problem and compute their Morse indices. In section 5 we undertake a systematic study from the topological viewpoint of the boundary-weighted barycenters which describe the topology of very negative sublevel sets of the functional $II$, and in section 6 we provide proofs of the main results. Finally we collect in the appendix various estimates of the bubbles and their interaction as well as  refined expansions of the functional and its gradient in the neighborhood of potential critical points at Infinity.

\vspace{4pt}

\noindent
\begin{center}
{\bf Acknowledgements}
\end{center}
The  first and the third author (M.A $\&$ C-B.N) have been supported  by the DFG project "Fourth-order uniformization type theorems for $4$-dimensional Riemannian manifolds".

\section{Notation and preliminaries}\label{eq:notpre}
In this  section, we fix our notation and give some useful preliminary results.
Throughout $\N$\;denotes the set of nonnegative integers, $\N^*$\;stands for the set of positive integers, and  for $n\in \N^*$, $\R^n$ the standard $n$-dimensional Euclidean space, $\R^n_+$ the open positive half-space of $\R^n$, and $\bar \R^n_+$ its closure.\\
 For $n\in \N^*$ and $r>0$, $B^{\R^n}_0(r)$ denotes the open ball of $\R^n$ of center $0$ and radius $r$, $\bar  B^{\R^n}_0(r)$ for its closure, $B^{\R^4_+}_0(r):=B^{\R^4}_0(r)\cap\bar \R^4_+$ , and $\bar B^{\R^4_+}_0(r):=\bar B^{\R^4}_0(r)\cap\bar \R^4_+$. $\mathbb{S}^n$ denotes the unit sphere of $\R^{n+1}$, and $\mathbb{S}^n_+$ the positive spherical cap, namely $\mathbb{S}^n_+:=\mathbb{S}^n\cap \R^{n+1}_+$.
\vspace{4pt}
For a Riemmanian metric $\bar g$ on $M$ and $\hat{g}$ its  induced metric on $\partial M$,  the ball  \;$B^{\bar g}_{p}(r)$\; (respectively \;$B^{\hat g}_{p}(r)$\;) is with respect to the normal geodesic (respectively Fermi) coordinates and  similarly  \;$d_{\bar g}(x,y)$\;  denotes the geodesic distance if $x \in \mathring{M}$ (respectively the distance inherited from the Fermi coordinates). Moreover $inj_{\bar g}(M)$\;(respectively $inj_{\hat g}(M)$\;) stands for the injectivity radius. $dV_{\bar g}$\;denotes the Riemannian measure associated to the metric\;$\bar g$, and $dS_{\bar g}$ the volume form on $\partial M$ with respect to the metric induced by $\bar g$ on $\partial M$.

We will write \;$Diag(M)$\; the diagonal subspace
of \;$M^2 = M\times M$.

\vspace{4pt}
\noindent
In this paper, $(M, g)$ always refers to a given underlying compact four-dimensional Riemannian manifold with boundary $\partial M$ and interior $\mathring{M}$, and $K:  M\longrightarrow \R$ a smooth function. Furthermore we assume that:
 $$ker \;P^{4, 3}_g\simeq \R,\;\;\;\kappa_{(P^4, P^3)}=4k\pi^2\;\; \text {for some}\;\;k\in \N^*,\text{ and} \;K>0\;\;\text{on}\;\; M.$$
\vspace{4pt}

\noindent
Now, we recall  $G$  the Green's function of the operator \;$(P_{g}^4(\cdot)+\frac{4}{k}Q_g, P^3_g(\cdot)+\frac{2}{k}T_g)$\;with homogeneous Neumann boundary condition satisfying the normalization $$\int_M G(x, y)Q_g(y)dV_g(y)+\oint_{\partial M}G(x, y)T_g(y)dS_g(y)=0,\;\;\forall x\in M.$$

\vspace{14pt}


\noindent
For $1\leq p\leq \infty$ and $k\in \N$, $\beta\in  ]0, 1[$, $L^p(M)$ and $L^p(\partial M)$, $W^{k, p}(M)$, $C^k(M)$, and $C^{k, \beta} (M)$ stand respectively for the standard $p$-Lebesgue space on $M$ and $\partial M$, (k, p)-Sobolev space, $k$-continuously differentiable space and $k$-continuously differential space of H\"older exponent $\beta$, all with respect to $g$ (if the definition needs a metric structure) and for precise definitions and properties, see  for example \cite{aubin} or \cite{gt}. Given a function $u: M\longrightarrow \R$\; such that \;$u\in L^1(M)$ and $u\in L^1(\partial M)$, we define \;$\ov{u}_{(Q, T)}$ as follows\;
$$
\ov{u}_{(Q, T)}:=\frac{1}{4k\pi^2}\left(\int_M uQ_gdV_g+\oint_{\partial M} u T_gdS_g\right).
$$
\vspace{4pt}

\noindent
Given a generic Riemannian metric $\tilde g$ on $M$ and a function \;$F(x, y)$\; defined on a open subset of\;$ M^2$ which is symmetric and\; with $F(\cdot, \cdot)\in C^2$ with respect to $\tilde g$, we define $\frac{\partial F(a, a)}{\partial a}:=\frac{\partial F(x, a)}{\partial x}|_{x=a}=\frac{\partial F(a, y)}{\partial y}|_{y=a}$, and $\D_{\tilde g} F(a_1, a_2):=\D_{\tilde g, x}F(x, a_2)|_{x=a_1}=\D_{\tilde g, y}F(a_2, y)|_{y=a_1}.$ Similarly, for a function \;$F(x, y)$\; defined on a open subset of\;$ (\partial M)^2$ which is symmetric and\; with $F(\cdot, \cdot)\in C^1$ with respect to $\tilde g$, we define $\frac{\partial F(a, a)}{\partial a}:=\frac{\partial F(x, a)}{\partial x}|_{x=a}=\frac{\partial F(a, y)}{\partial y}|_{y=a}$.
\vspace{4pt}

\noindent
 For  $l\in \N^*$ and $a\in M$, $O_a(1)$ stands for quantities bounded uniformly in $a$, $O_{l}(1)$ stands for quantities bounded uniformly in \;$l$\; and \;$o_l(1)$ stands for quantities which tends to $0$ as $l\rightarrow +\infty$.  For  $\epsilon$ positive and small, \;$a\in M$\; and \;$\l\in \R_+$ large, $\l\geq \frac{1}{\epsilon}$,\;$O_{a, \l, \epsilon}(1)$\; stands for quantities bounded uniformly in \;$a$, \;$\l$, and $\epsilon$. For $\epsilon$ positive and small, $(p, q)\in \N\times \N$ such that $2p+q=k$, $\bar \l:=(\l_1, \cdots, \l_{p+q})\in (\R_+)^{p+q}$, $\l_i\geq \frac{1}{\epsilon}$  for $i=1, \cdots, p+q$, and $A:=(a_1, \cdots, a_{p+q})\in M^p\times (\partial M)^q$ (where $ (\R_+)^p$ denotes the cartesian product of $p+q$ copies of $\R_+$, and the convention that $M^p\times (\partial M)^0:=M^p$ and $M^0\times (\partial M)^q=(\partial M)^q$ is used), $O_{A, \bar \l, \epsilon}(1)$ stands for quantities bounded uniformly in $A$, $\bar \l$, and $\epsilon$. Similarly for $\epsilon $ positive and small,  $(p, q)\in \N\times \N$ such that $2p+q=k$, $\bar \l:=(\l_1, \cdots, \l_p)\in (\R_+)^{p+q}$, $\l_i\geq \frac{1}{\epsilon}$ for $i=1, \cdots, p+q$, $\bar \alpha:=(\alpha_1, \cdots, \alpha_{p+q})\in \R^{p+q}$, $\alpha_i$ close to $1$ for $i=1, \cdots, p+q$, and $A:=(a_1, \cdots, a_{p+p})\in M^p\times (\partial M)^q$ with still the same convention as above (where $ \R^{p+q}$  denotes the cartesian product of $p+q$ copies of $\R$), $O_{\bar\alpha, A, \bar \l, \epsilon}(1)$ will mean quantities bounded from above and below independent of $\bar \alpha$, $A$, $\bar \l$, and $\epsilon$. For $x\in \R$, we will use the notation $O(x)$ to mean $|x|O(1)$ where $O(1)$ will be specified in all the contexts where it is used. Large positive constants are usually denoted by $C$ and the value of\;$C$\;is allowed to vary from formula to formula and also within the same line. Similarly small positive constants are also denoted by $c$ and their value may varies from formula to formula and also within the same line.
 \vspace{4pt}

 \noindent
 Now, for $(X, A)$ a pair of topological spaces and $q\in \N$, we denote by $H_{q}(X, A)$ its $q$-th relative homology group with $\Z_2$ coefficients. Here $H_q(X)=H_q(X,\emptyset)$. We use the notation $b_q(X)$ and $b_q(X, A)$ to denote the $q$-th betti number of $X$ and $(X, A)$ respectively. We will write $\chi(X)$ and $\chi(X, A)$ for the respective Euler characteristics of $X$ and $(X,A)$.
\vspace{4pt}

\noindent
We call $\bar k$ the number of negative eigenvalues (counted with multiplicity) of $P_g^{4, 3}$. We point out that $\bar k$ can be zero, but is always finite. If $\bar k\geq 1$, then we will denote by \;$E\subset \mathcal{H}_{\frac{\partial}{\partial n}}$ the direct sum of the eigenspaces corresponding to the negative eigenvalues of $P_g^{4, 3}$. The dimension of $E$ is of course $\bar k$. On the other hand, we have the existence of an $L^2$-orthonormal basis of eigenfunctions $v_1,\cdots, v_{\bar k}$ of $E$ satisfying
$$
P_g^{4, 3} v_i=\mu_i v_i\;\;\;\forall\;\;i=1\cdots \bar k,
$$
$$
\mu_1\leq \mu_2\leq \cdots\leq \mu_{\bar k}<0<\m_{\bar k+1}\leq\cdots,
$$
where $\mu_i$'s are the eigenvalues of $P_g^{4,3}$ counted with multiplicity. From the fact that $P_g^{4, 3}$ is self-adjoint and annihilates constants,
we have $E\subset \{u\in\mathcal{H}_{\frac{\partial}{\partial n}}:\;\;\ov{u}_{(Q, T)}=0\}$. We define also the positive definite (on $\{u \in \mathcal{H}_{\frac{\partial}{\partial n}}:\;\;\ov{u}_{(Q, T)}=0\}$) pseudo-differential operator $P^{4, 3,+}_g$ as follows
\begin{equation}\label{eq:operatorrev}
P_g^{4,3,+}u=P_g^{4,3}u-2\sum_{i=1}^{\bar k}\mu_i\left(\int_M u v_idV_g\right)v_i.
\end{equation}
Basically $P^{4, 3,+}_g$ is obtained from $P_g^{4, 3}$ by reversing the sign of the negative eigenvalues and we extend the latter definition to $\bar m=0$ for uniformity in the analysis and recall that in that case $P_g^{4, 3, +}=P_g^{4, 3}$. Using $P_g^{4, 3, +}$, we set for $t>0$
\begin{equation}\label{eq:jt}
\begin{split}
II_t(u):=<P_g^{4, 3, +}u, u>+2t\sum_{r=1}^{\bar m}\mu_r(u^{r})^2+4t\int_M \, Q_gudV_g+4t\oint_{\partial M}T_gudS_g-4\pi^2tk\ln \int_MKe^{4u}dV_g, \\\,\;\;u\in \mathcal{H}_{\frac{\partial}{\partial n}},
\end{split}
\end{equation}
with
\begin{equation}\label{eq:defnegativecr}
u^{r}:=\int_M uv_rdV_g,\;\;\; r=1, \cdots, \bar k, \;\;u\in \mathcal{H}_{\frac{\partial}{\partial n}}.
\end{equation}
Now, using \eqref{eq:operatorrev} and \eqref{eq:jt}, we obtain
\begin{equation}\label{eq:jt1}
\begin{split}
II_t(u):=<P_g^{4, 3}u, u>+2(t-1)\sum_{r=1}^{\bar m}\mu_r(u^{r})^2+4t\int_M Q_g u dV_g+4\oint_{\partial M}uT_gdS_g-4\pi^2tk\ln \int_M Ke^{4u}dV_g,\\ \;u\in  \mathcal{H}_{\frac{\partial}{\partial n}},
\end{split}
\end{equation}
and hence $II=II_1$. Furthermore, using \eqref{eq:defnegativecr}, we define
\begin{equation}\label{eq:defnegativec}
u^{-}=\sum_{r=1}^{\bar m}u^{r}v_r.
\end{equation}\vspace{6pt}

\noindent
We will use the notation $\langle \cdot, \cdot\rangle $ to denote the $L^2$ scalar product. On the other hand, it is easy to see that
\begin{equation}\label{eq:productpr}
\langle u, v\rangle _{P^{4, 3}}:=\langle P^{4, 3+}_gu, v\rangle , \;\,\;u, v\in\{w\in  \mathcal{H}_{\frac{\partial}{\partial n}}:\;\;\;\ov{u}_{(Q, T)}=0\}
\end{equation}
defines a inner product on $\{u\in  \mathcal{H}_{\frac{\partial}{\partial n}}:\;\;\;\ov{u}_{(Q, T)}=0\}$ which induces a norm equivalent to $W^{2, 2}$-norm (on $\{u\in  \mathcal{H}_{\frac{\partial}{\partial n}}:\;\;\;\ov{u}_{(Q, T)}=0\}$) and denoted by
\begin{equation}\label{eq:normpr}
||u||:=\sqrt{<u, u>_{P^{4, 3}}} \;\;\;u\in\{w\in  \mathcal{H}_{\frac{\partial}{\partial n}}:\;\;\;\ov{u}_{(Q, T)}=0\}.
\end{equation}
As above, in the general case, namely $\bar k\geq 0$, for $\epsilon$ small and positive, $\bar \beta:=(\beta_1, \cdots, \beta_{\bar m})\in \R^{\bar m}$  with $\beta_i$ close to $0$, $i=1, \cdots, \bar k$) (where $\R^{\bar k}$ is the empty set when $\bar k=0$), $(p, q)\in \N\times \N$ such that $2p+q=k$,  $\bar \l:=(\l_1, \cdots, \l_{p+q})\in (\R_+)^{p+q}$, $\l_i\geq \frac{1}{\epsilon}$  for $i=1, \cdots, p+q$, $\bar \alpha:=(\alpha_1, \cdots, \alpha_{p+q})\in \R^{p+q}$, $\alpha_i$ close to $1$ for $i=1, \cdots, p+q$, and $A:=(a_1, \cdots, a_{p+q})\in (\mathring{M})^p\times (\partial M)^q$, $w\in \mathcal{H}_{\frac{\partial}{\partial n}}$ with $||w||$ small, $O_{\bar\alpha, A, \bar \l, \bar \beta, \epsilon}(1)$ will stand for quantities bounded independent of $\bar \alpha$, $A$, $\bar \l$, $\bar \beta$, and $\epsilon$, and $O_{\bar\alpha, A, \bar \l, \bar \beta, w, \epsilon}(1)$ will stand for quantities bounded independent of $\bar \alpha$, $A$, $\bar \l$, $\bar \beta$,  $w$ and $\epsilon$.
\vspace{4pt}

\noindent
In the sequel also, \;$II^{c}$\; with $c\in \R$ will stand for \;$II^{c} :=\{u\in \mathcal{H}_{\frac{\partial}{\partial n}}:\;\;II(u)\leq c\}$.  Similarly also, given $c\in \R$, $II_t^c$ stands for $II^{c}_t :=\{u\in  \mathcal{H}_{\frac{\partial }{\partial n}}:\;\;II_t(u)\leq c\}$. We would like to emphasize that in some places in the literature, a different notation is used for the sublevel, precisely with the $c$ as subscript. Here we adopt the notation of P. Rabinowitz.
\vspace{4pt}

\noindent
Given a point $b\in \R^4$ and $\lambda$ a positive real number, we define $\delta_{b, \lambda}$ to be the {\em standard bubble}, namely
\begin{equation}\label{eq:standarbubble}
\delta_{b, \lambda}(y):=\ln \left(\frac{2\lambda}{1+\lambda^2|y-b|^2}\right),\;\;\;\;\;\;y\in \R^4.
\end{equation}
The functions $\delta_{b, \lambda}$ verify the following equation
\begin{equation}\label{eq:bubbleequation}
\D^2\delta_{b,\lambda}=6e^{4\delta_{b,\lambda}}\;\;\;\text{in}\;\;\;\R^4.
\end{equation}
Geometrically, equation $\eqref{eq:bubbleequation}$ means that the metric $g=e^{2\delta_{b, \lambda}} dx^2$ (after pull-back by the stereographic projection) has constant $Q$-curvature equal to $3$ (with $dx^2$ denoting the standard metric on $\R^4$). Furthermore, if $b\in \R^3=\partial \R^4_+$, then $\d_{b, \l}$ satisfies
\begin{equation}\label{eq:bubbleequationb}
\left\{
\begin{split}
\Delta^2 \delta_{b, \l}&=6e^{4\delta_{b, \l}}\;\;\;&\text{in}\;\;\R^4_+,\\
\partial_{x_4}\Delta \delta_{b, \l}&=0\;\;\;&\text{on} \;\;\;\R^3,\\
\partial_{x_4}\delta_{b, \l}&=0\;\;\;&\text{on}\;\;\R^3,
\end{split}
\right.
\end{equation}
with $\bar \R^4_+=\R^3\times \bar\R_+$ and a point $x\in \bar\R^4_+$ has the following representation $x=(x_1, \cdots, x_4)$. As above, equation $\eqref{eq:bubbleequationb}$ has also a geometric interpretation. Indeed, it is equivalent to the fact that the metric $g=e^{2\delta_{b, \lambda}} dx^2$ (after pull-back by the stereographic projection) has constant $Q$-curvature equal to $3$, zero $T$-curvature and vanishing mean curvature (with $dx^2$ denoting the standard metric on $\bar \R^4_+$).
Using the existence of conformal normal coordinates (see \cite{lp}, \cite{gun}, and \cite{marques}) and recalling that $H_g=0$, then for every $m$ large positive integer, we have that for $a \in M$, there exists a function $u_a\in C^{\infty}(M)$ such that the metric $g_a = e^{2u_a} g$ verifies
\begin{equation}\label{eq:detga}
det g_a(x)=1 +O_{a, x}((d_{g_a}(x, a))^m)\;\;\text{for}\;\;\; x\in B^{g_a}_a( \varrho_a).
\end{equation}
with $O_{a, x}(1)$ meaning bounded by a constant independent of $a$ and $x$, $0<\varrho_a< \max(\frac{inj_{g_a}(M)}{10}, \frac{inj_{\hat{g}_a}(\partial M)}{10})$. Moreover, we can take the family of functions $u_a$, $g_a$ and $\varrho_a$ such that
\begin{equation}\label{eq:varro0}
\text{the maps}\;\;\;a\longrightarrow u_a, \;g_a\;\;\text{are}\;\;C^1\;\;\;\text{and}\;\;\;\;\varrho_a\geq \varrho_0>0,
\end{equation}
for some small positive $\varrho_0$ satisfying $\varrho_0< \max(\frac{inj_g(M)}{10}, \frac{inj_{\hat{g}_a}(\partial M)}{10})$,and
\begin{equation}\label{eq:proua}
\begin{split}
&||u_a||_{C^4(M)}=O_a(1),\;\;\frac{1}{\ov C^2} g\leq g_a\leq \ov C^2 g, \;\;\;a\in M, \\\;&u_a(x)= O_a(d^2_{g_a}(a, x))=O_a(d_{g}^2(a, x)) \;\;\text{for}\;\;x\in\;\;B_a^{g_a}(\varrho_0)\supset B_a(\frac{\varrho_0}{2\ov C}),\;\;\;\;a\in M,\\&
u_a(a)=0,\;a\in M\;\;\;R_{g_a}(a)=0,\;\;a\in  M\;\;\text{and}\;\; H_{g_a}(a)=0\;\;a\in \partial M,
\end{split}
\end{equation}
for some large positive constant $\ov C$ independent of $a$. For the meaning of $O_a(1)$ in \eqref{eq:proua}, see section \ref{eq:notpre}. Furthermore, for $a\in M$, using the scalar curvature equation satisfied by $e^{-u_a}$, namely $-\D_{g_a} (e^{-u_a})+\frac{1}{6}R_{g_a}e^{-u_a}=\frac{1}{6}R_g(a)e^{-3u_a}$\; in $M$, and \eqref{eq:detga}-\eqref{eq:proua}, it is easy to see that the following holds
\begin{equation}\label{eq:lapemua}
\D_{g_a}(e^{-u_{a}})(a)=-\frac{1}{6}R_g(a).
\end{equation}
Similarly, for $a\in \partial M$, using the mean curvature equation satisfied by $e^{-u_a}$, namely $-\frac{\partial }{\partial n}_{g_a} (e^{-u_a})+H_{g_a}e^{-u_a}=H_g(a)e^{-2u_a}$\; on $\partial M$, and \eqref{eq:detga}-\eqref{eq:proua}, or just \eqref{eq:proua}, it is easy to see that the following holds
\begin{equation}\label{eq:lapemuaneu}
\frac{\partial }{\partial n}_{g_a}(e^{-u_{a}})(a)=0.
\end{equation}
For $a\in M$, and $r>0$, we set
\begin{equation}\label{eq:expballaga}
exp_a^{a}:=exp_{a}^{g_a}\;\;\;\;\text{and}\;\;\;B_a^{a}(r):=B_a^{g_a}(r).
\end{equation}
On the other hand, using the properties of $g_a$ (see \eqref{eq:detga}-\eqref{eq:proua})), it is easy to check that for every $u\in C^2(B_a^{a}(\varrho))$ with $0<\varrho<\frac{\varrho_0}{4}$ there holds
\begin{equation}\label{eq:gradlapga}
\begin{split}
&\n_{g_a}u(a)=\n_gu(a)=\n_4\hat u(0),\;\;\;\;\;\D_{g_a}u(a)=\D_4\hat u(0), \;\;\text {if}\;\; d_{g_{a}}(a, \partial M)\geq 4\varrho,\\
\;\;&\n_{\hat g_a}u(a)=\n_{\hat g}u(a)=\n_3\hat u(0), \;\;\frac{\partial }{\partial n}_{g_a}u(a)=\frac{\partial }{\partial n}_{g}u(a)=\partial_{x_4}\hat u(0), \;\;\text{if}\;\;a\in \partial M,
\end{split}
\end{equation}
where
\begin{equation}\label{eq:hatu}
\hat u(y)=u(exp_a^{a}(y)),\;\;y\in B_0^4(\varrho)\;\;\text {if}\;\; d_{g_{a}}(a, \partial M)\geq 4\varrho\;\;\text{and}\;\;y\in B_0^{\R^4_+}(\varrho)\;\text{if}\;\;a\in \partial M,
\end{equation}
with
$$
\hat g_a:=g_a|_{\partial M}\;\;\text {and}\;\;\hat g:=g|_{\partial M}.
$$
\vspace{4pt}

\noindent
Now for $0<\varrho<\frac{\varrho_0}{4}$, we define the  cut-off function $\chi_{\varrho} : \R_+ \rightarrow \R_+$ satisfying the following properties:
\begin{equation}\label{eq:cutoff}
\begin{cases}

\chi_{\varrho}(t)  = t \;\;&\text{ for } \;\;t \in [0,\varrho],\\

\chi_{\varrho}(t) = 2 \varrho \;\;&\text{ for } \;\; t \geq 2 \varrho, \\

\chi_{\varrho}(t) \in [\varrho, 2 \varrho] \;\;\;&\text{ for } \;\; t\in [\varrho, 2 \varrho].

\end{cases}
\end{equation}
Using the cut-off function $\chi_{\varrho}$, we define for $a\in M$ and $\lambda\in \R_+$  the function $\hat{\delta}_{a, \lambda}$ as follows
\begin{equation}\label{eq:hatdelta}
\hat{\delta}_{a, \lambda}(x):=\ln \left(\frac{2\lambda}{1+\lambda^2\chi_{\varrho}^2(d_{g_a}(x, a))}\right).
\end{equation}
Next, using the pull-back standard bubbles, namely \eqref{eq:hatdelta}, for $a\in  M$, $d_g(a, \partial M)\geq 4\ov C\varrho$, and $\lambda\in \R_+$, we define $\varphi_{a, \lambda}$ to be the solution of the following projected boundary value problem
\begin{equation}\label{eq:projbubble}
\left\{
\begin{split}
P^4_g \varphi_{a,\l} \, + \, \frac{4}{k} Q_g &=16\pi^2\, \frac{e^{4(\hat{\d}_{a,\l} \,  +  \, u_a)}}{\int_Me^{4 (\hat{\d}_{a,\l} \,  +  \, u_a)}dV_g}\;\;&\text{in}\;\;\mathring{M},\\
P^3_g \varphi_{a,\l} \, + \, \frac{2}{k} T_g &=0\;\;&\text{on}\;\;\partial M,\\
\frac{\partial \varphi_{a, \l}}{\partial n_g}&=0\;\;&\text{on}\;\;\partial M,\\
\ov{\varphi_{a,\l}}_{(Q, T)}=0 &=0.
\end{split}
\right.
\end{equation}
Similarly, for $a\in \partial M$, and $\lambda\in \R_+$, we define $\varphi_{a, \lambda}$ to be the solution of the following projected boundary value problem
\begin{equation}\label{eq:projbubble1}
\left\{
\begin{split}
P^4_g \varphi_{a,\l} \, + \, \frac{2}{k} Q_g &=8\pi^2\, \frac{e^{4(\hat{\d}_{a,\l} \,  +  \, u_a)}}{\int_Me^{4 (\hat{\d}_{a,\l} \,  +  \, u_a)}dV_g}\;\;&\text{in}\;\;\mathring{M},\\
P^3_g \varphi_{a,\l} \, + \, \frac{1}{k} T_g &=0\;\;&\text{on}\;\;\partial M,\\
\frac{\partial \varphi_{a, \l}}{\partial n_g}&=0\;\;&\text{on}\;\;\partial M,\\
\ov{\varphi_{a,\l}}_{(Q, T)}=0 &=0
\end{split}
\right.
\end{equation}
So differentiating with respect to $\l$ and $a$ (respectively) the relation $\ov{\varphi_{a, \l}}_{(Q, T)}=0$, we get (respectively)
\begin{equation}\label{eq:normlambda}
\ov{\frac{\partial \varphi_{a, \l}(x)}{\partial \l}}_{(Q, T)}=0,
\end{equation}
and
\begin{equation}\label{eq:norma}
\ov{\frac{\partial \varphi_{a, \l}(x)}{\partial a}}_{(Q, T)}=0.
\end{equation}
Now, we recall that $G$ is the unique solution of the following BVP
\begin{equation}\label{eq:defG4}
\left\{
\begin{split}
P^4_g G(a, \cdot) \, + \, \frac{4}{k} Q_g &=16\pi^2\, \d_a(\cdot)\;\;&\text{in}\;\;\mathring{M},\\
P^3_g G(a, \cdot) \, + \, \frac{2}{k} T_g &=0\;\;&\text{on}\;\;\partial M,\\
\frac{\partial G(a, \cdot)}{\partial n_g}&=0\;\;&\text{on}\;\;\partial M,\\
\ov{G(a, \cdot)}_{(Q, T)}=0 &=0.
\end{split}
\right.
\end{equation}
Using  \eqref{eq:defG4}, it is easy to see that the following integral representation formula holds
\begin{equation}\label{eq:G4integral}
u(x)-\ov{u}_{(Q, T)}=\frac{1}{16\pi^2}\left(\int_M G(x, y)P_gu(y)dV_g(y)+2\oint_{\partial M}G(x, y)P^3_gu(y)dS_g(y)\right), \;\;\;u\in C^4(M), \;x\in M,
\end{equation}
where $\ov{u}_{(Q, T)}$ is defined as in Section \ref{eq:notpre}. It is a well know fact that $G$ has a logarithmic singularity. In fact $G$ decomposes as follows
\begin{equation}\label{eq:decompG4}
 G(a,x)=S(a, x)+H(a, x),
\end{equation}
where
\begin{equation*}
S(x,y)=
  \left\{
    \begin{array}{ll}
      \ln  \left(\frac{1}{\chi_{\varrho}^2(d_{{g}_a}(a, x))}\right)\;\;\;\;\;\text{if}\;\;\;\;B^a_{x}(\varrho)\cap\partial M=\emptyset,&\\\\
      \ln  \left(\frac{1}{\chi_{\varrho}^2(d_{{g}_a}(a, x))}\right)+\ln  \left(\frac{1}{\chi_{\varrho}^2(d_{{g}_a}(a, \bar x))}\right)\;\;\;\;\text{otherwise},&
  \end{array}
  \right .
\end{equation*}
with $\bar x$ denoting the normal reflection of $x$ through $\partial M$ with respect to $g_a$ and
\begin{equation}\label{eq:regH4}
G\in C^{\infty}(M^2\setminus Diag(M)), \,\;\text{and } \;H\;\;\text{extend to a}\;\; C^{3, \beta}(M^2)\;\;\text{function},\;\forall \beta\in (0, 1).
\end{equation}
Now, using \eqref{eq:limitfsint}, \eqref{eq:limitfsbound}, and \eqref{eq:partiallimitint}, \eqref{eq:partiallimitbound} combined with the symmetry of $H$, it is easy to see that for every $(p, q)\in \N^2$ such that $2p+q=k$, there holds
\begin{equation}\label{eq:relationderivativeint}
\frac{\partial \mathcal{F}_{p, q}(a_1, \cdots, a_{p+q})}{\partial a_i}=\frac{\n_g\mathcal{F}^{A}_i(a_i)}{\mathcal{F}^{A}_i(a_i)}, \;\;\;i=1, \cdots, p.
\end{equation}
and
\begin{equation}\label{eq:relationderivativebound}
\frac{\partial \mathcal{F}_{p, q}(a_1, \cdots, a_{p+q})}{\partial a_i}=\frac{\n_{\hat g}\mathcal{F}^{A}_i(a_i)}{2\mathcal{F}^{A}_i(a_i)}, \;\;\;i=p+1, \cdots, p+q.
\end{equation}
Next, for $(p, q)\in \N^2$ such that $2p+q=k$ and $A\in (\mathring{M}^p)^*\times ((\partial M)^q)^*$, we set\begin{equation}\label{eq:deflA}
l_K(A):=
\begin{cases}
\sum_{i=1}^p\left(\frac{\D_{g_{a_i}} \mathcal{F}^{A}_i(a_i)}{\sqrt{\mathcal{F}^{A}_i(a_i)}}-\frac{2}{3}R_g(a_i)\sqrt{\mathcal{F}^{A}_i(a_i)}\right),\;\;&\text{if}\,\;q=0,\\\\
\sum_{i=p+1}^{p+q}\frac{1}{4(\mathcal{F}_i^A)^{\frac{3}{4}}(a_i)}\frac{\partial \mathcal{F}_i^A}{\partial n_{g_{a_i}}}(a_i),&\text{if}\;\;q\neq 0,
\end{cases}
\end{equation}
and use the properties of the metrics $g_{a_i}$ $i=1, \cdots, p$, the transformation rule of the conformal Laplacian under conformal change of metrics and direct calculations, to get
\begin{equation}\label{eq:auxiindexa1}
l_K(A)=4\mathcal{L}_K(A), \;\;\forall A\in Crit(\mathcal{F}_{p, q}).
\end{equation}
\vspace{4pt}

\noindent
For  $k\geq 2$ and $\bar k\geq 1$, we define $A^{\partial}_{k-1,\ov{k}}$ to be the colimit of the diagram

\begin{equation*}
\xymatrix{ \partial B_{1}^{\ov{k}}& B^{\partial}_{k-1}(M) \times \partial B_1^{\ov{k}}\ar[r]\ar[l]& B_{k-1}^{\partial}(M) \times B_1^{\ov{k}}}
\end{equation*}




\noindent
We denote by $B_l(\partial M)$ the space of formal barycenters of $\partial M$ of order $l$;
 \begin{equation}\label{eq:barytop}
B_{l}(\partial M):= \big \{\sum_{i=1}^{l}\alpha_i\d_{a_i}, a_i\in \partial M, \alpha_i\geq 0, i=1,\cdots, j,\;\,\sum_{i=1}^{l}\alpha_i \, = \, k \big \}.
\end{equation}
and set
\begin{equation*}
\xymatrix{ A_{i, \bar k}:= \partial B_1^{\ov{k}} &  B_i(\partial M) \times \partial B_1^{\ov{k}}\ar[r]\ar[l]&B_i(\partial M) \times B^{\ov{k}}_1}
\end{equation*}




\vspace{6pt}

\noindent
 Now we  state a  result describing the local behavior of blowing up sequences of solutions to some perturbations of (BVP) \eqref{eq:bvps}. In \cite{nd2}, it is proven that, for blowing up sequence of solutions to the following type of perturbations of (BVP) \eqref{eq:bvps}
\begin{equation}\label{eq:sepexact}
\left\{
\begin{split}
P^4_gu_l+2t_lQ_g&=2t_lKe^{4u_l}\;\;&\text{in}\;\;\mathring{M},\\
P^3_gu_l+t_lT_g&=0\;\;&\text{on}\;\;\partial M,\\
\frac{\partial u_l}{\partial n_g}&=0\;\;&\text{on}\;\;\partial M,
\end{split}
\right.
\end{equation}
where
\begin{equation}\label{eq:lbqexact}
t_l\longrightarrow 1\;\;\;\;\text{as}\;\;l\longrightarrow + \infty,
\end{equation}
the following lemma holds
\begin{lem}\label{eq:locdesblow}
Assuming that $(u_{l})$ is a sequence of solutions of \eqref{eq:sepexact} with $t_{l}$ satisfying \eqref{eq:lbqexact}, then there exists \;$(p, q)\in \N^2$ with $p+q\in \N^*$ such that up to a subsequence, there exists\;$p+q$\;converging sequences of points\; $(x_{i,l})_{l\in \N}$\;\;with limits $x_i$, $i=1,\cdots,p+q $,\;$x_i\in \mathring{M}$ for $i=1,\cdots, p$, $x_i\in \partial M$ for $i=p+1,\cdots, p+q$,\;\;$p+q$\; sequences $(\mu_{i,l})_{l\in \N} \;\;i=1,\cdots,p+q$\;of positive real numbers  converging to \;$0$\;such that the following hold: \\\\
a)
\\
\begin{equation}
\begin{split}
\hspace{-15pt}\frac{d_g(x_{i,l},x_{j,l})}{\mu_{i,l}}\longrightarrow +\infty \;\;\;\text{as}&\;\;l\rightarrow +\infty\;\;\;\; i\neq j \;\;i,j =1.\cdots,p+q,\;\;\;\;\\& \text{and}\\ \;\;\;\;t_l K(x_{i})\mu_{i, l}^{4}e^{4u_{l}(x_{i})}&e^{-4\ln 2}=3,\;\;i=1, \cdots, p+q.
\end{split}
\end{equation}
b)
\\

\begin{equation}
 \begin{split}
 &\text{For}\; \;i=1\cdots, p,\\
 &v_{i,l}(x)=u_{l}(exp_{x_{i}}(\mu_{i,l}x))-u_{l}(x_{i})+\ln 2\longrightarrow V_{0}(x):=\ln\left(\frac{2}{1+|x|^{2}}\right)\;\;\;\; in\; \;\;C^{4}_{loc}(\R^4)\;\;\;\text{as}\;\;l\rightarrow +\infty,\\&\text{and for}\;\; i=p+1, \cdots, p+q,\\  &v_{i,l}(x)=u_{l}(exp_{x_{i}}(\mu_{i,l}x))-u_{l}(x_{i})+\ln 2\longrightarrow V_{0}(x):=\ln\left(\frac{2}{1+|x|^{2}}\right)\;\;\;\; in\; \;\;C^{4}_{loc}(\bar\R^4_+)\;\;\;\text{as}\;\;l\rightarrow +\infty.
\end{split}
\end{equation}
c)
\\
\begin{equation}
\text{There exists}\;\;C>0\;\; \text{such that}\;\;
\inf_{i=1,...,p+q}d_{g}(x_{i,l},x)^{4}e^{4u_{l}(x)}\leq C \;\;\;\;\forall x\in M,\;\;\forall l\in \N.
\end{equation}
d)\\
\begin{equation}
 \begin{split}
t_l Ke^{4u_l}dV_g\rightarrow 8\pi^2\sum_{i=1}^p\d_{x_{i}}&+4\pi^2\sum_{i=p+1}^{p+q}\d_{x_{i}}  \;\;\;\text{in the sense of measure}\;\;\;\text{as}\;\;l\rightarrow +\infty,\\&
\text{and}
\\
\lim_{l\rightarrow +\infty}\int_{M}t_l Ke^{4u_l}dV_g&=4(2 p+q)\pi^2.
\end{split}
\end{equation}
e)\\
\begin{equation}
\begin{split}
u_l-\ov{u_l}_{(Q, T)}\rightarrow \sum_{i=1}^p G(x_{i}, \cdot)+\frac{1}{2}\sum_{i=p+1}^{p+q}G(x_{i}, \cdot)\;\;\;\text{in}\;\;\;C^4_{loc}(M\setminus\{x_1,\cdots, x_{p+q}\}),\\\;\ov{u_l}_{(Q, T)}\rightarrow -\infty\;\;\;\text{as}\;\;l\rightarrow +\infty.
\end{split}
\end{equation}
\end{lem}

\section{Blow up analysis and  deformation lemma}\label{s:bldeformlem}
In this section, we derive a Bahri-Lucia type deformation lemma which is a refined version of the classical Bahri-Lucia deformation Lemma. Indeed from the works of Bahri\cite{bah1}, \cite{bah2} and Lucia\cite{lu}, we have the following deformation Lemma which will improve using  refined blow-up analysis.
\begin{lem}\label{eq:deformlem}
 Assuming that $a, b\in \R$ such that $a<b$ and  there is no critical values of $II$ in $[a, b]$, then there are two possibilities\\
1) Either  $$II^a\;\text{is a deformation retract of}\;\; II^b.$$
2) Or there exists a sequence $t_l\rightarrow 1$ as $l\rightarrow +\infty$ and a sequence of critical point $u_l$ of $II_{t_l}$ verifying $a\leq II(u_l)\leq b$ for all $l\in \N^*$.
\end{lem}
Actually the proof of above Lemma provides a pseudogradient whose  the noncompact $\omega-$limit  set (i.e. the endpoints of noncompact orbits) is in one to one correspondence to the blow up set of some  approximation of type  \eqref{eq:sepexact}. Hence for a better understanding of such noncompact orbits, one has to describe  the behavior of the blowing up solutions near their blow up point. To that aim we prove the following Lemma:

\begin{pro}\label{eq:sharpest}
Assuming that $(u_l)_{l\in \N}$ is a bubbling sequence of solutions to \eqref{eq:sepexact} with $t_l$ satisfying \eqref{eq:lbqexact}, then up to a subsequence and keeping the notation in Lemma \ref{eq:locdesblow}, we have that the points $x_{i, l}$ are uniformly isolated, and the scaling parameters $\l_{i, l}:=\mu_{i, l}^{-1}$ are comparable, namely there exists $0<\eta_k<\frac{\varrho_0}{10}$ and $0<\varrho_k<\frac{\varrho_0}{10}$ two positive and small real numbers, where $\varrho_0$ is as in \eqref{eq:varro0}, and a large positive constant \;$\Lambda_k$ such that for $l$ large enough, there holds
\begin{equation}\label{eq:uniformiso}
\begin{split}
d_g(x_{i, l}, x_{j, l})\geq 4\ov C\eta_k,\;\;\;\;\forall\;i\neq j=1, \cdots, p+q,\;\;\;\;,\\\text{and}\;\;\;\;  \Lambda_k^{-1}\l_{j,l}\leq\l_{i,l} \leq \Lambda_k\l_{j,l},\;\;\forall\;i,j=1, \cdots, p+q.
\end{split}
\end{equation}
Furthermore, the interior concentrations (if any) are uniformly far from $\partial M$, namely if $p>0$, then for $l$ large enough, there holds
\begin{equation}\label{eq:uniformfarboun}
d_g(x_{i, l}, \partial M)\geq 4\ov C\varrho_k,\;\;\;\;\forall\;i=1, \cdots, p.
\end{equation}
Moreover, we have that the following estimate around the blow-up points holds (for $l$ large enough)
\begin{equation}\label{eq:sharpestimate}
u_l(x)+\frac{1}{4}\ln\frac{t_l K(x_{i})}{3}=\ln \frac{2\l_{i, l}}{1+\l_{i, l}^2(d_{g_{x_i}}(x, x_{i}))^2}+ O((d_{g_{x_i}}(x, x_{i}))),\;\;\;\;\forall\;x\in\;B^{x_i}_{x_{i}}(\eta_k), \;\;i=1, \cdots, p+q,
\end{equation}
where $g_{x_i}$, $\varrho$, and $\ov C$ are as in Section \ref{eq:notpre}, see \eqref{eq:detga} -\eqref{eq:proua}.
\end{pro}
Now, before proving Proposition \ref{eq:sharpest}, we would like first to show how it provides the refined Bahri-Lucia type deformation Lemma that we mentioned above. To do so, we start by defining the neighborhood of potential critical points at infinity of $II$ and for that we first fix $\L>\L_k$ to be a large positive constant. Next, for $(p, q)\in \N^2$ with $2p+q=k$, $\epsilon$ , $\varrho$ and $\eta$ small positive real numbers with $0<\varrho<\varrho_k$, and $0<\eta<\eta_k$, we define the $(p, q, \epsilon, \varrho, \eta)$-neighborhood of potential critical points at infinity of $II$ as follows
\begin{equation}\label{eq:ninfinity}
\begin{split}
V(p, q, \epsilon, \varrho, \eta):=\{u\in \mathcal{H}_{\frac{\partial }{\partial n}}\;:a_1, \cdots, a_p\in \mathring{M}, \;\;a_{p+1}, \cdots, a_{ p+q}\in \partial M,\;\;\;\;\l_1,\cdots, \l_{p+q}>0, \\||u-\ov{u}_{(Q, T)}-\sum_{i=1}^{p+q}\varphi_{a_i,\lambda_i}||+||\n^{P^{4, 3}} II(u-\ov{u}_{(Q, T)})||=O\left( \sum_{i=1}^{p+q}\frac{1}{\l_i}\right)\\\;\;\;\lambda_i\geq \frac{1}{\epsilon}, \;i=1, \cdots, p+q,\;\;\;\;\frac{2}{\Lambda}\leq \frac{\l_i}{\l_j}\leq \frac{\Lambda}{2}, \;i, j=1, \cdots, p+q\\\;\;d_g(a_i, a_j)\geq 4\ov C\eta\;\;\text{for}\;i\neq j, i, j=1, \cdots, p+q,\;\;\text{and}\;\;\;d_g(a_i, \partial M)\geq 4\ov C\varrho, \;\i=p+1, \cdots, p+q\},
\end{split}
\end{equation}
where $\ov C$ is as in \eqref{eq:proua}, $\L_k$, $\eta_k$, and $\varrho_k$ are given by Proposition \ref{eq:sharpest},  $\n^{P^{4, 3}}J$ is the gradient of $II$ with respect to $<\cdot, \cdot>_{P^{4, 3}}$, and $O(1):=O_{A, \bar \l, u, \epsilon}(1)$ meaning bounded uniformly in  $\bar\l:=(\l_1, \cdots, \l_{p+q})$, $A:=(a_1, \cdots, a_{p+q})$, $u$, $\epsilon$. Next, using the above set, Proposition \ref{eq:sharpest} and the method of the proof of Proposition 3.3 in \cite{no1}, we have that Lemma \ref{eq:locdesblow} and Proposition  \ref{eq:sharpest} imply the following Lemma.
\begin{lem}\label{eq:escape}
Let $\epsilon$, $\varrho$, $\eta$ be small positive real numbers with $0<\varrho<\varrho_k$, $0<\eta<\eta_k$, where $\varrho_k$ and $\eta_k$ are given by Proposition \ref{eq:sharpest}. Assuming that $u_l$ is a sequence of blowing up critical point of $II_{t_l}$ with $\ov{(u_l)}_{(Q, T)}=0, l\in \N$ and $t_l\rightarrow 1$ as $l\rightarrow +\infty$, then there exists \;$(p, q)\in \N^2$\; with \;$2p+q=k$\; and \;$l_{\epsilon,\varrho,  \eta}$\; a large positive integer such that for every \;$l\geq l_{\epsilon, \varrho, \eta}$, we have \;$u_l\in V(p,q, \varrho, \epsilon, \eta)$.
\end{lem}
On the other hand, as in \cite{no1} and by the same arguments, we have that Lemma \ref{eq:deformlem} and Lemma \ref{eq:escape} implies the following refined version of the classical Bahri-Lucia deformation lemma.
\begin{lem}\label{eq:deformlemr}
Assuming that $\epsilon$, $\varrho$, and $\eta$ are small positive real numbers with , $0<\varrho<\varrho_k$, $0<\eta<\eta_k$, where $\varrho_k$ and $\eta_k$ are given by Proposition \ref{eq:sharpest}, then for $a, b\in \R$ such that $a<b$, we have that if  there is no critical values of $II$ in $[a, b]$, then there are two possibilities:\\
1) Either  $$II^a\;\text{is a deformation retract of}\;\; II^b.$$
2) Or there exists a sequence $t_l\rightarrow 1$ as $l\rightarrow +\infty$ and a sequence of critical point $u_l$ of $II_{t_l}$ verifying $a\leq II(u_l)\leq b$ for all $l\in \N^*$,  \;$(p, q)\in \N^2$\; with \;$2p+q=k$,\; and \;$l_{\epsilon,\varrho,  \eta}$\; a large positive integer such that  $u_l\in V(p, q, \epsilon, \varrho, \eta)$ for all $l\geq l_{\epsilon, \varrho, \eta}$.
\end{lem}

\vspace{4pt}

\noindent
Next we come back to the proof of Proposition \ref{eq:sharpest} and for that we are going to divide the remainder of this section into three subsections. In the first one, we show that the blow-up points are uniformly isolated and the interior concentration  points (if any) are uniformly far from $\partial M$. In the second one, we prove that the Proposition \ref{eq:sharpest} holds with the classical form of the sup+inf estimate \eqref{eq:sharpestimate}, precisely (for non experts) with $O((d_{g_{x_i}}(x, x_{i})))$ replaced just by $O(1)$. The last subsection deals with the full version of sup+inf estimate \eqref{eq:sharpestimate}.

\subsection{ Blow up points are isolated and  interior ones are far from the boundary}

As already mentioned above, in this subsection, we show that the blow up points are uniformly isolated and that the interior blow-up points (if any) are uniformly far from $\partial M$. Precisely, we prove the following Lemma:
\begin{lem}\label{eq:isolated}
Assuming that $(u_l)_{l\in \N}$ is a bubbling sequence of solutions to BVP \eqref{eq:sepexact} with $t_l$ satisfying \eqref{eq:lbqexact}, then keeping the notations in Lemma \ref{eq:locdesblow}, we have that the points $x_{i, l}$ are uniformly isolated, namely there exists $0<\eta_k<\frac{\varrho_0}{10}$ (where $\varrho_0$ is as in \eqref{eq:varro0}) such that for $l$ large enough, there holds
\begin{equation}\label{eq:unifiso}
d_g(x_{i, l}, x_{j, l})\geq 4\ov C\eta_k ,\;\;\;\forall i\neq j=1, \cdots, p+q.
\end{equation}
Furthermore, there exists $0<\varrho_k<\frac{\varrho_0}{10}$ (where $\varrho_0$ is as in \eqref{eq:varro0}) such that if $p>0$, then for $l$ large enough, there holds
\begin{equation}\label{eq:uniffarboun}
d_g(x_{i, l}, \partial M)\geq 4\ov C \varrho_k, \;\;\;\forall i=1, \cdots, p.
\end{equation}

\end{lem}
\begin{pf}
We are going to use the method of \cite{nd7} to prove Lemma and hence we will be sketchy in many arguments. As in \cite{nd7}, we first fix $1<\nu<2$, and for $i=1, \cdots, p+q$, we set
$$
\bar u_{i, l}(r)=Vol_g(\partial B_{x_{i}}(r))^{-1}\int_{\partial B_{x_{i}}(r)}u_l(x)d\sigma_g(x),\;\;\;\forall\;0\leq r<inj_g(M),
$$
and
$$
\psi_{i,l}(r)=r^{4\nu}exp(4\bar u_{i,l}(r)),\;\;\;\forall\;0\leq r<inj_g(M).
$$
Furthermore, as in \cite{nd7}, we define $r_{i, l}$ as follows
\begin{equation}\label{eq:ril}
r_{i, l}:=\sup\{R_{\nu}\mu_{i, l}\leq r\leq \frac{R_{i, l}}{2}\;\;\text{such that}\;\;\psi_{i,l}^{'}(r)<0\;\;\text{in}\;[R_{\nu}\mu_{i, l}, r[\};
\end{equation}
where $R_{i, l}:=\min_{j\neq i}d_g(x_{i, l}, x_{j, l})$.
Thus, by continuity and the definition of $r_{i, l}$, we have that
\begin{equation}\label{eq:varr}
\psi_{i, l}^{'}(r_{i, l})=0.
\end{equation}
Now, as in \cite{nd7}, to prove \eqref{eq:unifiso}, it suffices to show that $r_{i, l}$ is bounded below by a positive constant in dependent of $l$. Thus, we assume by contradiction that (up to a subsequence) $r_{i, l}\rightarrow 0$ as $l\rightarrow +\infty$ and look for a contradiction. In order to do that, we use the integral representation formula for $(P^4_g, P^3_g)$ under homogeneous Neumann boundary condition and the integral method of \cite{nd4}, to derive the following estimate
$$
\psi^{'}_{i,l}(r_{i,l})\leq (r_{i,l})^{4\nu-1}exp(\bar u_{i,l}(r_{i,l}))\left (4\nu-8C+o_{l}(1)+O_l(r_{i,l})\right),
$$
with $C>1$. So from\;$1<\nu<2$, $C>1$ and $r_{i, l}\longrightarrow 0$ as $l\rightarrow +\infty$, we  deduce that for \;$l$\; large enough, there holds
\begin{equation}\label{eq:negativeb}
\psi_{i,l}^{'}(r_{i,l})<0.
\end{equation}
Thus, \eqref{eq:varr} and \eqref{eq:negativeb} lead to a  contradiction, thereby concluding the proof of \eqref{eq:unifiso}. Hence, the proof of the lemma is complete, since clary \eqref{eq:unifiso} implies \eqref{eq:uniffarboun}.
\end{pf}
\subsection{Harnack-type inequality  around blow-up points}
In this subsection, we present the weak form of Proposition \ref{eq:sharpest} that we mentioned above, namely we show that the difference of a bubbling sequence of solutions to BVP \eqref{eq:sepexact} with $t_l$ satisfying \eqref{eq:lbqexact} and the {\em pull back} bubble around a blow-up point is a $O(1)$. Indeed, we will prove the following Lemma:
\begin{lem}\label{eq:roughestimate}
Assuming that $(u_l)_{l\in \N}$ is a bubbling sequence of solutions to BVP \eqref{eq:sepexact} with $t_l$ satisfying \eqref{eq:lbqexact}, then keeping the notations in Lemma \ref{eq:locdesblow} and Lemma \ref{eq:isolated}, we have that for $l$ large enough, there holds
\begin{equation}\label{eq:classupinf}
u_l(x)+\frac{1}{4}\ln \frac{t_l K(x_{i})}{3}=\ln \frac{2\l_{i, l}}{1+\l_{i, l}^2(d_{g_{x_i}}(x, x_{i}))^2}+ O(1),\;\;\;\;\forall\;x\in\;B^{x_i}_{x_{i}}(\eta_k),
\end{equation}
up to choosing $\eta_k$ smaller than in Lemma \ref{eq:isolated}.
\end{lem}
\begin{rem}\label{eq:comparability}
We point out that the comparability of the scaling parameters $\l_{i,l}$'s follows directly from Lemma \ref{eq:roughestimate}.
\end{rem}
\begin{pf}
We are going to use the method of \cite{nd7}, hence we will be sketchy in many arguments. As in \cite{nd7}, thanks to Lemma \ref{eq:isolated}, we  will focus only on one blow-up point and called it $x$. We point out that $x$ may lie in $M$  with $d_g(x, \partial M)\geq 4\ov C\varrho_k$ or $x\in \partial M$. Thus, we are in the situation where there exists  a sequence $x_l\in M$ such that $x_l\rightarrow x$ with $x_l$ local maximum point for $u_l$ and $u_l(x_l)\rightarrow +\infty$.
Now, we set  $\hat g=e^{2u_x}g$ and choose $\eta_1$ such that
$20\eta_1<\min\{inj_g(M), inj_{\hat{g}}(\partial M), \varrho_0, \varrho_{k},  d\}$ with $4d\leq r_{i,l}$ where $r_{i,l}$ is as in the proof of Lemma \ref{eq:isolated}. Next, we let $\hat w$ be the unique solution of the following boundary value problem
\begin{equation}\label{eq:bvpauxii}
\left\{
\begin{split}
P^4_{\hat g}\hat w&=P^4_{\hat g}u_x\;\;&\text{in}\;\; \mathring{M},\\
P^3_{\hat g}\hat w&=P^3_{\hat g}u_x\;\;&\text{on}\;\;\partial M,\\
\frac{\partial \hat w}{\partial n_{\hat g}}&=0\;\;&\text{on}\;\;\partial M,\\
\ov{\hat w}_{(Q, T)}&=\ov{u}_x.
\end{split}
\right.
\end{equation}
Using standard elliptic regularity theory and \eqref{eq:proua}, we derive
\begin{equation}\label{eq:noinflu}
\hat w(y)=O(d_g(y, x)) \;\;\;\text{in}\;\;\;\; B_{x}^{\hat g}(2\eta_1).
\end{equation}
On the other hand, using the conformal covariance properties of the Paneitz operator and of the Chang-Qing one, see \eqref{eq:law}, we have that $\hat u_l:=u_l-\hat w$ satisfies
\begin{equation*}
\left\{
\begin{split}
P^4_{\hat g}\hat u_l+2\hat Q_l&=2t_lKe^{4\hat u_l}\;\;&\text{in}\;\;\mathring{M},\\
P^3_{\hat g}\hat u_l+\hat T_l&=0\;\;&\text{on}\;\;\partial M,\\
\frac{\partial \hat u_l}{\partial n_{\hat g}}&=0\;\;&\text{on}\;\;\partial M.
\end{split}
\right.
\end{equation*}
with
$$
\hat Q_l=t_le^{-4\hat w}Q_g+\frac{1}{2}P^4_{\hat g}\hat w\;\;\;\text{and}\;\;\;\hat T_l=t_le^{-3\hat w}T_g+P^3_{\hat g}\hat w.
$$
Next, as in \cite{nd7}, we are going to establish the classical sup+inf-estimate for $\hat u_l$ (and even the full version which will be done in the next Lemma), since thanks \eqref{eq:noinflu} all terms coming from $\hat w$ can be absorbed on the right hand side of \eqref{eq:classupinf}.
Now, we are going to rescale the functions $\hat u_l$ around the points $x$.... In order to do that, we define $\varphi_l: B^{\R^4}_0(2\eta_1\mu_l^{-1})\longrightarrow B^{\hat g}_{x}(2\eta_1)$ if $x\in M$, and $\varphi_l: B^{\R^4_+}_0(2\eta_1\mu_l^{-1})\longrightarrow B^{\hat g}_{x}(2\eta_1)$ if $x\in \partial M$  by the formula $\varphi_l(z):=\mu_lz$ and $\mu_l$ is the corresponding scaling parameter given by Lemma \ref{eq:locdesblow}. Furthermore, as in \cite{nd7}, we define the following rescaling of $\hat u_l$
$$
v_l:=\hat {u_l}\circ\varphi_l+\ln  \mu_l+\frac{1}{4}\ln \frac{t_l K(x)}{3}.
$$
Using the Green's representation formula for $(P^4_{\hat g}, P^3_{\hat g})$ under homogenous Neumann boundary condition with respect to $\hat g$, the method of \cite{nd4} and standard doubling argument to deal with the situation $x\in \partial M$, we get, by
\begin{equation}\label{eq:estint}
v_l(z)+2\ln |z|=O(1), \;\;\text{for}\;\;z\in \bar B^{\R^4}_0(\frac{\eta_1}{\mu_l})\setminus B^{\R^4}_0(-\ln \mu_l),
\end{equation}
in case of interior blow-up, and in case of boundary blow-up, we obtain
\begin{equation}\label{eq:estboun1}
v_l(z)+2\ln |z|=O(1), \;\;\text{for}\;\; z\in \bar B^{\R^4_+}_0(\frac{\eta_1}{\mu_l})\setminus B^{\R^4_+}_0(-\ln \mu_l).
\end{equation}
Now, we are going to show that the estimate \eqref{eq:estint} holds also in $\bar B^{\R^4}_0(-\ln  \mu_l)$ (in case of interior blow-up), and  that the estimate \eqref{eq:estboun1} holds in $\ov B^{\R^4_+}_0(-\ln \mu_l)$ as well (in case of boundary blow-up). To do so, we use Lemma \ref{eq:locdesblow}, the same arguments as in \cite{nd7}, and standard doubling argument (when $x\in \partial M$), to get
\begin{equation}\label{eq:euclideanversion}
v_l(z)+2\ln |z|=O(1), \;\;\text{for}\;\;\;z\in \bar B^{\R^4}_0(-\ln  \mu_l),\;\;\text{if}\;\;\;x\in M
\end{equation}
and
\begin{equation}\label{eq:euclideanversionb}
v_l(z)+2\ln  |z|=O(1), \;\;\text{for}\;\;\;z\in \bar B^{\R^4_+}_0(-\ln  \mu_l), \;\;\;\text{if}\;\;\;x\in \partial M.
\end{equation}
Now, combining \eqref{eq:estint} and \eqref{eq:euclideanversion} when $x\in M$, and \eqref{eq:estboun1} and \eqref{eq:euclideanversionb} when $x\in \partial M$, we obtain
\begin{equation}\label{eq:eucball}
v_l(z)+2\ln |z|=O(1), \;\;\text{for}\;\;\;z\in \bar B^{\R^4}_0(\frac{\eta_1}{\mu_l}),\;\;\text{if}\;\;\;x\in \mathring{M}
\end{equation}

\begin{equation}\label{eq:euchalf}
v_l(z)+2\ln |z|=O(1), \;\;\text{for}\;\;\;z\in \bar B^{\R^4_+}_0(\frac{\eta_1}{\mu_l}), \;\;\text{if}\;\;\;x\in \partial M.
\end{equation}
Thus scaling back, namely using $y=\mu_lz$ and the definition of $v_l$, we obtain the desired $O(1)$-estimate.  Hence the proof of the Lemma is complete.
\end{pf}
\subsection{Refined estimate around blow-up points}
As already mentioned above, in this subsection, we show formula \eqref{eq:sharpestimate}. Precisely, we prove the following Lemma:
\begin{lem}\label{eq:strongerestimate}
Assuming that $(u_l)_{l\in \N}$ is a bubbling sequence of solutions to BVP \eqref{eq:sepexact} with $t_l$ satisfying \eqref{eq:lbqexact}, then keeping the notations in Lemma \ref{eq:locdesblow}, Lemma \ref{eq:isolated}, and Lemma \ref{eq:roughestimate}, we have that the following estimate holds ($l$ large enough)
$$
u_l(x)+\frac{1}{4}\ln\frac{t_l K(x_{i})}{3}=\ln \frac{2\l_{i, l}}{1+\l_{i, l}^2(d_{g_{x_i}}(x, x_{i}))^2}+ O(d_{g_{x_i}}(x, x_{i})),\;\;\;\;\forall\;x\in\;B_{x_{i}}^{x_i}(\eta_k).
$$
\end{lem}
\begin{pf}
We are going to use the method of \cite{nd7}, hence we will be sketchy in many arguments. Now, let $V_0$ be the unique solution of the following conformally invariant integral equation
$$
V_0(z)=\frac{3}{4\pi^2}\int_{\R^4}\ln \frac{|y|}{|z-y|}e^{4V_0(y)}dy+\ln 2,\;\;\; V_0(0)=\ln  2, \n V_0(0)=0.
$$
Next, we set $w_l(z)=v_l(z)-V_0(z)$ for $z\in B^{\R^4}_0(\eta_1\mu_l^{-1})$ when $x\in M$, and $w_l(z)=v_l(z)-V_0(z)$ for $z\in \bar B^{\R^4}_0(\eta_1\mu_l^{-1})$ when $x\in \partial M$, and use Lemma \ref{eq:roughestimate} \;to infer that
\begin{equation}\label{eq:weakes}
|w_l|\leq C\;\;\;\,\text{in}\;\;\;\; B^{\R^4}_0(\eta_1\mu_l^{-1})\;\;\;\text{if}\;\:\;x\in \mathring{M}.
\end{equation}
and
\begin{equation}\label{eq:weaksb}
|w_l|\leq C,\:\;\;\;\text{in}\;\;\;\; B^{\R^4_+}_0(\eta_1\mu_l^{-1})\;\;\;\text{if}\;\,\;x\in \partial M.
\end{equation}
On the other hand, it is easy to see that to achieve our goal, it is sufficient to show
\begin{equation}\label{eq:3.2}
|w_l|\leq C\mu_l|z|\,\,\,\,\text{in}\,\,\,\, B^{\R^4}_0(\eta_1\mu_l^{-1}),\;\;\;\text{if}\;\:\;x\in \mathring{M}.
\end{equation}
and
\begin{equation}\label{eq:3.21}
|w_l|\leq C\mu_l|z|\,\,\,\,\text{in}\,\,\,\, B^{\R^4_+}_0(\eta_1\mu_l^{-1}),\;\;\;\text{if}\;\:\;x\in \partial M.
\end{equation}
To show \eqref{eq:3.2} and \eqref{eq:3.21}, we first set
$$
\L_l:=\max_{z\in \Omega_l}\frac{|w_l(z)|}{\mu_l(1+|z|)}
$$
with
$$
\Omega_l=\ov B^{\R^4}_0(\eta_1\mu_l^{-1}), \;\;\;\text{if}\;\:\;x\in \mathring{M}.,
$$
and
$$
\Omega_l=\ov B^{\R^4_+}_0(\eta_1\mu_l^{-1}), \;\;\;\text{if}\;\:\;x\in \partial M.
$$
 We remark that to show \eqref{eq:3.2} and \eqref{eq:3.21}, it is equivalent to prove that $\L_l$ is bounded.  Now, let us suppose that $\L_l\rightarrow +\infty$ as $l\rightarrow +\infty$, and look for a contradiction. To do so, we will use the method of  \cite{nd4} combined with a standard doubling argument. For this, we first choose  a sequence of points $z_l\in \Omega_l$ such that $\L_l=\frac{|w_l(z_l)|}{\mu_l(1+|z_l|)}$. Next, up to a subsequence, we have that either $z_l\rightarrow z^{*}$ as $l\rightarrow+\infty$ (with $z^{*}\in \R^4$) or $|z_l|\rightarrow +\infty$ as $l\rightarrow+\infty$. Now, we make the following definition
$$
\bar w_l(z):=\frac{w_l(z)}{\L_l\mu_l(1+|z_l|)},
$$
and have
\begin{equation}\label{eq:3.3}
|\bar w_l(z)|\leq \left(\frac{1+|z|}{1+|z_l|}\right),
\end{equation}
and
\begin{equation}\label{eq:contradictionc}
|\bar w_l(z_l)|=1.
\end{equation}
Now, we consider the case where the points $z_l$ escape to infinity.\\\\
{\em Case 1 }: $|z_l|\rightarrow +\infty$\\
In this case, using \eqref{eq:weakes}, \eqref{eq:weaksb}, the Green's representation formula for $(P^4_{\hat g}, P^3_{\hat g})$ under homogeneous Neumann boundary condition with respect to $\hat g$, and the method of \cite{nd7}, we obtain
$$
\bar w_l(z_l)=\frac{3}{4\pi^2}\int_{\Omega_l}\ln \frac{|\xi|}{|z_l-\xi|}\left(\frac{O(1)(1+|\xi|)^{-7}}{(1+|z_l|)}+\frac{O(1)(1+|\xi|)^{-7}}{\L_l(1+|z_l|)}\right)d\xi+o(1).
$$
Now, using the fact that $|z_l|\rightarrow+\infty$ as $l\rightarrow +\infty$, one can easily check that
$$
\bar w_l(z_l)=\frac{3}{4\pi^2}\int_{\Omega_l}\ln  \frac{|\xi|}{|z_l-\xi|}\left(\frac{O(1)(1+|\xi|)^{-7}}{(1+|z_l|)}+\frac{O(1)(1+|\xi|)^{-7}}{\L_l(1+|z_l|)}\right)d\xi + o(1).
$$
Hence, we reach a contradiction to \eqref{eq:contradictionc}. \\
Now, we are going to show that, when the points $z_l\rightarrow z^{*}$ as $l\rightarrow +\infty$, we reach a contradiction as well.\\\\
{\em Case 2}: $z_l\rightarrow z^{*}$\\
In this case, using the assumption $z_l\longrightarrow z^*$, the Green's representation formula for $(P^4_{\hat g}, P^3_{\hat g})$ under homogeneous Neumann boundary condition with respect to $\hat g$, and the method of \cite{nd7}, we obtain that up to a subsequence
\begin{equation}\label{eq:lisupinf}
\bar w_l\rightarrow w \;\;\;\text{in}\;\;C^1_{loc}(\R^4) \;\;\;\text{as}\;\;\; l\rightarrow +\infty,\;\;\:\text{if}\;\;\;x\in  \mathring{M},
\end{equation}
\begin{equation}\label{eq:lisupinfb}
\bar w_l\rightarrow w \;\;\;\text{in}\;\;C^1_{loc}(\bar\R^4_+) \;\;\;\text{as}\;\;\; l\rightarrow +\infty, \;\;\partial_{x_4}w=0\;\;\text{on}\;\R^3,\;\;\;\text{if}\;\;\;x\in\partial  M,
\end{equation}
and
\begin{equation}\label{eq:barwl1b}
\begin{split}
\bar w_l(z)=\frac{3}{\pi^2}\int_{\Omega_l} \ln \frac{|\xi|}{|z-\xi|}\frac{K\circ\varphi_l(\xi)}{K\circ\varphi_l(0)}e^{4\vartheta_l(\xi)}\bar w(\xi)d\xi+\frac{3}{\L_l\mu_l(1+|z_l|)\pi^2}\int_{\Omega_l} \ln \frac{|\xi|}{|z-\xi|}O(\mu_l(1+|\xi|)^{-7})d\xi\\+\frac{O(1)+O(|z|)}{\L_l(1+|z_l|)},
\end{split}
\end{equation}
where $e^{4\theta_l}:=\int_0^1e^{4(sv_l+(1-s)V_0)}ds$.
Thus, appealing to \eqref{eq:lisupinf}, \eqref{eq:lisupinfb}, and \eqref{eq:barwl1b}, we infer that $w$ satisfies
\begin{equation}\label{eq:wlim1}
w(z)=\frac{3}{\pi^2}\int_{\R^4}\ln \frac{|\xi|}{|z-\xi|}e^{4V_0(\xi)}w(\xi)d\xi,\;\;\;\text{if}\;\;\;x\in \partial{M},
\end{equation}
and
\begin{equation}\label{eq:wlim2}
w(z)=\frac{3}{\pi^2}\int_{\R^4}\ln \frac{|\xi|}{|z-\xi|}e^{4V_0(\xi)}w(\xi)d\xi,\;\;\partial_{x_4}w=0\;\;\text{on}\;\R^3\;\;\;\;\text{if}\;\;\;x\in \partial M.
\end{equation}
Now, using \eqref{eq:3.3}, we have that $w$ satisfies the following asymptotic
\begin{equation}\label{eq:continf}
|w(z)|\leq C(1+|z|).
\end{equation}
On the other hand, from the definition of $v_l$,  it is easy to see that
\begin{equation}\label{eq:vanord1}
w(0)=0, \;\;\;\text{and}\;\;\;\;\n w(0)=0.
\end{equation}
So, using \eqref{eq:wlim1}-\eqref{eq:vanord1}, Lemma 3.7 in \cite{nd7}, and a standard doubling argument, we obtain
$$
w=0.
$$
 However, from \eqref{eq:contradictionc}, we infer that $w$ satisfies also
 \begin{equation}\label{eq:nonzero}
 |w(z^*)|=1.
 \end{equation}
So we reach a contradiction in the second case also. Hence the proof of the Lemma is complete.
\end{pf}
\vspace{6pt}

\noindent
\begin{pfn}{ of Proposition \ref{eq:sharpest}}\\
Proposition \ref{eq:sharpest} follows directly from Lemma \ref{eq:isolated}, Lemma \ref{eq:roughestimate}, Remark \ref{eq:comparability}, and Lemma \ref{eq:strongerestimate}.
\end{pfn}






\section{ A Morse lemma at infinity }\label{s:morseleminf}
In this section, we characterize the critical points at infinity of $II$  and establish a Morse type Lemma for them. To do so, we will divide this section into two subsections. In the first one, we perform a finite-dimensional Lyapunov-Schmidt type reduction. In the second one, we combine the latter finite-dimensional reduction and the construction of a suitable pseudogradient at infinity to achieve our goal. \\
For all we will parameterize the neighborhood of potential critical points at infinity. Namely following the ideas of Bahri-Coron \cite{bc}, and using Lemma \ref{eq:intbubbleest} and Lemma \ref{eq:outbubbleest}, we have that for every $\varrho$ and  $\eta$ small positive real numbers with $0<\varrho<\varrho_k$, and $0<\eta<\eta_k$ where $\varrho_k$ and $\eta_k$ are given by Proposition \ref{eq:sharpest}, there exists $\epsilon_k=\epsilon_k(\varrho, \eta)>0$ such that for every $(p, q)\in \N^2$ with $2p+q=k$, there holds
\begin{equation}\label{eq:mini}
\begin{split}
\forall\;0<\epsilon\leq \epsilon_k,\;\;\forall u\in V(p, q, \epsilon, \varrho, \eta), \text{the minimization problem }\\\min_{B_{\epsilon, \varrho, \eta}^{p, q}}||u-\ov{u}_{(Q, T)}-\sum_{i=1}^{p+q}\alpha_i\varphi_{a_i, \l_i}-\sum_{r=1}^{\bar k}\beta_r(v_r-\ov{(v_r)}_{(Q, T)})||
\end{split}
\end{equation}
has a unique solution, up to permutations, where $B^{p, q}_{\epsilon, \varrho, \eta}$ is defined as follows
\begin{equation}
\begin{split}
{B_{\epsilon, \varrho, \eta}^{p, q}:=\{(\bar\alpha, A, \bar \l, \bar \beta)\in \R^{p+q}\times \mathring{M}^p\times (\partial M)^q\times (0, +\infty)^{p+q}\times \R^{\bar k}}:|\alpha_i-1|\leq \epsilon, \l_i\geq \frac{1}{\epsilon}, i=1, \cdots, p+q, \\d_g(a_i, a_j)\geq 4\ov C\eta, i\neq j=1, \cdots, p+q, \;d_g(a_i, \partial M)\geq 4\ov C\varrho,\;\; |\beta_r|\leq R, r=1, \cdots, \bar k\}.
\end{split}
\end{equation}
 Moreover, using the solution of \eqref{eq:mini}, we have that every $u\in V(p, q, \epsilon, \varrho, \eta)$ can be written as
\begin{equation}\label{eq:para}
u-\ov{u}_{(Q, T)}=\sum_{i=1}^{p+q}\alpha_i\varphi_{a_i, \l_i}+\sum_{r=1}^{\bar k}\beta_r(v_r-\ov{(v_r)}_{(Q, T)})+w,
\end{equation}
where $w$ verifies the following orthogonality conditions
\begin{equation}\label{eq:ortho}
\begin{split}
 \ov w_{(Q, T)}=<\varphi_{a_i, \l_i}, w>_{P^{4, 3}}=<\frac{\partial\varphi_{a_i, \l_i}}{\partial \l_i}, w>_{P^{4, 3}}=<\frac{\partial\varphi_{a_i, \l_i}}{\partial a_i}, w>_{P^{4, 3}}=<v_r, w>=0, i=1, \cdots, p+q,\\r=1, \cdots, \bar k
\end{split}
\end{equation}
and the estimate
\begin{equation}\label{eq:estwmin}
||w||=O\left(\sum_{i=1}^{p+q}\frac{1}{\l_i}\right),
\end{equation}
where here $O\left(1\right):=O_{\bar \alpha, A, \bar \l, \bar\beta , w, \epsilon}\left(1\right)$, and for the meaning of $O_{\bar \alpha, A, \bar \l, \bar\beta , w, \epsilon}\left(1\right)$, see Section \ref{eq:notpre}. Furthermore, the concentration points $a_i$,  the masses $\alpha_i$, the concentrating parameters $\l_i$ and the negativity parameter $\beta_r$ in \eqref{eq:para} verify also
\begin{equation}\label{eq:afpara}
\begin{split}
d_g(a_i, a_j)\geq 4\ov C\eta\;\;\text{for}\;\;i\neq j=1, \cdots, p+q, \;\;d_g(a_i, \partial M)\geq 4\ov C\varrho \;\;\text{for}\;\;i=p+1, \cdots, p+q, \\\;\;\frac{1}{\L}\leq\frac{\l_i}{\l_j}\leq\L\;\;\text{for}\;\;i, j=1, \cdots, p+q, \;\;\l_i\geq\frac{1}{\epsilon}\;\;\text{for}\;\;i=1, \cdots, p+q,\;\;\\\text{and}\;\;\;\sum_{r=1}^{\bar k}|\beta_r|+\sum_{i=1}^{p+q}|\alpha_i-1|\sqrt{\ln \l_i}=O\left(\sum_{i=1}^{p+q}\frac{1}{\l_i}\right),
\end{split}
\end{equation}
with still $O\left(1\right)$ as in \eqref{eq:estwmin}.




\subsection{Finite-dimensional reduction near  infinity}\label{s:finitedim}
In this subsection, we perform a finite-dimensional Lyapunov-Schmidt type reduction by exploiting the stability property of the standard bubbles as in \cite{no1}. Indeed, in doing so, we first derive the following proposition...
\begin{pro}\label{eq:expansionJ1}
Assuming that $(p, q)\in \N^2$  such that $2p+q=k$, $0<\varrho<\varrho_k$, $0<\eta<\eta_k$, and $0<\epsilon\leq \epsilon_k$, where $\varrho_k$ and $\eta_k$ are given by Proposition \ref{eq:sharpest}, and $\epsilon_k$ is given by  \eqref{eq:mini}, and $u=\ov{u}_{(Q, T)}+\sum_{i=1}^{p+q}\alpha_i\varphi_{a_i, \l_i}+\sum_{r=1}^{\bar k}\beta_r (v_r-\ov{(v_r)}_{(Q, T)})+w\in V(p, q, \epsilon, \varrho \eta)$ with $w$, the concentration points $a_i$,  the masses $\alpha_i$, the concentrating parameters $\l_i$  ($i=1, \cdots, p+q$), and the negativity parameters $\beta_r$ ($r=1, \cdots, \bar k$) verifying \eqref{eq:ortho}-\eqref{eq:afpara}, then we have
\begin{equation}\label{eq:exparoundbubble}
II(u)=II(\sum_{i=1}^{p+q}\alpha_i\varphi_{a_i, \l_i}+\sum_{r=1}^{\bar k}\beta_r (v_r-\ov{(v_r)}_{(Q, T)})-f(w)+Q(w)+o(||w||^2),
\end{equation}
where
\begin{equation}\label{eq:linear}
f(w):=16k\pi^2\frac{\int_M Ke^{4\sum_{i=1}^{p+q}\alpha_i\varphi_{a_i, \l_i}+4\sum_{r=1}^{\bar k}\beta_r v_r}wdV_g}{\int_M Ke^{4\sum_{i=1}^{p+q}\alpha_i\varphi_{a_i, \l_i}+4\sum_{r=1}^{\bar k}\beta_r v_r}dV_g},
\end{equation}
and
\begin{equation}\label{eq:quadratic}
Q(w):=||w||^2-32k\pi^2\frac{\int_M Ke^{4\sum_{i=1}^{p+q}\alpha_i\varphi_{a_i, \l_i}+4\sum_{r=1}^{\bar k}\beta_r v_r}w^2dV_g}{\int_M Ke^{4\sum_{i=1}^{p+q}\alpha_i\varphi_{a_i, \l_i}+4\sum_{r=1}^{\bar k}\beta_r v_r}dV_g}.
\end{equation}
Moreover, setting
\begin{equation}\label{eq:eali}
\begin{split}
E_{a_i, \l_i}:=\{w\in \mathcal{H}_{\frac{\partial}{\partial n}}: \;\;<\varphi_{a_i, \l_i}, w>_{P^{4, 3}}=<\frac{\partial\varphi_{a_i, \l_i}}{\partial \l_i}, w>_{P^{4, 3}}=<\frac{\partial\varphi_{a_i, \l_i}}{\partial a_i}, w>_{P^{4, 3}}=0,\\\,\;\,\;\ov{w}_{(Q, T)}=<v_k, w>=0, \;k=1, \cdots, \bar k,\,\;\text{and}\;\;||w||=O\left(\sum_{i=1}^{p+q}\frac{1}{\l_i}\right)\},
\end{split}
\end{equation}
and
\begin{equation}\label{eq:eal}
A:=(a_1, \cdots, a_{p+q}), \;\;\bar \l=(\l_1, \cdots, \l_{p+q}), \;\;E_{A, \bar \l}:=\cap_{i=1}^{p+q} E_{a_i, \l_i},
\end{equation}
we have that, the quadratic form $Q$ is positive definite in $E_{A, \bar \l}$. Furthermore, the linear par $f$ verifies that, for every $w\in E_{A, \bar \l}$, there holds
\begin{equation}\label{eq:estlinear}
f(w)=O\left( ||w||\left(\sum_{i=1}^{p}\frac{|\n_g \mathcal{F}^{A}_i(a_i)|}{\l_i}+     \sum_{i= p + 1}^{p+q}\frac{|\n_{\hat{g}} \mathcal{F}^{A}_i(a_i)|}{\l_i}   +   \sum_{i=1}^{p+q}|\alpha_i-1| \, + \sum_{i={p+1}}^{p+q} \frac{  |\frac{\partial \mathcal{F}_i^A }{\partial n_g}(a_i)|}{\l_i} +\sum_{r=1}^{\bar k}|\beta_r|  + \sum_{i=1}^{p+q}  \frac{\ln \l_i}{\l^2_i}\right)\right).
\end{equation}
where for $i=1, \cdots, p+q$, $\mathcal{F}_i^A$ is as in Lemma \ref{eq:energyest}, $o(1)=o_{ \bar\alpha, A, \bar \beta, \bar\l, w, \epsilon}(1)$, $O\left(1\right):=O_{\bar\alpha, A, \bar \beta, \bar\l, w, \epsilon}\left(1\right)$, and for the meaning of $o_{ \bar\alpha, A, \bar \beta, \bar\l, w, \epsilon}(1)$ and $O_{\bar\alpha, A, \bar \beta, \bar\l, w, \epsilon}\left(1\right)$, see Section \ref{eq:notpre}.
\end{pro}
\begin{pf}
 It follows from Lemma \ref{eq:gradientlambdaest}, Lemma \ref{eq:intbubbleest}, Lemma \ref{eq:positive}, Lemma \ref{eq:positiveg}, and the same strategy as in the proof of Proposition 6.1 in \cite{no1}.
\end{pf}

\vspace{4pt}

\noindent
As in \cite{no1}, we have that Proposition \ref{eq:expansionJ1} implies the following direct corollaries.
\begin{cor}\label{eq:cexpansionj1}
Assuming that $(p, q)\in \N^2$  such that $2p+q=k$, $0<\varrho<\varrho_k$, $0<\eta<\eta_k$, and $0<\epsilon\leq \epsilon_k$, where $\varrho_k$ and $\eta_k$ are given by Proposition \ref{eq:sharpest}, and $\epsilon_k$ is given by \eqref{eq:mini}, and $u:=\sum_{i=1}^{p+q}\alpha_i\varphi_{a_i, \l_i}+\sum_{r=1}^{\bar k}\beta_r(v_r-\ov{(v_r)}_{(Q, T)})$ with the concentration  points $a_i$,  the masses $\alpha_i$, the concentrating parameters $\l_i$ ($i=1, \cdots, p+q$)  and the negativity parameters $\beta_r$ ($r=1, \cdots, \bar k$) satisfying \eqref{eq:afpara}, then there exists a unique $\bar w(\bar \alpha, A, \bar \l, \bar \beta)\in E_{A, \bar \l}$ such that
\begin{equation}\label{eq:minj}
II(u+\bar w(\bar \alpha, A, \bar \l, \bar \beta))=\min_{w\in E_{A, \bar \l}, u+w\in V(p, q, \epsilon, \varrho, \eta)} II(u+w),
\end{equation}
where $\bar \alpha:=(\alpha_1, \cdots, \alpha_{p+q})$, $A:=(a_1, \cdots, a_{p+q})$, $\bar \l:=(\l_1, \cdots, \l_{p+q})$ and $\bar \beta:=(\beta_1, \cdots, \beta_{\bar k})$.\\
Furthermore, $(\bar \alpha, A, \bar \l, \bar \beta)\longrightarrow \bar w(\bar \alpha, A, \bar \l, \bar \beta)\in C^1$ and satisfies the following estimate
\begin{equation}\label{eq:linminqua}
\frac{1}{C}||\bar w(\bar \alpha, A,  \bar \l, \bar \beta)||^2\leq |f(\bar w(\bar \alpha, A, \bar \l, \bar \beta))|\leq C|| \bar w(\bar \alpha, A, \bar \l, \bar \beta)||^2,
\end{equation}
for some large positive constant $C$ independent of $\bar \alpha$, $A$, $\bar \l$, $\bar \beta$, and $\epsilon
$, hence
\begin{equation}\label{eq:estbarw}
||\bar w(\bar \alpha, A, \bar \l, \bar \beta)||=O\left(\sum_{i=1}^{p}\frac{|\n_g \mathcal{F}^{A}_i(a_i)|}{\l_i}    + \sum_{i= p +1}^{p + q}\frac{|\n_{\hat{g}} \mathcal{F}^{A}_i(a_i)|}{\l_i}        +  \sum_{i=p+1}^{p+q} \frac{ |\frac{\partial \mathcal{F}_i^A (a_i)}{\partial n_g}|  }{\l_i}    + \sum_{i=1}^{p+q}|\alpha_i-1|\, + \sum_{r=1}^{\bar k}|\beta_r|+\sum_{i=1}^{p+q}\frac{\ln \l_i}{\l_i^{2}}\right),
\end{equation}
where $O\left(1\right):=O_{\bar\alpha, A, \bar \beta, \bar\l, \epsilon}\left(1\right)$ and for its meaning see Section \ref{eq:notpre}.
\end{cor}
\vspace{4pt}

\noindent
\begin{cor}\label{eq:c1expansionj1}
Assuming that $(p, q)\in \N^2$  such that $2p+q=k$, $0<\varrho<\varrho_k$, $0<\eta<\eta_k$, and $0<\epsilon\leq \epsilon_k$, where $\varrho_k$ and $\eta_k$ are given by Proposition \ref{eq:sharpest}, and $\epsilon_k$ is given by \eqref{eq:mini}, and $u_0:=\sum_{i=1}^{p+q}\alpha_i^0\varphi_{a_i^0, \l_i^0}+\sum_{r=1}^{\bar k}\beta_r^0(v_r-\ov{(v_r)}_{(Q, T)})$ with the concentration  points $a_i^0$,  the masses $\alpha_i^0$, the concentrating parameters $\l_i^0$ ($i=1, \cdots, p+q$) and the negativity parameters $\beta_r^0$ ($r=1, \cdots, \bar k$) satisfying \eqref{eq:afpara}, then there exists an open neighborhood  $U$ of $(\bar\alpha^0, A^0, \bar \l^0, \bar \beta^0)$  (with $\bar \alpha^0:=(\alpha^0_1, \cdots, \alpha^0_{p+q})$, $A^0:=(a_1^0, \cdots, a^0_{p+q})$,  $\bar \l:=(\l_1^0, \cdots, \l_{p+q}^0)$ and  $\bar \beta^0:=(\beta_1^0, \cdots, \beta_{\bar k}^0)$) such that for every $(\bar \alpha, A, \bar \l, \bar \beta)\in U$ with $\bar \alpha:=(\alpha_1, \cdots, \alpha_{p+q})$, $A:=(a_1, \cdots, a_{p+q})$,  $\bar \l:=(\l_1, \cdots, \l_{p+q})$, $\bar \beta:=(\beta_1, \cdots, \beta_{\bar k})$, and the $a_i$,  the $\alpha_i$, the $\l_i$ ($i=1, \cdots, p+q$)  and the $\beta_r$ ($r=1, \cdots, \bar k$) satisfying \eqref{eq:afpara}, and $w$ satisfying \eqref{eq:afpara}  with $\sum_{i=1}^{p+q}\alpha_i\varphi_{a_i, \l_i}+\sum_{r=1}^{\bar k}\beta_r(v_r-\ov{(v_r)}_{(Q, T)})+w\in V(p, q, \epsilon, \varrho, \eta)$, we have the existence of a change of variable
\begin{equation}\label{eq:changev}
w\longrightarrow V
\end{equation}
from a neighborhood of $ \bar w(\bar \alpha, A, \bar \l, \bar \beta)$ to a neighborhood of $0$ such that
\begin{equation}\label{eq:expjv}
\begin{split}
&II(\sum_{i=1}^{p+q}\alpha_i\varphi_{a_i, \l_i}+\sum_{r=1}^{\bar k}\beta_r(v_r-\ov{(v_r)}_{(Q, T)})+w)=\\&II(\sum_{i=1}^{p+q}\alpha_i\varphi_{a_i, \l_i}+\sum_{r=1}^{\bar k}\beta_r(v_r-\ov{(v_r)}_{(Q, T)})+\bar w (\bar \alpha, A, \bar \l, \bar \beta))\\&+\frac{1}{2}\partial^2 II(\sum_{i=1}^{p+q}\alpha_i^0\varphi_{a_i^0, \l_i^0}+\sum_{r=1}^{\bar k}\beta_r^0(v_r-\ov{(v_r)}_{(Q, T)})+\bar w(\bar \alpha^0, A^0, \bar \l^0, \bar \beta^0))(V, V),
\end{split}
\end{equation}
\end{cor}
\vspace{4pt}
Thus, as in \cite{no1},  in $V(p, q, \epsilon, \varrho, \eta)$ we have a splitting of the variables $(\bar \alpha, A, \bar \l, \bar \beta)$ and $V$ and one can decrease the functional $II$ in the variable $V$  without touching the variable $(\bar \alpha, A, \bar \l, \bar \beta)$  by considering just the flow
\begin{equation}\label{eq:vflow}
\frac{dV}{dt}=-V.
\end{equation}
So, since $II$ is invariant by translations by constants, then the variational study of $II$ in $V(p, q, \epsilon, \varrho,  \eta)$ is equivalent to the one of the following finite-dimensional functional \begin{equation}\label{eq:finitedj}
II(\bar \alpha, A, \bar \l, \bar \beta):=II(\sum_{i=1}^{p+q}\alpha_i\varphi_{a_i, \l_i}+\sum_{r=1}^{\bar k}\beta_r(v_r-\ov{(v_r)}_{(Q, T)})+\bar w(\bar \alpha, A, \bar \l, \bar \beta)),
\end{equation}
where $\bar \alpha=(\alpha_1, \cdots, \alpha_{p+q})$, $A=(a_1, \cdots, a_{p+q})$, $\bar \l=(\l_1, \cdots, \l_{p+q})$ and $\bar \beta=\beta_1, \cdots, \beta_{\bar k}$ with the concentration  points $a_i$,  the masses $\alpha_i$, the concentrating parameters $\l_i$ ($i=1, \cdots, p+q$)  and the negativity parameters $\beta_r$ ($r=1, \cdots, \bar k$) satisfying \eqref{eq:afpara}, and $\bar w(\bar \alpha, A, \bar \l, \bar \beta)$ is as in Corollary \ref{eq:cexpansionj1}. Hence the goal of this subsection is achieved.

\subsection{Construction of a pseudogradient near  infinity}\label{s:cpseudo}
In this subsection, we construct a pseudogradient for the finite-dimensional functional $II(\bar \alpha, A, \bar \l, \bar \beta)$ given by \eqref{eq:finitedj}, and use it to characterize the critical points at infinity of $II$ and establish a Morse type Lemma for them. Indeed, we have:
\begin{pro}\label{eq:conspseudograd}
Assuming that $(p, q)\in \N^2$  such that $2p+q=k$, $0<\varrho<\varrho_k$, $0<\eta<\eta_k$, and $0<\epsilon\leq \epsilon_k$, where $\varrho_k$ and $\eta_k$ are given by Proposition \ref{eq:sharpest}, and $\epsilon_k$ is given by \eqref{eq:energyest}, then there exists a pseudogradient $W$ of $II(\bar \alpha, A, \bar \l, \bar \beta)$ such that \\
1) For every $u:=\sum_{i=1}^{p+q}\alpha_i\varphi_{a_i, \l_i}+\sum_{r=1}^{\bar k}\beta_r (v_r-\ov{(v_r)}_{(Q, T)})\in V(p, q, \epsilon, \varrho, \eta)$ with the concentration  points $a_i$, the masses $\alpha_i$, the concentrating parameters $\l_i$ ($i=1, \cdots, p+q$)  and the negativity parameters $\beta_r$ ($r=1, \cdots, \bar k$) satisfying \eqref{eq:afpara}, we have that if $q=0$, then there holds
\begin{equation}\label{eq:pseudoexact}
<-\n II(u), W>\geq c\left(\sum_{i=1}^{p}\frac{1}{\l_i^2}+\sum_{i=1}^{p}\frac{|\n_g\mathcal{F}^{A}_i(a_i)|}{\l_i}+\sum_{i=1}^{p}|\alpha_i-1|+\sum_{i=1}^{p}|\tau_i|+\sum_{r=1}^{\bar k}|\beta_r|)\right),
\end{equation}
and if $q\neq 0$, then there holds
\begin{equation}\label{eq:pseudoperturbq}
\begin{split}
&<-\n II(u), W>\geq\\& c\left(\sum_{i=1}^{p+q}\frac{1}{\l_i}+\sum_{i=1}^p\frac{|\n_g\mathcal{F}^{A}_i(a_i)|}{\l_i}+\sum_{i=p+1}^{p+q}\frac{|\n_{\hat g}\mathcal{F}^{A}_i(a_i)|}{\l_i}+\sum_{i=1}^{p+q}|\alpha_i-1|+\sum_{i=1}^{p+q}|\tau_i|+\sum_{r=1}^{\bar k}|\beta_r|\right),
\end{split}
\end{equation}
where $c$ is a small positive constant independent of $A:=(a_1, \cdots, a_{p+q})$, $\bar\alpha=(\alpha_1, \cdots, \alpha_{p+q})$, $\bar\l=(\l_1, \cdots, \l_{p+q})$, $\bar \beta=(\beta_1, \cdots, \beta_{\bar k})$ and $\epsilon$.\\\\
2) Furthermore, for every $u:=\sum_{i=1}^{p+q}\alpha_i\varphi_{a_i, \l_i}+\sum_{r=1}^{\bar k}\beta_r (v_r-\ov{(v_r)}_{(Q, T)})+\bar w(\bar \alpha, A, \bar \l, \bar \beta)\in V(p, q, \epsilon, \varrho, \eta)$ with the concentration  points $a_i$, the masses $\alpha_i$, the concentrating parameters $\l_i$ ($i=1, \cdots, p+q$)  and the negativity parameters $\beta_r$ ($r=1, \cdots, \bar k$) satisfying \eqref{eq:afpara}, and $\bar w(\bar \alpha, A, \bar \l, \bar \beta)$ is as in \eqref{eq:minj},  we have that if $q=0$, then there holds
\begin{equation}\label{eq:pseudoperturb2}
<-\n II(u), W+\frac{\partial \bar w(W)}{\partial (\bar\alpha, A, \bar \l, \bar \beta)}>\geq c\left(\sum_{i=1}^{p}\frac{1}{\l_i^2}+\sum_{i=1}^{p}\frac{|\n_g\mathcal{F}^{A}_i(a_i)|}{\l_i}+\sum_{i=1}^{p}|\alpha_i-1|+\sum_{i=1}^{p}|\tau_i|+\sum_{r=1}^{\bar k}|\beta_r|)\right),
\end{equation}
and if $q\neq 0$, then there holds
\begin{equation}\label{eq:pseudoperturbq2}
\begin{split}
&<-\n II(u), W+\frac{\partial \bar w(W)}{\partial (\bar\alpha, A, \bar \l, \bar \beta)}>\geq\\& c\left(\sum_{i=1}^{p+q}\frac{1}{\l_i}+\sum_{i=1}^p\frac{|\n_g\mathcal{F}^{A}_i(a_i)|}{\l_i}+\sum_{i=p+1}^{p+q}\frac{|\n_{\hat g}\mathcal{F}^{A}_i(a_i)|}{\l_i}+\sum_{i=1}^{p+q}|\alpha_i-1|+\sum_{i=1}^{p+q}|\tau_i|+\sum_{r=1}^{\bar k}|\beta_r|\right),
\end{split}
\end{equation}
where $c$ is a small positive constant independent of $A:=(a_1, \cdots, a_{p+q})$, $\bar\alpha=(\alpha_1, \cdots, \alpha_{p+q})$, $\bar\l=(\l_1, \cdots, \l_{p+q})$, $\bar \beta=(\beta_1, \cdots, \beta_{\bar k})$ and $\epsilon$.\\
3) Moreover \;the pseudogradient \;$W$\; is a bounded vector field and the region where the concentration rates $\l_i$'s are not bounded along the flow lines of $W$ are the one where  $A:= (a_1, \cdots, a_{p+q} )$  converges along the flow lines of $W$ to a critical point $B$ of $\mathcal{F}_{p, q}$ satisfying $l_K(B)<0$.
\end{pro}
\begin{pf}\\
{\bf Case: $q=0$}\\
In this case, it follows from \eqref{eq:relationderivativeint}, \eqref{eq:deflA}, \eqref{eq:auxiindexa1}, Lemma \ref{eq:energyest}, Lemma \ref{eq:gradientlambdaest}, Corollary \ref{eq:cgradientlambdaest},  Proposition \ref{eq:conspseudograd} and the arguments  to derive it go along  with the ones in the proof of Proposition 8.1 in \cite{no1}.

\noindent
{\bf Case: $q\neq 0$} \\
We will provide the proof for the case $k \geq 2$. The case $k=1$ can be treated similarly.\\
Let $\zeta$ be a cut off function satisfying
$$
0 \leq \zeta \leq 1, \, \, \,  \zeta(s) = 0   \mbox{  if  }|s| \leq 1/2, \quad \quad \zeta(s) = 1 \mbox{ if } |s| \geq 1.
$$
The claimed global pseudogradient will be constructed as a convex combination of local ones. To derive these local pseudogradients, we make use of the following vector fields: \\
To move the concentration  rates we will  use of the vector fields
$$
   W_{\l_i}:= - \frac{1}{\a_i}\l_i \frac{\partial \varphi_i}{\partial  \l_i},
   \qquad
   W_{\tau_i}:= -  \frac{\tau_i}{|\tau_i|} \zeta \Big(\frac{\l_i}{\ln\l_i}|\tau_i| \Big)
\frac{1}{\a_i} \l_i \frac{\partial \varphi_i}{\partial \l_i}.
$$
To move the concentration points we make use of the vector field: \\
For $i = 1, \cdots, p$

$$
 W_{a_i} := \frac{\n \mathcal{F}^A_i(a_i)}{|\n
\mathcal{F}^A_i(a_i)|}
\zeta \Big(  \frac{ \l_i}{ \ln \l_i } | \n_{g}  \mathcal{F}^A_i(a_i)|   \Big)
\frac{1}{\l_i} \frac{\partial \varphi_i}{\partial a_i} .
$$
and for $i=p+1, \cdots, p+q$ we define

$$
W_{a_i} := \frac{\n_{\hat{g}} \mathcal{F}^A_i(a_i)}{|\n
\mathcal{F}^A_i(a_i)|}
\zeta \Big(  \frac{ \l_i}{ \ln \l_i } | \n_{g}  \mathcal{F}^A_i(a_i)|   \Big)
 \left [  \frac{1}{\l_i} \frac{\partial \varphi_i}{\partial a_i}    \, - \, \frac{3}{8} \,  \l_i \frac{ \partial \varphi_i}{\partial \l_i} \right ]
$$
We divide the set $V(p,q,\e,, \varrho, \eta )$ into  subsets. To that aim we define for $C > 0$ large to be chosen later the set
$$
\mathcal{V}_1(p,q,\e) :=  \big \{  u \in V(p,q,\e, \varrho, \eta), \mbox{ s.t. } \exists \, i \in \{1, \cdots, p+q    \big \}; \,   |\tau_i| \geq \frac{B}{\l_i}  \}
$$
 and setting  $F :=  \{    i \in \{1, \cdots, p+q  \} \mbox{ such thta }  |\tau_i| \geq \frac{M}{\l_i}  \}$, we define   in $\mathcal{V}_1(p,q,\e)$  the vector field
 $$
 W_1 := \sum_{ i \in F}  W_{\tau_i}.
 $$
Using \eqref{eq:gradientlambdaest} and taking $C$ large, we derive that for some  positive constant $c_1$ there holds :
$$
<- \n II(\ov{u}), W_1> \geq   c_1 \big (  \sum_{i \in F}  \frac{1}{\l_i}     \big ),
$$
hence the  claimed estimate holds  in this set since all $\l_i$ are comparable and according to  Lemma \eqref{l:Keu} we have that $\sum_i \tau_i = O( \sum \frac{1}{\l_i})$.\\
Now  for $\d$ a small positive number, we define the following subset
\begin{equation}
\begin{split}
\mathcal{V}_2(p,q,\e) := \big  \{  u \in V(p,q,\e, \varrho, \eta), \mbox{ s.t. } \forall  \, i \in \{1, \cdots, p+q     \}; \,   |\tau_i| \leq  \frac{ 2 C}{\l_i} \mbox{ and } \exists \, i \in \{ 1, \cdots, p   \}   \\
\mbox{ such that }  \,  |\n_g \mathcal{F}^A_i(a_i)| \geq \d
 \mbox{ or } \exists  \, i \in \{p+1, \cdots, p+q  \} \,  \mbox{ such that }  \,  |\n_{\hat{g}} \mathcal{F}^A_i(a_i)| \geq \d.   \big  \}
\end{split}
\end{equation}
Setting
$$
W_2 := \sum_{i=1}^{p+q} W_{a_i},
$$
we derive using \eqref{eq:gradientaest} that for some positive constant $c_2$ we have that:
$$
< - \n II(\ov{u}), W_2> \, \geq  \, c_2 \, \left \{    \sum_{i=1}^p \frac{|\n_g  \mathcal{F}_i^A(a_i)|}{\l_i}  \, + \, \sum_{i= p + 1}^{p + q } \frac{|\n_{\hat{g}}  \mathcal{F}_i^A(a_i)|}{\l_i}   \right \},
$$
hence the claimed estimate holds  also in this region. \\
Now we define the following region
\begin{equation}
\begin{split}
\mathcal{V}_3(p,q,\e) := \big  \{  u \in V(p,q,\e, \varrho, \eta), \mbox{ s.t. } \forall  \, i \in \{1, \cdots, p+q     \}; \,   |\tau_i| \leq  \frac{ 2C}{\l_i} \mbox{ and } \forall \, i \in \{ 1, \cdots, p   \}   \\
\mbox{ such that }   |\n_g \mathcal{F}^A_i(a_i)| < 2 \d
 \mbox{ and  } \forall \,  i \in \{p+1, \cdots, p+q  \}  \mbox{ such that }   |\n_{\hat{g}} \mathcal{F}^A_i(a_i)| < 2  \d.   \big  \}
\end{split}
\end{equation}
We observe that in this region the concentration points are in an arbitrarily small $\d-$neighborhood of some critical point $B$ of  $\mathcal{F}_{p,q}$ and we subdivide it in two subsets
$$
\mathcal{V}_3^-(p,q, \e) \, =  \left \{   u \in  \mathcal{V}_3(p,q, \e) \mbox{ such that }   \mathcal{L}_K (B) < 0,   \right  \}
$$
 and
 $$
\mathcal{V}_3^+(p,q, \e) \, =  \left \{   u \in  \mathcal{V}_3(p,q, \e) \mbox{ such that }   \mathcal{L}_K (B) >  0,   \right  \}
$$
 In the set $\mathcal{V}_3(p,q, \e)$ we define the vector field
 $$
 W_3 := \, - sign(\mathcal{L}_K (B)) \sum_{i=1}^{p+q}  W_{\l_i}.
 $$
Using \eqref{eq:cgradientlambdaest} we derive that for some positive constants $c_3 $ there holds
$$
<-\n II(\ov{u}, W_3)> \, \geq c \sum_{i=1}^{p+q}  \frac{1}{\l_i},
$$
hence the claimed estimate holds also in this region and the global pseudogradient  $W$ is a convex combination of $W_i, \, i =1, 2,3.$ \\
The proof of the second claim  \eqref{eq:pseudoperturbq2} follows from \eqref{eq:pseudoperturbq} using the fact that
$$
||\ov{w}||^2  \, = \, o\left(\sum_{i=1}^{p+q}\frac{1}{\l_i}+\sum_{i=1}^p\frac{|\n_g\mathcal{F}^{A}_i(a_i)|}{\l_i}+\sum_{i=p+1}^{p+q}\frac{|\n_{\hat g}\mathcal{F}^{A}_i(a_i)|}{\l_i}+\sum_{i=1}^{p+q}|\alpha_i-1|+\sum_{i=1}^{p+q}|\tau_i|+\sum_{r=1}^{\bar k}|\beta_r|\right)
$$

\noindent
Now we observe that the only region where the $\l_i$'s are not bounded along the flow lines of $W$ is the region $\mathcal{V}_3^-$ where the concentration points converge to some  critical point $B$ of $\mathcal{F}_{p,q}$ such that $\mathcal{L}_K (B) < 0.$
\end{pf}
\vspace{6pt}

\noindent
Now we define the notion of \emph{ critical points at infinity } for $II$
\begin{df}
A { \it critical point at infinity} for  $II$, with respect to the pseudogradient $W$  is an accumulation of some non compact orbits of $W$ such that the flow lines enter and  remain for ever in some $V(p,q,\e_k, \varrho, \eta)$ for $\e_k \to 0.$
\end{df}

\noindent
As a corollary of Proposition \ref{eq:conspseudograd}, we derive  the following characterization of the critical points at infinity of $II$.
\begin{cor}\label{eq:loccritinf}
 1) The  critical points at infinity of $II$ are  uniquely described by a number of interior masses $p\in \N$ and boundary masses $q\in \N$ with $2p+q=k$ and with respect to which they correspond to the "configurations" $\alpha_i=1$, $\l_i=+\infty$,$\tau_i=0$ $i=1, \cdots, p+q$, $\beta_r=0$, $r=1, \cdots, \bar k$, $A$ is a critical point of $\mathcal{F}_{p, q}$ such that $\mathcal{L}_K (A) < 0 $ and $V=0$ and we denote them by $z^{\infty}$ with $z$ being the corresponding critical point of $\mathcal{F}_{p, q}$ .\\\\
 2) The $II$-energy of a critical point at infinity $z^{\infty}$ denoted by \;$\mathcal{E}_{II}(z^{\infty})$ is given by
 \begin{equation}\label{eq:infinitycriticallevel}
 \mathcal{E}_{II}(z^{\infty})=-\frac{20}{3}k\pi^2-4k\pi^2\ln(\frac{k\pi^2}{6})-8\pi^2\mathcal{F}_{p, q}(z_1, \dots, z_{p+q})
 \end{equation}
 where $z=(z_1, \cdots, z_{p+q})$.
\end{cor}
\vspace{4pt}

\noindent
Furthermore, using Lemma \ref{eq:energyest}, Proposition \ref{eq:conspseudograd}, \eqref{eq:estbarw}, Corollary \ref{eq:c1expansionj1}, and classical Morse lemma, we derive the following Morse type reduction near a  critical point at infinity of $II$.
\begin{lem}\label{eq:morselem}(Morse type reduction  near infinity )\\
Assuming that $(p, q)\in \N^2$  such that $2p+q=k$, $0<\varrho<\varrho_k$, $0<\eta<\eta_k$, and $0<\epsilon\leq \epsilon_k$, where $\varrho_k$ and $\eta_k$ are given by Proposition \ref{eq:sharpest}, and $\epsilon_k$ is given by \eqref{eq:mini} , and $u_0:=\sum_{i=1}^{p+q}\alpha_i^0\varphi_{a_i^0, \l_i^0}+\sum_{r=1}^{\bar k}\beta_r^0(v_r-\ov{(v_r)}_{(Q, T)})+\bar w((\bar\alpha^0, A^0, \bar \l^0, \bar \beta^0))\in \mathcal{V}^{-}_{3}(p, q, \epsilon, \varrho, \eta) $ (where $\bar \alpha^0:=(\alpha^0_1, \cdots, \alpha^0_{p+q})$, $A^0:=(a_1^0, \cdots, a^0_{p+q})$,  $\bar \l:=(\l_1^0, \cdots, \l_{p+q}^0)$ and  $\bar \beta^0:=(a_1^0, \cdots, \beta_{\bar k}^0)$) with the concentration  points $a_i^0$,  the masses $\alpha_i^0$, the concentrating parameters $\l_i^0$ ($i=1, \cdots, p+q$) and the negativity parameters $\beta_r^0$ ($r=1, \cdots, \bar k$) satisfying \eqref{eq:afpara} and furthermore $A^0\in Crit(\mathcal{F}_{p, q})$, then there exists an open neighborhood  $U$ of $(\bar\alpha^0, A^0, \bar \l^0, \bar \beta^0)$ such that for every $(\bar \alpha, A, \bar \l, \bar \beta)\in U$ with $\bar \alpha:=(\alpha_1, \cdots, \alpha_{p+q})$, $A:=(a_1, \cdots, a_{p+q})$,  $\bar \l:=(\l_1, \cdots, \l_{p+q})$, $\bar \beta:=(\beta_1, \cdots, \beta_{\bar k})$, and the $a_i$,  the $\alpha_i$, the $\l_i$ ($i=1, \cdots, p+q$)  and the $\beta_r$ ($r=1, \cdots, \bar k$) satisfying \eqref{eq:afpara}, and $w$ satisfying \eqref{eq:afpara}  with $\sum_{i=1}^{p+q}\alpha_i\varphi_{a_i, \l_i}+\sum_{r=1}^{\bar k}\beta_r(v_r-\ov{(v_r)}_{(Q, T)})+w\in \mathcal{V}^{-}_{3}(p, q, \epsilon, \varrho, \eta)$, we have the existence of a change of variable
\begin{equation}\label{eq:morsevinf}
\begin{split}
&\alpha_i\longrightarrow s_i ,i=1, \cdots, p+q,\\& A\longrightarrow \tilde A=(\tilde A_{-}, \tilde A_{+})\\&\l_1\longrightarrow \theta_1,\\&\tau_i\longrightarrow \theta_i, i=2, \cdots, p+q,\\&\beta_r\longrightarrow \tilde \beta_r\\ &V\longrightarrow \tilde V,
\end{split}
\end{equation}
such that
\begin{equation}
II(\sum_{i=1}^{p+q}\alpha_i\varphi_{a_i, \l_i}+\sum_{r=1}^{\bar k}\beta_r(v_r-\ov{(v_r)}_{(Q, T)}+w)=-|\tilde A_{-}|^2+|\tilde A_{+}|^2+\sum_{i=1}^{p+q}s_i^2-\sum_{r=1}^{\bar k}\tilde \beta_r^2+\theta_1^2-\sum_{i=2}^{p+q}\theta_i^2+||\tilde V||^2
\end{equation}
where $\tilde A=(\tilde A_{-}, \tilde A_{+})$ is the Morse variable of the map \;$\mathcal{E}_{II}: (\mathring{M}^p)^*\times ((\partial M)^q)^*\longrightarrow \R$ which is defined by the right hand side of \eqref{eq:infinitycriticallevel} and $\mathcal{V}_3(p, q, \epsilon)$ is a  neighborhood of ``true`` critical points at infinity of $II$ defined in the proof of Proposition \eqref{eq:conspseudograd} . Hence a critical point at infinity \;$x^{\infty}$ of $II$ has Morse index at infinity  \;$M_{\infty}(x^{\infty})=i_{\infty}(x)+\bar k$.
\end{lem}
\vspace{4pt}



\section{The  boundary-weighted barycenters}\label{s: boundary-barycenters}

This section is devoted to the investigation from algebraic topological viewpoint of the  boundary-weighted barycenters  sets $B_l^{\partial }(M)$, which  are used in this paper to describe the  homotopy type  of  very negative sublevels of the Euler-Lagrange functional $II$ associated to our variational problem. In particular we will compute its Euler characteristic and compute  its Betti numbers.\\
Throughout  section $M$ denotes a compact Riemannian manifold of   dimension $m \geq 2$ with Boundary $\partial M$ and interior $\mathring{M}.$



\subsection{The Euler characteristic of $B_l^\partial (M)$}
Throughout this section and for the sake of simplicity, we normalize the sum  weights to be 1 instead $k$ in the rest of the paper. For the sake of clarity we rewrite the definitions of the basic spaces under this new normalization. That is  we define
\begin{eqnarray*}\label{bar}
B_{p,q}(M, \partial M):= \left  \{ \sum_{i=1}^p \a_i a_i  + \sum_{j=1}^q \beta_j b_j, \alpha_i,\beta_j\in [0,1],\quad a_i\in \mathring{M}, b_j\in\partial M,\ \hbox{ and } \sum_{i=1}^p\alpha_i+\sum_{j=1}^{q}\beta_j=1  \right  \}.
\end{eqnarray*}
and observe that the set $B_{p,q}(M, \partial M)$ is a subspace of the space   $B_{p+q}(M)$ of weighted barycenters of order $p+q$ and is endowed with  the induced topology. \\
Furthermore we  also define for $l \in \N^*$ the space of weighted-barycenters  of order $l$ as
\begin{equation}\label{boundarybar}
B_l^\partial (M) : = \bigcup_{2p+q\leq l}B_{p,q}(M,\partial M).
\end{equation}
We observe that it is a stratified subspace of the barycenter space $B_l(M)$. Moreover  the closure of each stratum $B_{p,q}(M,\partial M)$ in $B_{p+q}(M)$ is contained in
$\bigcup B_{p-i,q+i}(M,\partial M)$ and
we  have inclusions $B_{l-1}^\partial (M)\subset B_l^\partial (M)$.\\
In this section  will often drop writing $M$ from the notation for simplicity.\\
Before going into details in the investigation of the space of weighted-barycenters  we point that,
unlike  the usual barycenter spaces $ B_l(M)$, the subspaces $B_l^\partial (M)$ are not homotopy invariant.
In particular if $M$ is contractible, it is not the case in general
that $B_l^\partial (M)$ is
 also contractible. We will illustrate this with a simple non-trivial example
 which is the closed disk $D^2$ with connected boundary $S^1$. The space
$B_2^\partial (M)$ is made out of two points on boundary or a single point in interior.
This has the topology of the union
$B_2(\partial M)\cup B_1(M)$, with
$B_2(\partial M)\cap B_1(M) = B_1(\partial M)=\partial M$ (see \eqref{pushout1} for a
 schematic description). When $M$ is a disk, $B_2^\partial (D)$ is obtained by
 gluing a disk $D=B_1(D)$ to the circle $S^1=B_1(\partial D)$ sitting inside $B_2(S^1)$. It is well-known
 that $B_2(S^1)\cong S^3$ (\cite{kk}, Corollary 1.4 (b)). Thus
 $B_2^\partial (D)$ is obtained up to homotopy from $S^3$ by collapsing out a circle $S^1$.
 To know what this quotient $S^3/S^1$ is, we need understand the nature of the inclusion
 $B_1(\partial M)=S^1\hookrightarrow S^3=B_2(\partial M)$ (this is by definition a knot).
 We don't know the precise nature of this knot (we suspect it is the trefoil knot, see \cite{mostovoy})
 but since all we need is the homology, we can use Lefshetz duality
 $$\tilde H^*(S^3/S^1)\cong H_{3-*}(S^3\setminus S^1)$$
 The homology of the complement of (any) knot is independent of the embedding and is given by
 $$H_1(S^3\setminus S^1) \cong{\mathbb Z}\ \ \ \ ,\ \ \ \ H_i(S^3\setminus S^1)=0 ,\ i> 1$$
 This calculation is obtained by analyzing the Mayer-Vietoris sequence of the union $S^3=(S^3\setminus S^1)\cup T$ where $T$
 is a tubular neighborhood of $S^1$ homeomorphic to a torus. The upshot is that
 \begin{equation}\label{disk}\tilde H_*(S^3/S^1) = \tilde H_*(B_2^\partial (D))\cong \begin{cases}{\mathbb Z}, &*=3\\
 {\mathbb Z},&*=2\end{cases}
 \end{equation}
 and is zero in all other degrees. We will later check our main theorem against this calculation.

\noindent
The objective of this section is to compute the Euler characteristic of the space $\chi (B_l^\partial )$ for $l \in N^*$.

\begin{thm}\label{main1}
Suppose $M$ is a compact even dimensional manifold with boundary $\partial M$. Then
$$\chi (B_{2l-1}^\partial ) = \chi (B_{l-1}(M))\ \ \ \hbox{and}
\ \ \ \chi (B_{2l}^\partial ) = \chi (B_l(M)).$$
In particular
$$
\chi (B_{l}^\partial (M))= 0  \quad \mbox{ if }  \quad \chi (M)=0.
$$
\end{thm}
We list some useful known facts about the Euler characteristic of familiar constructions:
\begin{itemize}
\item
The Euler characteristic $\chi$ satisfies the inclusion-exclusion principle on closed subsets; i.e. $\chi (A\cup B) = \chi (A)+\chi (B)-\chi (A\cap B)$.
This has a generalization. We say $Y$ is the colimit of the diagram
\begin{equation*}
\xymatrix{B&X\ar[r]\ar[l]&A}
\end{equation*}
if the diagram maps to $Y$ and if $Z$ is any space to which this diagram maps to, there must be an arrow $Y\longrightarrow Z$ factoring this diagram as in the figure
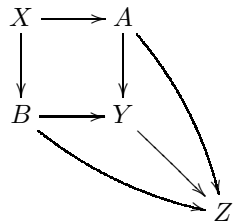
\begin{equation*}
\xymatrix{
X\ar[r]\ar[d]&A\ar[d]\ar@/^/[ddr]\\
B\ar[r]\ar@/_/[drr]& Y\ar[dr]\\
&&Z}
\end{equation*}
i.e. $Y$ is the ``smallest" space to which the diagram can be mapped.
If $Y$ is such a colimit, then
\begin{equation}\label{diag}
\chi (Y) = \chi (A)+\chi (B) -\chi (X)
 \end{equation}
 This equality is true because there is a Mayer-Vietoris sequence associated to this diagram. In particular if all maps in the diagram are inclusions, then $X=A\cap B$, $Y=A\cup B$ and we recover the inclusion-exclusion principle for $\chi$.
\item $\chi (B_l(M))$ has been computed in \cite{kk} and is given by
\begin{equation}\label{chibar}
\chi (B_l(M)) = 1 - {1\over l!}(1 - \chi )\cdots (l-\chi)
\end{equation}
where $\chi := \chi (M)$. In particular, observe that when $M$ is closed odd dimensional, $\chi (M)=0$ and hence necessarily $\chi (B_k(M)) = 0$... This fact will be used throughout.
\item
If $X*Y$ is the join of $X$ and $Y$, then
$$\chi (X*Y) = \chi (X) + \chi (Y) - \chi (X)\chi (Y)$$
A cute way to see this is to notice that $X*Y$ is the colimit of the obvious projection maps in the diagram
$\xymatrix{CX\times Y&X\times Y\ar[l]\ar[r]&X\times CY}$
where $CX$ is the cone on $X$ (a contractible space).
Now apply \eqref{diag}.
\item Let $\vee $ denote the one point union; i.e. $X\vee Y=X\sqcup Y/x_0\simeq y_0$. Then $*$ is (up to homotopy) distributive with respect to $\vee$; i.e. $X*(Y\vee Z)\simeq (X*Y)\vee (X*Z)$.
\end{itemize}
We need a definition

\begin{df} Define $B^p_q(M)$ the subspace of $B_{p+q}(M)$ given by
 $$B_q^p(M) := \{\sum t_ix_i\in B_{p+q}(M)\ \hbox{with at most $p$ of the $x_i$'s in the interior of $M$}\}$$
\end{df}
This is a closed subset of $B_{p+q}(M)$ and since points in the interior are allowed to move into the boundary, we have inclusions $B_{q+i}^{p-i}\subset B_q^p$.
Naturally $B_q^p$ is a subspace of $B_l^\partial$ if $2p+q\leq l$.
For example $B_2^1\subset B_4^\partial$ but $B_2^2$ is not a subspace of $B_4^\partial$. Note also that

\begin{lem}
$B_q^p$ is the closure of $B_{p,q}$ in $B_{p+q}(M)$.
\end{lem}
Before giving a proof of Theorem \ref{main1}, we look closely at the first few cases.
When $k=1$, there isn't much to prove
$$B_1^\partial(M) = B_{0,1}(M,\partial M)= B_1(\partial M)=\partial M$$
If $M$ is even dimensional, then $\partial M$ is closed odd dimensional and its
Euler characteristic must vanish; i.e.
$\chi ( B_1^\partial(M)) = 0$. \\
The case $l=2$ was stated  in \eqref{disk}. We indicated that $B_2^\partial (X)$ is the colimit of
\begin{equation}\label{pushout1}
\xymatrix{ B_1(M)&\partial M\ar[r]\ar[l]&B_2(\partial M)},
\end{equation}
so that
$\chi (B_2^\partial) = \chi (B_2(\partial M)) + \chi (M)-\chi (\partial M)$.
If $M$ is even dimensional, $\chi (\partial M)$ is trivial and so is
$\chi (B_2(\partial M))$ according to \eqref{chibar}. Thus
$\chi (B_2^\partial) = \chi (M)$
in accordance with Theorem \ref{main1}.

 \begin{df}Let $P$ be a poset which we view as a category with morphisms pointing upward. Here $p<q$ in $P$ means $p\rightarrow q$. We will assume that $P$ is a lower semilattice meaning that for any $p,q\in P$ there is a greatest lower bound. A diagram of spaces over $P$ is a functor from $P$ into the category of topological spaces. Given a diagram of spaces, we define the \textit{colimit} of this diagram to be the ``smallest" space to which the diagram can map to (we refer to \cite{strom} for the precise definitions). When all maps in the diagram are inclusions (this is our case), suffices to say that the colimit is constructed by taking the union of all spaces in the diagram glued over common intersections.
\end{df}
It turns out that $B_l^\partial$ is the colimit of a poset diagram.
Let's consider the first few cases. When $l=2$, this is the poset with two edges given
in \eqref{pushout1} which takes the form
$$\xymatrix{B_2^0&&B_0^1\\
&\partial M\ar[ur]\ar[ul]}$$
For $l=3$, $B_3^\partial$ consists of those barycenters with only one point in interior and one point on the boundary, or with $3$ points on boundary and no points in the interior. This is the union $B_3(\partial M)\cup B_1^1$ over the intersection $B_2(\partial M)$, thus it is the colimit of
$$\xymatrix{B_3^0 &B_2^0=B_2(\partial M)\ar[r]\ar[l]& B_1^1}$$
Similarly $B_4^\partial$ is the colimit of the following diagram
$$\xymatrix{B_4(\partial M)&&B_2^1&&B_0^2=B_2(M)\\
&B_3(\partial M)\ar[ul]\ar[ur]&&B_1^1\ar[ur]\ar[ul]\\
&&B_2(\partial M)\ar[ur]\ar[ul]
}$$
This is again interpreted the following way: $B_4^\partial$ is the union
$B_4(\partial M)\cup B_2^1\cup B_2(M)$ with all three subspaces
overlapping according to the arrows in the diagram; i.e. $B_4(\partial M)\cap B_2^1 = B_3(\partial M)$, etc.\\
In general we have

\begin{pro}\label{colimit} The space $B_{2l}^\partial (M)$ is the colimit of the following diagram of spaces
$$
\xymatrix@C=0.5em{B_{2l}^0=B_{2l}(\partial M)&&B_{2l-2}^1&&\cdots&&B_2^{l-1}&&B_0^l=B_l(M)\\
&B_{2l-1}^0\ar[ul]\ar[ur]&&\cdots&&B_3^{l-2}\ar[ur]&&B_1^{l-1}\ar[ur]\ar[ul]\\
&&&&\vdots&&&&\\
&&&B_{l+1}^0\ar[ur]\ar[ul]&&B_{l-1}^1\ar[ur]\ar[ul]\\
&&&&B_l^0\ar[ur]\ar[ul]}$$
Similarly, the diagram whose colimit is $B_{2l-1}^\partial (M)$ is obtained by truncating the top row.
\end{pro}
The proof is self-evident. Since all vertical maps are inclusions, the
colimit of our diagram is obtained by taking the union of all spaces
in the top row whose pairwise intersections are given by the row underneath....
\begin{rem}
We can similarly define the boundary weighted barycenter spaces $B_l^\partial (M,r)$
consisting of points in the interior with integer weight $r\geq 1$, not just $r=2$. This is a modification of \eqref{boundarybar} where now $rp+q\leq l$. Here too we can give a complete description of the space as in Proposition \ref{colimit}.
\end{rem}

\begin{lem} \label{formula}
$\chi (B_{2l}^\partial ) =\sum_{i=0}^l\chi (B_{2l-2i}^i) -
\sum_{i=0}^{l-1}\chi (B_{2l-1-2i}^i)$.
\end{lem}

\begin{pf}
Look at the first top two rows of the colimit diagram in Proposition \ref{colimit}. By the inclusion-exclusion principle $\chi (B_{2l}^0\cup B_{2l-2}^1) = \chi (B_{2l}^0)+\chi (B_{2l-2}^1) - \chi (B_{2l-1}^0)$. The next space in the top row $B_{2l-4}^2$ intersects with this union along the subspace
$B_{2l-3}^1$ (that's the point) so that
\begin{eqnarray*}
\chi (B_{2l}^0)\cup B_{2l-2}^1\cup B_{2l-4}^2) &=& \chi (B_{2l}^0\cup B_{2l-2}^1) + \chi (B_{2l-4}^2) -
\chi (B_{2l-3}^1)\\ &=&
\chi (B_{2l}^0)+\chi (B_{2l-2}^1) + \chi (B_{2l-4}^2) -
\chi (B_{2l-3}^1)- \chi (B_{2l-1}^0)
\end{eqnarray*}
The lemma follows by induction.
\end{pf} \\
We therefore need to compute $\chi (B_q^p)$ for various $p,q$.  When $p=0$ or $q=0$, the computation is obvious:
$\chi (B_q^0) = 0$ if $M$ is even dimensional while
$\chi (B_0^p) = \chi (B_p(M))$.
It remains to determine $\chi (B_q^p)$ for non-zero $q,p$. This is done by induction. The inductive
hypothesis is ensured by the following Lemma.

\begin{lem}\label{bqp} For $q\geq 1$, $B_q^p$ is the colimit of the diagram
$$\xymatrix{B_p(M)*B_q(\partial M)&
B_1^{p-1}*B_q(\partial M)\ar[r]\ar[l]& B_{q+1}^{p-1}}$$
\noindent
In particular when $M$ is even-dimensional with boundary then\
$\chi (B_q^p) = \chi (B_p)$
\end{lem}

\begin{pf} Note that in the diagram above, only the lefthanded map
$B_1^{p-1}*B_q(\partial M)\longrightarrow B_p(M)*B_q(\partial M)$
is an inclusion.
An element in $B_q^p$ consists again of barycenters with at least $q$-points in $\partial M$. Let's now look at a typical element in $B_q^p$ which we write as $\sum t_ix_i$ with at most $p$ of the $x_i$'s in the interior of $M$.
When exactly $p$ such points are in the interior, this is an element of the join
$B_p(M)*B_q(\partial M)$ since it can be written as
$$t \left(\sum_{i=1}^p {t_i\over\sum t_i}a_i\right) + (1-t)\left(\sum_{j=1}^q {s_j\over\sum s_j}b_j\right)\ \ \ \ ,\ \ \ t = \sum t_i=1-\sum s_j\ ,\ a_i\in \dot{M}, b_j\in\partial M$$
On $B_q^p(M)-B_{q+1}^{p-1}$, the map $B_p(M)*B_q(\partial M)\longrightarrow B_q^p(M)$ is one-to-one.
When in $\sum t_ix_i$ one of the $x_i$'s goes to the boundary, $\sum t_ix_i$ approaches an element in $B_{q+1}^{p-1}$. This means that the subspace
$B_1^{p-1}(M)*B_q(\partial M)$ maps to a quotient.
Rephrasing this in terms of colimits,
we obtain the diagram of the proposition. \\
Using the formula for the Euler characteristic of a colimit \eqref{diag} and the formula for $\chi (X*Y)$, we obtain readily that
$$\chi (B_q^p) = \chi (B_p) + \chi (B_{q+1}^{p-1}) - \chi (B_1^{p-1})$$
(we have applied here repeatedly that $\chi (B_q(\partial M))=0$ if $\chi (\partial M))=0$).
We can then proceed by induction. The computation $\chi (B_q^p)=\chi (B_p)$ for all $q\geq 1$ will be immediate if we establish  that $\chi (B_q^1)=\chi (M)$ for all $q\geq 1$. But $B_q^1$ is the colimit of the diagram
$$\xymatrix{
B_{q+1}(\partial M)&\partial M*B_q(\partial M)\ar[r]\ar[l]&M* B_q(\partial M)
}$$
from which we get that $\chi (B_p^1)=\chi (M)$ as desired.
\end{pf} \\
{\bf Proof of Theorem \ref{main1}}: \\
 This is a direct consequence of Lemmas \ref{formula} and \ref{bqp}.


\subsection{The Homology of $B_k^\partial(M)$}

Throughout $\tilde H (X)$ means \textit{reduced} homology; that is for $X$ connected,
$\tilde H_*(X)$ (any coefficients) is $H_*(X)$ if $*>0$ while $\tilde H_0(X)=0$. We also
make the convention that $H_*(X)=0$ if $*<0$. We have the easy observation

\begin{lem}\label{relative} For homological degree $\ast > 0$ and field coefficients
$$H_*(B_l^\partial/ B_{l-1}^\partial)\cong
H_{*-1}(B_{l-1}^\partial)\oplus
H_*(B_l^\partial)
$$
\end{lem}

\begin{pf} The subspace $B_{l-1}^\partial$ is contractible in $B_{l}^\partial$. Indeed
choose $x_0$ some chosen basepoint in $\partial M$ (call it the ``conepoint") and define
\begin{equation}\label{wn}
W_l \, = \, \left \{\sum_1^l \alpha_ix_i\in B_{l}^\partial\setminus B_{l-1}^\partial\ |\ x_i=x_0\ \hbox{for some $i$} \right \}\bigcup B_{l-1}^\partial
\end{equation}
Then $W_l$ is contractible into itself via the contraction
$(t,\sum_1^l \alpha_ix_i)\longmapsto tx_0+\sum (1-t)\alpha_ix_i$. This says that
$B_{l-1}^\partial$ is contractible in $B_l^\partial$ and
the homology long exact sequence for the pair $(B_{l}^\partial, B_{l-1}^\partial)$
breaks down into short exact sequences and so
the claim is immediate with field coefficients.
\end{pf}

We start our homology computation of $H_*(B_l^\partial (M))$ by looking back at the diagram in Proposition \ref{colimit}. This diagram is organized in rows and we let $T_i$ the colimit of the first $i$ rows, $0\leq i\leq l$. We have
$T_0= B_l^0$, $T_1=B_{l+1}^0\cup B_{l-1}^1$, and more generally
$$T_i =  \left \{\sum t_jx_j\in B_{l+i}^\partial ,\ \hbox{where at most $l$ of the $x_j$'s in interior of $M$} \right \}.$$
Obviously $T_k=B_{2l}^\partial$ and $T_{l-1}=B_{2l-1}^\partial$. It is not true however
that $T_{l-2}$ is $B_{2l-2}^\partial$.

As for the case of $B_{l-1}^\partial$ inside $B_l^\partial$, the subspace
$T_i$ is contractible in $T_{i+1}$ so that (with field coefficients)
$$H_*(T_i^\partial/ T_{i-1}^\partial)\cong
H_{*-1}(T_{i-1}^\partial)\oplus
H_*(T_i^\partial).
$$
Consider the corresponding rows in $i$ and $i-1$
$$
\xymatrix@C=0.5em{B_{l+i}^0&&B_{l+i-2}^1&&\cdots&&\cdots&&B_{l-i}^i\\
&B_{l+i-1}^0\ar[ul]\ar[ur]&&B_{l+i-3}^{1}\ar[ul]\ar[ur]&&\cdots&&B_{l-i+1}^{i-1}\ar[ur]\ar[ul]}
$$
The consecutive quotients ${\overline{T}_{i}}:= T_i/T_{i-1}$ for $0\leq i\leq l$ are given according to
\begin{equation}\label{wedge}
\overline{T}_i = {B_{l+i}^0\over B_{l+i-1}^0} \vee
\bigvee_{0\leq j\leq i-1} {B^j_{l+i-2j}\over (B^{j-1}_{l+i-2j+1}\cup B^{j}_{l+i-2j-1})} \vee
{B^i_{l-i}\over B^{i-1}_{l-i+1}}
\end{equation}

\subsubsection{The spaces $\overline{B}_q^p := B_q^p/(B_{q-1}^p\cup B_{q+1}^{p-1})$}
We can think of this quotient as the space of formal barycenters $\sum_{i=0}^{p+q} t_ix_i$ with precisely $p$ points in the interior and $q$ points on the boundary, with the topology that points on the boundary never leave that boundary and if either $t_i\rightarrow 0$ or if a point from the interior approaches the boundary, the whole configuration approaches the basepoint.\\
\noindent
To describe $\overline{B}_q^p$, the following preliminary lemma is needed. We make use of some notation:
\begin{itemize}
 \item $X*\emptyset = X$. An element of $X*Y$ is a segment starting at $X$ when $t=0$ and ending at $Y$ when $t=1$. We view $X$ (resp. $Y$) as a subspace of $X*Y$ corresponding to when $t=0$ (resp. $t=1$).
\item The ``half smash" product of two based spaces is one of
$$X\rtimes Y := {X\times Y\over x_0\times Y}\ \ \ \ \ ,\ \ \ \ X\ltimes Y:= {X\times Y\over X\times y_0}$$
\item The ``unreduced suspension" of a space $\Sigma X$, is the join $X*S^0$ where $S^0 = \{1,-1\}$. Equivalently $\Sigma  X= [0,1]\times X/(0,x)\sim (0,x'), (1,x)\sim (1,x')$. An element of the suspension is written as an equivalence class $[t,x]$.
\end{itemize}

\begin{lem}\label{joinstuff} Let $(X,A)$ and $(Y,B)$ be two connected CW pairs. Then $X*Y/A*Y\simeq X/A* Y$ and
$$(X*Y)/(X*B\cup A*Y)\simeq
\begin{cases} X/A *Y/B &,\  A\neq\emptyset, B\neq\emptyset,\\
X/A \ltimes \Sigma Y &,\  A\neq\emptyset, B=\emptyset,\\
\Sigma (X\times Y)&,\ A=\emptyset, B=\emptyset.
\end{cases}
$$
\end{lem}

\begin{pf}
Here $X*Y/A*Y$ is obtained from $X/A*Y$ by collapsing $x_0*Y$, where $x_0$ is the natural basepoint of quotient $X/A$. Since $x_0*Y$ is a cone hence contractible, collapsing out $x_0*Y$ from $(X/A)*Y$ doesn't change homotopy type.\\
\noindent
On the other hand, let $(X,x_0)$ and $(Y,y_0)$ be pointed spaces.
Consider the subspace $X*y_0\cup x_0*Y$ of the join. This is the union of two cones along the segment $x_0*y_0$. This space is $1$-connected (Van-Kampen) and has trivial homology in positive degrees (Mayer-Vietoris). This means (by Whitehead theorem and since we are working with CW complexes) that this subspace is contractible.
We thus have the equivalence $X*Y\simeq X*Y/X*y_0\cup x_0*Y$. When $A,B$ are non-empty, we therefore have
$$(X*Y)/(X*B\cup A*Y) = \left[(X/A)*(Y/B)\right]/x_0\times (Y/B)\cup (X/A)*y_0
\simeq (X/A)*(Y/B)$$
\noindent
Consider next the case $A=\emptyset, B=\emptyset$ (third case). Elements of the join are classes $[x,t,y]$ with identifications at $t=0$ and $t=1$. When we pass to the quotient $X*Y/X\cup Y$, we obtain the identification space consisting of $[x,t,y]$ with everything collapsed out to $x_0$ (resp. $y_0$) when $t=0$ (resp. when $t=1$). This is by definition the unreduced suspension $\Sigma (X\times Y)$. This can be seen directly from the depiction of the join construction in the figure below.
\begin{figure}[htb]
\begin{center}
\epsfig{file=join.eps,height=1.4in,width=6in,angle=0.0}
\end{center}
\end{figure}
\noindent
The second case is proven similarly but with a little trick. Let's look at the figure again and notice that $X*Y/Y$ is $X\times\Sigma Y$ with a copy of $X$ (in the right of the middle diagram) collapsed out. This is of the homotopy type of $X\ltimes \Sigma Y$; that is
$$X*Y/Y \simeq X\ltimes \Sigma Y$$
from which we deduce the series of equivalences
$X*Y/X\cup A*Y \simeq (X/A)*Y/(X/A)\simeq (X/A)\ltimes\Sigma Y$
and the claim follows.
\end{pf}

\noindent
We turn back to $\overline{B}_q^p$. According to Lemma \ref{bqp} (or by inspection) we have the homeomorphism
\begin{equation}\label{first}
{B_q^p\over B_{q+1}^{p-1}}= {B_p(M) * B_q(\partial M)\over B_1^{p-1} * B_{q}(\partial M)},
\end{equation}
and similarly
\begin{eqnarray*}
\overline{B}_q^p = {B_q^p\over B_{q-1}^p\cup B_{q+1}^{p-1}}&=& {B_p(M) * B_q(\partial M)\over B_1^{p-1} * B_{q}(\partial M)\cup B_p(M) * B_{q-1}(\partial M)},\\
&\simeq& {B_p(M)\over B_1^{p-1}} * {B_q(\partial M)\over B_{q-1}(\partial M)},\\
&\simeq& {B_p(M)\over B_1^{p-1}} * \left({B_q(\partial M)\vee \Sigma B_{q-1}(\partial M)}\right),
\end{eqnarray*}
where for the last equivalence we have used the standard fact that if
 $A$ is contractible in $X$, then $X/A\simeq X\vee\Sigma A$ (and
$B_{q-1}(\partial M)$ is contractible in $B_q(\partial M)$ as already mentioned).
By the distributivity property $X*(Y\vee Z)\simeq (X*Y)\vee (X*Z)$, we can write further
\begin{equation}\label{decompose}
\overline{B}_q^p\simeq
\left({B_p(M)\over B_1^{p-1}} * B_q(\partial M)\right)\vee
\left({B_p(M)\over B_1^{p-1}} *
\Sigma B_{q-1}(\partial M)\right).
\end{equation}
We must therefore understand the quotient
$B_p(M)/ B_1^{p-1}$.

\begin{lem}\label{second} There is a homotopy equivalence
$${B_p(M)\over B_1^{p-1}}\simeq B_{p}(M/\partial M).$$
\end{lem}

\begin{pf}
When $p=1$, $B_1^{p-1} = \partial M$ and $B_p(M)/ B_1^{p-1}= M/\partial M$.
When $p=2$, $B_2(M)$ is the symmetric join $\sj{2}(M) := M*M/\mathfrak S_2$ and
$B_1^1(M) = (\partial M*M)\cup (M*\partial M)/\mathfrak S_2$.
The quotient is
$B_2(M)/B_1^1(M) = \sj{2}(M/\partial M)/A$
where $A$ is the subspace of barycenters $tx+(1-t)y$ with $x=*$ the basepoint;
that is $A$ is up to homotopy
the cone on $(M/\partial M)$ which is contractible. We've just shown that
$B_2(M)/B_1^1(M) \simeq \sj{2}(M/\partial M)$.
The general case is similar. Note that $M/\partial M$ has a preferred basepoint
$x_0$. It is now clear that
$${B_p(M)\over B_1^{p-1}}= {B_{p}(M/\partial M)\over W_p}$$
where $W_p$ is as defined in \eqref{wn} with ``conepoint" $x_0$. This subspace is
contractible so that up to homotopy
$B_{p}(M/\partial M)/W_p\simeq B_p(M)$
and the claim follows.
\end{pf}

\noindent
We can now put everything previous together to get our main calculation. Let $\sigma$ be
the suspension operator acting on homology so that if $x$ is a homology class, then
$\deg (\sigma (x))=\deg (x)+1$. Recall in our notation that $B_0(Y)=\emptyset$ and that
$\emptyset *X=X$.

\begin{thm}\label{main} Let $M$ be a compact manifold with boundary $\partial M$.
The reduced homology of ${\mathcal B}_l^\partial(M)$ with field coefficients is given by
\begin{eqnarray*}
\tilde H_*({\mathcal B}_{2l}^\partial(M))&\cong&
\tilde H_*(B_{2l}(\partial M))\oplus
 \bigoplus_{0<i\leq l}
\tilde H_*\left(B_i({M/\partial M})*B_{2l-2i}(\partial M)\right)\\
&\cong&\tilde H_*(B_{2l}(\partial M))\oplus \tilde H_*(B_l(M/\partial M)\\
&&\oplus\ \sigma \bigoplus_{i=1}^{l-1}\tilde H_*(B_i({M/\partial M}))\otimes
 \tilde H_*(B_{2l-2i}(\partial M))
\end{eqnarray*}
where $\sigma$ is the suspension operator.
Similarly
\begin{eqnarray*}
\tilde H_*({ B}_{2l-1}^\partial(M))&\cong&
\tilde H_*(B_{2l-1}(\partial M))\oplus
 \bigoplus_{0<i\leq l-1}
\tilde H_*\left(B_i({M/\partial M})*B_{2l-2i-1}(\partial M)\right)\\
&\cong&H_*(B_{2l-1}(\partial M))\oplus
 \sigma \bigoplus_{i=1}^{k-1}\tilde H_*(B_i({M/\partial M}))\otimes
 \tilde H_*(B_{2l-2i-1}(\partial M))
\end{eqnarray*}
The formula is valid for $l=1$ by setting $\bigoplus_1^0=0$.
\end{thm}

\begin{pf} We go back to $T_i$ the colimit of the bottom $i$ rows of the
colimit diagram in Proposition \ref{colimit}, for $0\leq i\leq k$. Since $T_{i-1}$ contractible in $T_i$, we have
$${\overline{T}_{i}}:= T_i/T_{i-1} = T_i\vee\Sigma T_{i-1},$$
where $\Sigma T_{i-1}$ is the suspension of $T_{i-1}$. It is known that
$\tilde H_*(\Sigma X)=\sigma \tilde H_*(X)$.
On the other hand there is the wedge decomposition \eqref{wedge}
$$\overline{T}_i = {B_{l+i}^0\over B_{l+i-1}^0} \vee
\bigvee_{0< j<i} {B^j_{l+i-2j}\over (B^{j-1}_{l+i-2j+1}\cup B^{j}_{l+i-2j-1})} \vee
{B^i_{l-i}\over B^{i-1}_{l-i+1}}.
$$
For the same reason as earlier indicated,
$\displaystyle {B_{l+i}^0\over B_{l+i-1}^0}\simeq B_{l+i}(\partial M)\vee\Sigma B_{l+i-1}(\partial M)$, while
$$ {B^i_{l-i}\over B^{i-1}_{l-i+1}}\simeq {B_i(M)\over B_1^{i-1}}*B_{l-i}(\partial M)
\simeq B_i(M/\partial M)*B_{l-i}(\partial M)
$$ according to \eqref{first}, Lemma \ref{joinstuff} and Lemma \ref{second}....
Putting everything together gets us
\begin{eqnarray}\label{decompose}
{\overline{T}_{i}}= T_i\vee\Sigma T_{i-1}
&\simeq&B_{l+i}(\partial M)\vee\Sigma B_{l+i-1}(\partial M)\vee B_i(M/\partial M)*B_{l-i}(\partial M)\nonumber\\
&&\bigvee_{0< j<i} (B_j(M/\partial M)*B_{l+i-2j}(\partial M)) \vee
(B_j(M/\partial M)*\Sigma B_{l+i-2j-1}(\partial M))
\end{eqnarray}
Using induction starting with $B_1^\partial (M)=\partial M$ we can infer
from the decomposition of $\overline{T}_i$ above
that the integral homology (not just with field coefficients) is given by
\begin{eqnarray}\label{homti}
\tilde H_*(T_i)&\cong&
\tilde H_*\left(B_{l+i}(\partial M)\vee \bigvee_{0< j\leq i} (B_j(M/\partial M)*B_{l+i-2j}(\partial M)\right)\\
&\cong&
\tilde H_*(B_{l+i}(\partial M))\oplus
 \bigoplus_{j=1}^{i}
\tilde H_*\left(B_j({M/\partial M})*B_{l+i-2j}(\partial M)\right)
\end{eqnarray}
The argument behind this deduction is to first apply homology to both sides of \eqref{decompose} and to use the cancelation property for finite abelian groups\footnote{For finite complexes, the homology groups are finitely generated abelian groups. For those groups, if $A\times G\cong B\times G$, then $A\cong B$ (cancellation)}
\noindent
We pointed out that for connected based spaces $X*Y\simeq \Sigma X\wedge Y$ so that (with field coefficients now) $\tilde H_*(X*Y)\cong \sigma \tilde H_*(X)\otimes \tilde H_*(Y)$.  The theorem is  immediate replacing $T_l$ by $B_{2l}^\partial$ and $T_{l-1}$ by $B_{2l-1}^\partial$.
\end{pf}

\begin{rem} The case $l=1$ in the theorem gives
$$\tilde H_*({\mathcal B}_{2}^\partial(M))\cong
\tilde H_*(B_{2}(\partial M))\oplus \tilde H_*(M/\partial M)$$
When $M=D$ is the disk with boundary $S^1$,
$$\tilde H_*(B_2^\partial (D))\cong \tilde H_*(B_2(S^1))\oplus \tilde H_*(B_{1}(S^2))
\cong
\tilde H_*(S^3)\oplus \tilde H_*(S^2)\cong
\begin{cases}{\mathbb Z}, &*=3\\
 {\mathbb Z},&*=2\end{cases}
 $$
and this is precisely the computation in  \eqref{disk}.
\end{rem}

\noindent
 One first consequence of Theorem \ref{main} is that we can compute the connectivity and homological
dimension of our spaces $B_l^\partial$. A space $X$ is said to be $r$-connected for some $r\geq 0$ if all the homotopy groups of $X$ vanish up to $r$. A space is $0$-connected if it is connected and it is $1$-connected if it is simply connected.

\begin{lem}\label{funda} Suppose $\partial M$ is $1$-connected, then $B_l^\partial (M)$ is also $1$-connected
\end{lem}

\begin{pf}
Given a diagram as in Proposition \ref{colimit}, it is well-known that if all spaces $B_q^p$ making up the diagram are $1$-connected, then the colimit is $1$-connected (this is a consequence of the Van-Kampen theorem).
When $p=0$, $B_q^0=B_q(\partial M)$ and this is at least $2q-1$ connected according to \cite{kk} so in particular $1$-connected. We can now use induction on $p$ and the colimit diagram in Lemma \ref{bqp} to complete the proof. More precisely recall that if $X$ and $Y$ are connected, then $X*Y$ is simply connected. Since by induction $B_{q+1}^{p-1}$ is connected, it follows that all constituent subspaces in the colimit diagram in Lemma \ref{bqp} are $1$-connected and hence so is $B_q^p$. The proof is complete.
\end{pf}

\begin{lem}
Suppose $M$ is a manifold of $\dim M=m$. Then
 $H_*(\mathcal B_l^\partial)=0$ for $\ast \geq lm$. If in addition both $M$ and $\partial M$ are $r$-connected with $r\geq 1$, then $B_l^\partial (M)$ is $l+r-2$-connected if $l$ is even and $\min \{l+2r-2,2l+r-2\}$-connected if $l\geq 3$ is odd.
\end{lem}

\begin{pf} The first claim is by inspection knowing that the homological dimension of $\mathcal B_n(X)$ is at most $nh(X)+n-1$, where $h(X)$ is the homological dimension of $X$. On the other hand it is shown in \cite{kk} that if $M$ is $r$-connected with $r\geq 1$, then
$\mathcal B_n(M)$ is $2n+r-2$ connected. If both $M$ and $\partial M$ are $r$-connected\footnote{There is in general no relationship between the connectivity of $M$ and $\partial M$.}), then so is $M/\partial M$. This means that
$\tilde H_i({\mathcal B_k}(\partial M))= \tilde H_i({\mathcal B_l}(M/\partial M)) = 0$ if $i\leq 2l+r-2$. We obtain our connectivity bound by computing the homological connectivity of $B_l^\partial (M)$ from Theorem \ref{main} and then invoking Lemma \ref{funda} (\footnote{Here we recall that if $X$ is simply connected with $\tilde H_i(X)=0$ for $i\leq r$, then $\pi_i(X)=0$ for $i\leq r$}). Similar calculation when $l$ is odd.
\end{pf}

\noindent
A second main consequence of Theorem \ref{main} is that we can recover the Euler characteristic computation of $B_l^\partial(M)$ when $M$ even dimensional.\\
\centerline{\bf Alternative Proof of  Theorem \ref{main1}} \\
From \eqref{homti}
\begin{eqnarray*}
\chi (T_i)-1& =&\chi(B_{l+i}(\partial M))-1+
 \sum_{j=1}^{i}(\chi (B_j(M/\partial M)*B_{l+i-2j}(\partial M))-1),\\
 &=&\chi(B_{l+i}(\partial M))-1,\\
 &&+ \sum_{j=1}^{i}\left[\chi (B_j(M/\partial M))
 +\chi (B_{l+i-2j}(\partial M)) - \chi (B_j(M/\partial M))
 \chi (B_{l+i-2j}(\partial M))-1\right].
 \end{eqnarray*}
Suppose $M$ is even-dimensional compact with connected boundary. We can use the fact that $\chi(B_l(\partial M))=0$ since $\chi (\partial M)=0$. From the previous formula
\begin{eqnarray*}
\chi (T_i)&=&-i+ \sum_{j=1}^{i} \chi (B_j (M/\partial M))
\end{eqnarray*}
From which we get
\begin{eqnarray}\label{intermediate}
\chi (B_{2l}^\partial)=\chi (T_k)&=&-l+ \sum_{j=1}^{l} \chi (B_j (M/\partial M)),\\
\chi (B_{2l-1}^\partial)=\chi (T_{l-1})&=&-l-1+ \sum_{j=1}^{l-1} \chi (B_j (M/\partial M)).\nonumber
\end{eqnarray}
According to \eqref{diag} again, $\chi (M/\partial M) = \chi (M)-\chi (\partial M)+1=\chi (M)+1$. Set $\chi:=\chi (M)$. Combining with \eqref{chibar} yields
\begin{eqnarray*}
\chi (B_j(M/\partial M)) &=& 1 - {1\over j!}(-\chi )(1-\chi)\cdots (j-1-\chi),\\
&=&1-{1\over j!}(1-\chi )(1-\chi)\cdots (j-\chi) + {1\over (j-1)!}(1-\chi)\cdots (j-1-\chi),\\
&=&\chi (B_j(M)) - \chi (B_{j-1}(M)) + 1.
\end{eqnarray*}
Adding those up for $0<j\leq l$ and plugging into \eqref{intermediate}
yields immediately our formula
$\chi (B_{2l}^\partial)=\chi (B_l(M))$ and that
$\chi (B_{2l-1}^\partial) = \chi (B_{l-1}(M))$. This is precisely the content
of Theorem \ref{main1}.
\hfill\phantom{}

\subsection{Barycenter spaces of disconnected spaces}

This section which is of independent interest discusses $B_l(X)$ when
$X=A\sqcup B$ is a disjoint union. This case is not treated in \cite{kk} and
it is needed in case our manifold $M$ has disconnected boundary.\\
\noindent
Obviously $B_1(A\sqcup B)= A\sqcup B$ is disconnected. For $n\geq 2$, $B_n(A\sqcup B)$ is
however connected. By convention set $B_0(Y)=\emptyset$ as before.
An element of the join $X*Y$ is written $[x,t,y]$.
As for Lemma \ref{joinstuff}, there are inclusions
$X\subset X*Y$ and $Y\subset X*Y$ where $X=\{[x,0,y]\}$ and $Y=\{[x,1,y]\}$. We will
write these subspaces as $\emptyset* Y$ and $X*\emptyset$.

\begin{pro}\label{main3} $B_n(A\sqcup B)$ is the colimit of the lattice diagram
$$
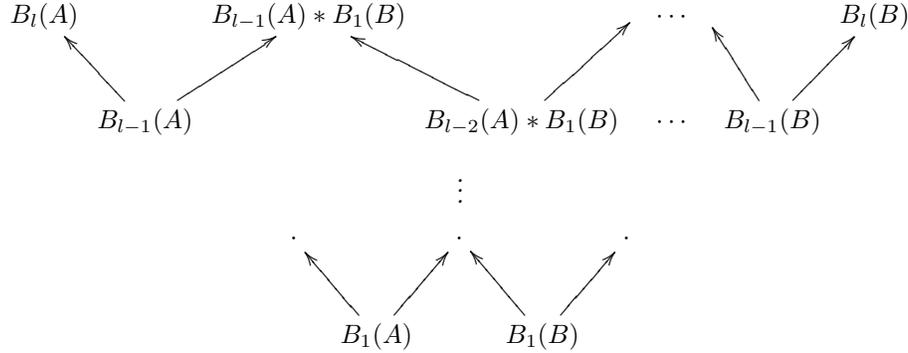
\xymatrix@C=0.1em{B_l(A)&&B_{l-1}(A)*B_{1}(B)&&&\cdots&&&B_l(B)\\
&B_{l-1}(A)\ar[ul]\ar[ur]&&B_{l-2}(A)*B_{1}(B)\ar[ul]\ar[ur]&&\cdots&&
B_{l-1}(B)\ar[ul]\ar[ur]}$$
$$\vdots$$
$$\xymatrix@C=1em{\cdot&&\cdot&&\cdot\\
&B_{1}(A)\ar[ur]\ar[ul]&&B_1(B)\ar[ur]\ar[ul]&}$$
All arrows are inclusions.
\end{pro}

\noindent
This diagram expresses the fact that points can be either in $A$ or $B$. Points in $A$ alone or $B$ alone ``interact" with each other, while points in $A$ do not interact with points in $B$.
Proposition \ref{main3} gives a homotopy decomposition which we will use to get to the homology of $B_n(A\sqcup B)$(\footnote{In the algebraic topological jargon, a homotopical decomposition of a space $X$ is a diagram of spaces whose  ``homotopy colimit" is weakly equivalent to $X$. When the diagram is indexed over a poset and all maps are ``nice inclusions" (i.e. cofibrations), the homotopy colimit and the colimit are equivalent. This is our case.}).\\

\noindent
Next and throughout, $\overline B_l(X):= B_l(X)/B_{l-1}(X)$. Note that because $B_{l-1}(X)$ is contractible in $B_l(X)$, $\ob{l}(X)\simeq B_l(X)\vee\Sigma B_{l-1}(X)$. This is a fact we will use routinely below.

\begin{lem}\label{ob} Let $X=A\sqcup B$ and $l\geq 3$. Then
\begin{eqnarray*}
\ob{l}(A\sqcup B)&\simeq&
\ob{k}(A)\ \vee\ \ob{l-1}(A)\ltimes\Sigma B\ \vee\
\ob{l-2}(A)*\ob{2}(B)\ \vee\ \nonumber\\
&&\cdots\ \vee\ \ob{2}(A)*\ob{l-2}(B)\ \vee\
\Sigma A\rtimes\ob{l-1}B\ \vee\ \ob{l}(B).
\end{eqnarray*}
\end{lem}

\begin{pf}
In the diagram of Proposition \ref{main3}, the colimit of the first $l$ rows is $B_l(A\sqcup B)$ while the colimit of the first $l-1$ rows is $B_{l-1}(A\sqcup B)$. By collapsing out this latter subspace we obtain the wedge
\begin{equation}\label{bbark}
\ob{l}(A\sqcup B) = \ob{l}(A)\vee
\bigvee_{i+j=l\atop 1\leq i\leq l-1} {B_i(A)*B_j(B)\over B_{i-1}(A)*B_j(B)\cup B_i(A)*B_{j-1}(B)}\vee
\ob{l}(B),
\end{equation}
where we convene that $B_0=\emptyset$.
But
\begin{eqnarray*}
{B_i(A)*B_j(B)\over B_{i-1}(A)*B_j(B)\cup B_i(A)*B_{j-1}(B)}&\simeq&
\overline{B}_i(A)*\overline{B}_j(B)\ \ \ \ \ \ \ \ \ \ \ \ \ \ \ \ \ \ (\hbox{Lemma \ref{joinstuff}, $i,j>1$}),
\end{eqnarray*}
while for $j=1$ (similarly if $i=1$) we have
\begin{eqnarray*}
{B_{l-1}(A)*B\over B_{l-1}(A)\cup B_{l-2}(A)*B}&\simeq&
\overline{B}_{l-1}\ltimes \Sigma B.
\end{eqnarray*}
Combining with \eqref{bbark} gets us the answer.
\end{pf}

\noindent
All homology below is reduced and is taken with field coefficients.

\begin{lem}\label{l=2}
In the case $l=2$,
$B_2(A\sqcup B)$ has the homology of
$$B_2(A)\vee \Sigma (A\times B)\vee B_2(B).$$
\end{lem}

\begin{pf}
The corresponding colimit diagram for $B_2(A\sqcup B)$ has just two rows
$$
\xymatrix@C=0.5em{B_2(A)&&B_{1}(A)*B_{1}(B)&&B_2(B)\\
&B_{1}(A)\ar[ul]\ar[ur]&&B_{1}(B)\ar[ul]\ar[ur]
}$$
so that as we've already seen
\begin{eqnarray*}
\ob{2}(A\sqcup B)&\simeq& \ob{2}(A)\ \vee\
\Sigma (A\times B)\ \vee\ \ob{2}(B)\\
&\simeq& B_2(A)\vee
\Sigma A\ \vee\
\Sigma (A\times B)\ \vee\ B_2(B)\vee
\Sigma B.
\end{eqnarray*}
But the left hand side has another description. Since $A$ and $B$ are disjoint in $B_2(A\sqcup B)$, take a path there $ta+(1-t)b, t\in I$ with $a\in A$ and $b\in B$.
The union $I\cup A\cup B$ is connected and contractible in $B_2(A\sqcup B)$.... It follows
that $\ob{2}(A\sqcup B)\simeq B_2(A\sqcup B)\vee\Sigma (A\vee B)$. Comparing with the above,
we have the equivalence
$$B_2(A\sqcup B)\vee\Sigma A\vee\Sigma B\simeq
B_2(A) \vee \Sigma (A\times B) \vee B_2(B)\vee
\Sigma A\vee
\Sigma B.$$
In homotopy theory, cancelation of wedge summands doesn't hold in general. However passing
to homology with field coefficients, we will obtain direct summands that we can cancel. This is saying precisely that $H_*(B_2(A\sqcup B))\cong H_*(B_2(A)\vee
\Sigma (A\times B)\vee B_2(B) )$
as desired.
\end{pf}

\noindent
We are now in a position to state the main result of this section.

\begin{thm} For $A,B$ two disjoint connected spaces and $l\geq 2$,
$B_k(A\sqcup B)$ has the same reduced homology as
\begin{eqnarray*}
&&B_l(A)\vee \Sigma B_{l-1}(A)\vee B_l(B)\vee \Sigma B_{l-1}B\\
&&\vee\bigvee_{i=1}^{l-1} B_{l-i}(A)*B_i(B)\ \vee\ \bigvee_{i=2}^{l-1}
\Sigma B_{l-i}(A)*B_{i-1}(B)
\end{eqnarray*}
(if the lower index of the wedge is greater than the upper index, which happens when $l=2,3$, then the corresponding wedge term doesn't exist).
\end{thm}

\begin{pf} This is an inductive argument based on the decomposition of Lemma \ref{ob} and the fact that $\ob{l}(X)\simeq B_l(X)\vee\Sigma B_{l-1}(X)$. Note that the formula agrees with the decomposition in Lemma \ref{l=2} since $\Sigma (X\times Y)\simeq X*Y\vee \Sigma X\vee\Sigma Y$ sot that
$$H_*(B_2(A\sqcup B))\cong H_*(B_2(A))\oplus H_*(B_2(B))\oplus H_*(\Sigma A)\oplus H_*(\Sigma B)\oplus H_*(A*B)$$
which is precisely the statement of the theorem when $l=2$.\\
To proceed further we use the same idea as in Lemma \ref{l=2}. We need use here that
$$H_*(X\ltimes \Sigma Y) \cong H_*((X\wedge \Sigma Y)\vee \Sigma Y)
\cong H_*(X*Y)\oplus H_*(\Sigma Y)
$$
According to the decomposition in Lemma \ref{ob}, $H_*(\ob{l}(A\sqcup B))$ is therefore given by
\begin{eqnarray}\label{obk}
&&H_*(\ob{l}(A\sqcup B)\nonumber\\
&&\cong H_*(\ob{l}(A))\oplus H_*(\ob{l-1}(A)* B)\oplus H_*(\Sigma B)\oplus\nonumber\\
&&\bigoplus_{i=2}^{l-2} H_*(\ob{l-i}(A)* \ob{i}(B))
\oplus H_*(A*\ob{l-1}B)\oplus H_*(\Sigma A)\oplus H_*(\ob{l}(B))\nonumber\\
&&\cong H_*(\ob{l}(A))\oplus H_*(\ob{l-1}(A))\oplus
\bigoplus_{i=1}^{l-1} H_*(\ob{l-i}(A)* \ob{i}(B))
\oplus H_*(\ob{l-1}(B))\oplus H_*(\ob{l}(B))
\end{eqnarray}
Here we should recall that
\begin{eqnarray*}
\ob{i}(A)*\ob{j}(B)&\simeq&(B_i(A)\vee\Sigma B_{i-1}(A))*(B_j(B)*\Sigma B_{j-1}(B))\\
&\simeq& B_i(A)*B_j(B)\vee B_i(A)*\Sigma B_{j-1}(B)\vee
\Sigma B_{i-1}(A)*B_j(B)\\
&&\vee \Sigma B_{i-1}(A)*\Sigma B_{j-1}(B)
\ \ \ \ \ \ (\hbox{distributivity property})
\end{eqnarray*}
We can combine this with \eqref{obk} to obtain (using reduced homology)
\begin{eqnarray*}
\bigoplus_{i=1}^{l-1} H_*(\ob{l-i}(A)* \ob{i}(B))
&\cong&H_*(B_{l-1}(A)* B)\oplus H_*(\Sigma B_{l-2}(A)*B) \\
&&\bigoplus_{i=2}^{l-2} H_*(B_{l-i}(A)* B_i(B))\oplus
H_*(\Sigma B_{l-i-1}(A)* B_i(B))\\
&&\bigoplus_{i=2}^{l-2} H_*(\Sigma B_{l-i}(A)* B_{i-1}(B))\oplus
H_*(\Sigma^2 B_{l-i-1}(A)* B_{i-1}(B))\\
&&\oplus H_*( A* B_{l-1}(B))\oplus H_*(\Sigma A* B_{l-2}(B))
\end{eqnarray*}
If we now suppose that $H_*(B_{l-1}(A\sqcup B))$ has the homology stated in the statement
of the theorem, then by identifying the expression for $H_*(\ob{l}(A\sqcup B))$ in
\eqref{obk} with $H_*(B_{l}(A\sqcup B))\oplus
H_*(\Sigma B_{l-1}(A\sqcup B))$ and using induction, we end up with the expression of the homology $H_*(B_l(A\sqcup B))$ stated in the theorem.
\end{pf}\\

\noindent
We look at particular cases below.


\begin{rem} Let $A=S^1 = B$ so that $X$ is the disjoint union of two circles. Then
$B_3(S^1)\simeq S^5$, $B_2(S^1)\simeq S^3$ and $B_3(S^1\sqcup S^1)$ has the homology of
\begin{eqnarray*}
&&B_3(A)\vee\Sigma B_2(A)\vee B_2(A)*B\vee \Sigma A*B\vee A*B_2(B)\vee B_3(B)\vee\Sigma B_2(B)\\
&&\simeq\ \bigvee^4S^5\vee\bigvee^3S^4
\end{eqnarray*}
In general if $X$ is a disjoint union of circles, then $B_n(X)$ will be a wedge of spheres, the top dimensional spheres being of dimension $2n-1$.
\end{rem}


\section{Proof of the main results}\label{s:proofm}

In this section, we present the proof of our main  results. The first type of these are based on our strong Morse inequalities and the related Hopf-Poincrar\'e index type formula, while second type of results takes advantage for the precise knowledge of  lack of compactness, to derive that, under appropriate conditions of the function $K$ blow up of subcritical respectively supercritical approximation can rule out.

\subsection{Proof of  Theorem \ref{eq:morsepoincare1} to Theorem \ref{eq:Cm}}
The basic idea of the proof is  to extend the classical Morse theory to our problem. To do so we first characterize the topology of the very high and the very negative ones of $II$. We treat first the very high sublevels of $II$ and for that we start with the following lemma.
\begin{lem}\label{eq:energybddinf}
Assuming that $(p, q)\in \N^2$  such that $2p+q=k$, $0<\varrho<\varrho_k$, $0<\eta<\eta_k$, and $0<\epsilon\leq \epsilon_k$, where $\varrho_k$ and $\eta_k$ are given by Proposition \ref{eq:sharpest}, and $\epsilon_k$ is given by \eqref{eq:mini},  then there exists $C_0^k:=C^k_0(\varrho, \eta)$ such that  for every  $0<\epsilon\leq \epsilon_k$ where $\epsilon_k$ is as in \eqref{eq:mini}, there holds
\begin{equation*}
V(p, q,\epsilon, \varrho,  \eta)\subset II^{C_0^k}\setminus II^{-C_0^k}.
\end{equation*}
\end{lem}
\begin{pf}
It follows directly from \eqref{eq:para}-\eqref{eq:afpara}, Lemma \ref{eq:energyest}, Lemma \ref{eq:gradientlambdaest}, and Proposition \ref{eq:expansionJ1}.
\end{pf}

\vspace{4pt}
\noindent
Next, combining Proposition \ref{eq:escape} and the latter Lemma, we have the following Corollary.
\begin{cor}\label{eq:energybddcrit}
There exists a large positive constant \;$C_1^k$\; such that
\begin{equation*}
Crit(II)\subset II^{C_1^k}\setminus II^{-C_1^k}.
\end{equation*}
\end{cor}
\begin{pf}
It follows, via a contradiction argument, from the the fact that $II$ is invariant by translation by constants, Proposition \ref{eq:escape}, and  Lemma \ref{eq:energybddinf}.
\end{pf}
\vspace{6pt}

\noindent
Now, we are ready to characterize the topology of very high sublevels of $II$. Indeed, as in \cite{no1} and for the same reasons, we have that Lemma \ref{eq:deformlemr}, Lemma \ref{eq:energybddinf} and Corollary \ref{eq:energybddcrit} imply the following one which describes the topology of very high sublevels of the Euler-Lagrange functional $II$.
\begin{lem}\label{eq:tophigh}
Assuming that $(p, q)\in \N^2$  such that $2p+q=k$, $0<\varrho<\varrho_k$, and $0<\eta<\eta_k$, where $\varrho_k$ and $\eta_k$ are given by Proposition \ref{eq:sharpest}, then there exists a large positive constant  $L^k:=L^k(\varrho, \eta)$ with \;$L^k>2\max\{C_0^k, C_1^k\}$\;  such that for every \;$L\geq L^k$, we have  that \;$II^L$\; is a deformation retract of $\mathcal{H}_{\frac{\partial}{\partial n}}$, and hence it has the homology of a point, where $C^k_0$ is as in Lemma \ref{eq:energybddinf} and $C^k_1$ as in Lemma \ref{eq:energybddcrit}.
\end{lem}
\vspace{4pt}

\noindent
Next, we turn to the study of the topology of very negative sublevels of $II$ when $k\geq 2$ or $\bar k\geq 1$. Indeed, as in \cite{no1} and for the same reasons, we have that the well-know topology of very negative sublevels in the {\em nonresont} case (see  \cite{nd2}), Lemma \ref{eq:escape}, Lemma \ref{eq:energybddinf} and Corollary \ref{eq:energybddcrit} imply the following Lemma which gives the homotopy type of the very negative sublevels of the Euler-Lagrange functional $II$.
\begin{lem}\label{eq:topnegative}
Assuming that $k\geq 2$ or $\bar k\geq 1$, $(p, q)\in \N^2$  such that $2p+q=k$, $0<\varrho<\varrho_k$, and  $0<\eta<\eta_k$, where $\varrho_k$ and $\eta_k$ are given by Proposition \ref{eq:sharpest}, then there exists a large positive constant $L_{k, \bar k}:=L_{k, \bar k}(\varrho, \eta)$ with  $L_{k, \bar k}>2\max\{C_0^k, C_1^k\}$ such that for every \;$L\geq L_{k, \bar k}$, we have  that \;$II^{-L}\;$  has the same homotopy type as  $B_{k-1}^{\partial}(M)$ if $k\geq 2$ and $\bar k=0$,  as\;$A_{k-1, \bar k}^{\partial}$ if $k\geq 2$ and $\bar k\geq 1$ and  as $S^{\bar k-1}$ if $k=1$ and $\bar k\geq 1$, where $C^k_0$ is as in Lemma \ref{eq:energybddinf} and $C^k_1$ as in Lemma \ref{eq:energybddcrit} .
\end{lem}
\vspace{4pt}

\noindent
Now, we  present the proof of Theorem \ref{eq:morsepoincare1}-Theorem \ref{eq:Cm}.\\\\
\begin{pfn}{ of Theorem \ref{eq:morsepoincare1}-Theorem \ref{eq:Cm}}\\
We will provide only the proof of Theorem \ref{eq:morsepoincare2}, its corollary and Theorem \ref{eq:Cm}. The proofs of the remaining statements concerning the critical case $\kappa = 4 \pi^2$ are similar and even simpler since they involve only single boundary blow up points.\\
For the sake of simplicity we assume that the Paneitz operator  $P_g^{4,3}$is non negative. The argument extends virtually to the general case, see \cite{no1}. \\
Arguing by contradiction we assume that the the functional $II$ does not admit a critical point. Thanks to Lemma \eqref{eq:energybddinf}, we may choose $L$ large enough such that all critical point at Infinity are contained in the strip $(II^{L}, II^{-L}).$ \\
Now we  define
$$
M(t):= \sum_{i=1}^{4k -1} m_i^k t^i; \quad \quad P(t):= \sum_{i=0}^{\infty} b_i(II^L, II^{-L}) t^{i},
$$
where
$$
b_i:= rank \,  H_i(II^{L}, II^{-L}).
$$
We notice that it follows from the exact sequence of the pair $(II^L, II^{-L})$ that
$$
H_0(II^{L}, II^{-L}) \, \simeq \, H_1(II^{L}, II^{-L})  \mbox{ and } H_i(II^{L}, II^{-L}) \, \simeq  \, H_{i-2} (II^{-L}), \, \forall i \geq 2.
$$
Hence it follows that Lemma \ref{eq:topnegative} that
$$
P(t)\, = \,  \sum_{i=2}^{4 k -5} c_{i-1}^{k-1} \,  t^{i}.
$$
Moreover it follows from   our  Morse Lemma \ref{eq:morselem} that  strong  Morse inequalities hold. Namely we have that:
\begin{equation}\label{eq:mineq}
 M(t) \, - \, P(t) \, = \, (1+t) R(t),
\end{equation}
 where $R(t) :=  \sum_{i \geq 1} n_i t^{i}$ is a polynomial in $t$ with non negative integer coefficients. \\
 Equating the coefficients of $t$ in the polynomials on the left and right hand side of \eqref{eq:mineq} we obtain a solution of the system \eqref{eq:mp3} hence contradicting the assumption of Theorem \ref{eq:morsepoincare2}.\\
 Now choosing  $t= -1$ in equation \eqref{eq:mineq} we derive that:
 $$
 \sum_{A \in \mathcal{F}_{\infty}} (-1)^{i_{\infty}(A)} \, = \, \chi(II^L, II^{-L}) \, = \, 1 - \chi(B^{\partial}_{k-1}),
 $$
 which violates the condition of corollary \eqref{eq:existence2}.\\
 The proof of Theorem \ref{eq:Cm} follows for similar arguments using the  Morse subcomplex related to the critical points at Infinity whose Morse indices are less or equal some $ l \leq 4 k -1.$
\end{pfn}
\vspace{4pt}

\subsection{Proof of Theorem \ref{eq:algtopg}}
Our next theorem  is based on the construction of a solution of supercritical   non resonant approximation. Such a solution is built using the top homology class of  the boundary-weighted barycenters $B^{\partial }_{k}(M)$.  For such a solution we derive an accurate estimate of its  Morse index. We then use such a spectral information to rule out its blow up, proving that it should converge to a solution of our equation.

\begin{pfn}{ of Theorem  \ref{eq:algtopg} }\\
We consider the following superapproximation of the resonant case
\begin{equation}\label{eq:Pepsilon}
(\mathcal{P}_{\e})
\begin{cases}
 P^4_g u \, + \, 2 Q_g  \, = \, 2(\kappa + \e ) \frac{K e^{4 u}}{\int_M K e^{4 u}}  \mbox{ in } \mathring{M}, \\
P_g^3 u \, +  \, T_g \,  \,  = \, \, 0  \quad \mbox{ on }  \partial M,\\
\frac{\partial u}{ \partial n_g} \, \, = \, \, 0 \quad  \mbox{ on } \partial M,
\end{cases}
\end{equation}
where $ \e$ is a small positive number and $ \kappa:= \kappa_{(P^4,P^3)} = 4 \pi^2 k,$ where $k \in N^*.$ \\
Regarding the problem $(\mathcal{P}_{\e})$ we prove the following claim: \\
{ \textbf{Claim 1}:  \it For a sequence of $\e_k \to 0$ , the problem $(\mathcal{P}_e)$ admits a solution $u_{\e}$ whose Morse index $Morse(u_e)$ satisfies
$$
4 k  \, + \ov{k} \leq  \, Morse(u_{\e}) \,  \leq   \, 4 k + 1 + \ov{k}.
$$ }
For the sake of simplicity of notation we provide the proof only in the case $\ov{k} = 0$. The arguments extend virtually to the  case $\ov{k} \neq 0.$ \\
To prove the above claim we argue as follows: 
We embedded the $B_k(\partial M)$ into the space of variation $\mathcal{H}_{\frac{\partial}{\partial n}}$ through the map:
$$f_k(\lambda) : B_k(\partial M)\longrightarrow \mathcal{H}_{\frac{\partial}{\partial n}}$$ as follows
\begin{equation}\label{eq:fplambda}
f_k(\lambda)(\sum_{i=1}^k \alpha_i\delta_{a_i}):=\sum_{i=1}^k \alpha_i \varphi_{a_i, \lambda},
\end{equation}
with the $\varphi_{a_i, \l}$'s defined by \eqref{eq:projbubble}.\\
Now  we notice that it follows  from Lemma \ref{l:A1}, Lemma \ref{l:A2} and Lemma \ref{eq:topnegative} that  $f_k(\l)$ maps for $\l$ and $L$ large
$$
f_k(\l): B_k(\partial M) \to II^{-L}_{\e},
$$
where $II_{ \e}$ is the Euler-Lagrange functional associated to  $(P_{\e}).$
 Hence $\mathcal{M}_k(\l):= f_k(\l)(B_k(\partial M))$ is a stratified set of top dimension $4k-1$.
 Now observe that $\mathcal{M}_k(\l)$ is contractible in $\mathcal{H}_{ \frac{\partial}{\partial n}}$ (by taking its suspension for example), hence the image of such a contraction $\mathcal{U}_k$ is a stratified set of top dimension $4k$. We now deform $\mathcal{U}_k$ using the pseudogradient flow obtained as a convex combination of the Bahri-Lucia pseudogradient of Lemma \ref{eq:deformlem} and the pseudogradient constructed in Proposition \ref{eq:conspseudograd}.\\
By transversality arguments, we may assume that such a deformation avoids the stable manifolds of critical points of $II_{\e}$ whose indices are grater or equal $4 k + 2$.\\
Now using the compactness of the variational problem in the non resonant case and a theorem of Bahri-Rabinowitz \cite{BR}, we derive that $\mathcal{U}_k(\l)$ retracts by deformation onto $II_{\e}^{-L} \cup \Sigma$, where $\Sigma$ is the union of the unstable manifolds of some critical points of $II_{e}$ caught by the flow. Since $\mathcal{U}_k(\l)$ is contractible whereas $II_{\e}^{-L}$ is not, $\Sigma$ is not empty. Moreover from the above transversality arguments, the Morse indices of its critical points are upper bounded by $4k +1$. Hence
$$
\Sigma = \bigcup  \big (W_u(x);  x \mbox{ is a critical point of } II_{e} \mbox{ whose Morse index } m(x) \leq 4 k + 1 \big ).
$$
Now using the exact homology sequence of the pair $(II_{\e}^{-L} \cup \Sigma)$ we derive that
$$\begin{CD}
\cdots@>>> H_{4k}(II_{\e}^{-L} \cup \Sigma) @>>> H_{4k}(II_{\e}^{-L} \cup \Sigma, II_{\e}^{-L}) @>>> H_{4k -1 }(II_{\e}^{-L})\\ @>>> H_{4k-1}(II_{\e}^{-L} \cup \Sigma)@>>>\cdots
\end{CD}$$
Since $II_{\e}^{-L} \cup \Sigma$ is retract by deformation of $\mathcal{U}_k$ which is contractible, we derive that
$$
H_{4k}(II_{\e}^{-L} \cup \Sigma, II_{\e}^{-L}) \, = H_{4k-1}(II_{\e}^{-L}) \neq 0.
$$
It follows that $\Sigma$ contains at least a critical point of $II_{e}$ whose Morse index is either $4k$ or $4 k + 1$.\\
To conclude the proof of the theorem we prove the following claim: \\
 \textbf{Claim 2:}
$$
u_{\e_k}  \to u_{\infty} \mbox{ in } C^{4,\a}(M),
$$
where $u_{\infty}$ is a solution  of equation \eqref{eq:bvps}.  \\
To prove the claim it is enough to rule out the blow up of $u_{\e_k}$.  Arguing by contradiction we assume that
$$
|| u_{\e_k} - \ov{(u_{e_k})}_{Q,T} || \to + \infty.
$$
Arguing as in Lemma \ref{eq:escape} we derive that
$$
u_{\e_k} \in V(p,q,\d_k, \varrho, \eta),
$$
for some $ \d_k \to 0.$ \\
Testing the equation $\mathcal{P}_\e$ by $  \sum_{i=1}^{p+q}  \frac{\l_i}{\a_i} \frac{\partial \varphi_{a_i,\l_i}}{\partial \l_i}$, very much like in Corollary  \ref{eq:cgradientlambdaest}, we derive that:
$$
\e_k  \, = \, 2 \pi^2 \sum_{i=p+1}^{p+q}   \frac{ \frac{\partial \mathcal{F}^A_{i}}{\partial n_{g_i}}(a_i)  }{\l_i \mathcal{F}^A_i(a_i)} \, + \left (   \e_k^2 \, + \sum_k \frac{1}{\l_i^2} + |\a_i -1|^2 \, + |\tau_i|^2  \right).
$$
Now  expanding  $II_{\e_k}(u_{\e_k})$, just like in Lemma \ref{eq:energyest} we obtain:
\begin{equation*}
\begin{split}
II_{\e_k}( u_{\e_k}) \,  = \,   -\frac{20}{3}k\pi^2 \, - \, 4k\pi^2 (1 + \e_k) \, \ln (\frac{k\pi^2}{6}) \,  - \, 8\pi^2( 1 + \e_k ) \, \mathcal{F}_{p, q}(a_1, \dots, a_{p+q})   \\ 16 \pi^2 \left ( 2\sum_{i=1}^p(\alpha_i-1)^2 \ln \l_i  \, +\sum_{i=p+1}^{p+q}(\alpha_i-1)^2 \ln\l_i \right ) \,  -  \, 4 \pi^2 \left [ \sum_{i=1}^p 2 \sigma_i^2 \, +\, \sum_{i= p+1}^{p+q}  \sigma_i^2     \right ] \\
 -  2 \pi^2 \sum_{i=p+1}^{p+q}\frac{1}{\l_i} \, \frac{1}{\mathcal{F}_i^A(a_i)}\frac{\partial \mathcal{F}_i^A}{\partial n_{g_{a_i}}}(a_i)   + O\left(   \e_k^2  \, +  \, \sum_{k=1}^{p+q}|\alpha_k-1|^3+|\sigma_k|^3 + \frac{1}{\l_k^2}\right),
\end{split}
\end{equation*}
Denoting by $A$ the critical point of $\mathcal{F}_{p,q}$ to which the concentration points converge, we derive from the above expansion that

\begin{equation}
Morse(u_{\e_k}) \, = \
\begin{cases}
5 p + 4 q - Morse(\mathcal{F}_{p,q}, A) \, - 1 \, \mbox{ if }  \,  \mathcal{L}_K(A)  < 0, \\
5 p + 4 q - Morse(\mathcal{F}_{p,q}, A) \,  \, \mbox{ if }  \, \mathcal{L}_K(A)  > 0.
\end{cases}
\end{equation}
Therefore we see that the Morse index of such a blowing solution is superabound by $4 k$, where $k = 2p + q$. Taking into account the Morse estimates of  the Claim $1$, we derive that a necessary condition for $u_{\e_k}$ to blow up  is that:
$$
Morse(u_{e}) = 4 k, \quad p = 0 \quad q = k
$$
and  $A$ is local minimum of  $\mathcal{F}_{0,k}$  satisfying that $ \mathcal{L}_K(A) > 0.$ Hence we reach a contradiction with the assumption of Theorem \ref{eq:algtopg}. Therefore the theorem is fully proven.
\end{pfn}


\begin{pfn}{ of Theorem  \ref{th:thmodd1} }\\
We consider the following subcritical approximation
\begin{equation}\label{eq:Pepsilon2}
(\mathcal{P}^{\e})
\begin{cases}
 P^4_g u \, + \, 2 Q_g  \, = \, 2(\kappa - \e ) \frac{K e^{4 u}}{\int_M K e^{4 u}}  \mbox{ in } M, \\
P_g^3 u \, +  \, T_g \,  \,  = \, \, 0 \mbox{ on }  \partial M,\\
\frac{\partial u}{ \partial n_g} \, \, = \, \, 0 \mbox{ on } \partial M,
\end{cases}
\end{equation}
and denote by $II^{\e}$ its associated Euler-Lagrange functional.\\
Thanks to Moser-Trudinger inequality the functional $II^{\e}$ achieves its minimum, say $u^{\e}.$ \\
We claim that $u^{\e}$ converges to $u^0$ a  critical point of $II$.  Since otherwise it would blow up and hence due to Lemma \ref{eq:escape} and corollary \ref{eq:loccritinf} that $u^{\e}$ has to concentrate at a maximum  point $a \in \partial M$ of $K$ such that $$ \frac{\partial K}{\partial n_g}(a) < 0 $$ which contradicts the assumption of the theorem. Hence Theorem \ref{th:thmodd1} is fully proven.
\end{pfn}



\section{Appendix}\label{eq:appendix}

In this section we collect various estimates of the projected bubble $\varphi_{a,\l}$ and its derivative as well as useful estimates of the functional $II$ and its gradient. These estiamtes are standard and use elementary arguments. Hence we to keep the paper  at reasonable length we omit their proofs.
\subsection{Bubble estimates}

Using the conformal invariance of the Paneitz operator, Chang-Qing operator and the conformal Neuman operator  (recalling $H_g=0$), the properties of the metric $g_{a}$ (see \eqref{eq:detga}-\eqref{eq:proua}), the BVP satisfied by the $\varphi_{a, \l}$'s (see \eqref{eq:projbubble} and \eqref{eq:projbubble1}), and the Green's representation formula \eqref{eq:G4integral}, we derive the following two Lemmas.

\begin{lem}\label{eq:intbubbleest}
Assuming that $\epsilon$ is positive and small, $0<\varrho<\varrho_k$ where $\varrho_k$ is as in Proposition \ref{eq:sharpest}, and $\l\geq \frac{1}{\epsilon}$, then\\
1)\;If $d_g(a, \partial M)\geq 4\ov C\varrho$, then
$$
\varphi_{a, \l}(\cdot)=\hat \d_{a, \l}(\cdot)+\ln \frac{\l}{2}+H(a, \cdot)+\frac{1}{4\l^2}\D_{g_{a}}H(a, \cdot)+O\left(\frac{1}{\l^3}\right),
$$
and if $a\in \partial M$, then
$$
\varphi_{a, \l}(\cdot)=\hat \d_{a, \l}(\cdot)+\ln \frac{\l}{2}+\frac{1}{2}H(a, \cdot)+\frac{1}{8\l^2}\D_{g_{a}}H(a, \cdot)+O\left(\frac{1}{\l^3}\right),
$$
where\;$O(1)$ means $O_{a, \l, \epsilon}(1)$ and for it meaning see Section \ref{eq:notpre}.\\\\
2)\;If $d_g(a, \partial M)\geq 4\ov C\varrho$, then
$$
\l\frac{\partial \varphi_{a, \l}(\cdot)}{\partial \l}=\frac{2}{1+\l^2\chi_{\varrho}^2(d_{g_{a}}(a, \cdot))}-\frac{1}{2\l^2}\D_{g_{a}}H(a, \cdot)+O\left(\frac{1}{\l^3}\right),
$$
and if $a\in \partial M$, then
$$
\l\frac{\partial \varphi_{a, \l}(\cdot)}{\partial \l}=\frac{2}{1+\l^2\chi_{\varrho}^2(d_{g_{a}}(a, \cdot))}-\frac{1}{4\l^2}\D_{g_{a}}H(a, \cdot)+O\left(\frac{1}{\l^3}\right),
$$
where $O(1)$ is as in point 1).\\\\
3)If $d_g(a, \partial M)\geq 4\ov C\varrho$, then
$$
\frac{1}{\l}\frac{\partial \varphi_{a, \l}(\cdot)}{\partial a}=\frac{\chi_{\varrho}(d_{g_{a}}(a, \cdot))\chi_{\varrho}^{'}((d_{g_{a}}(a, \cdot))}{d_{g_{a}}(a, \cdot)}\frac{2\l exp_a^{-1}(\cdot)}{1+\l^2\chi_{\varrho}^2(d_{g_a}(a, \cdot))}+\frac{1}{\l}\frac{\partial  H(a, \cdot)}{\partial a}+O\left(\frac{1}{\l^3}\right),
$$
and if $a\in \partial M$, then$$
\frac{1}{\l}\frac{\partial \varphi_{a, \l}(\cdot)}{\partial a}=\frac{\chi_{\varrho}(d_{g_{a}}(a, \cdot))\chi_{\varrho}^{'}((d_{g_{a}}(a, \cdot))}{d_{g_{a}}(a, \cdot)}\frac{2\l exp_a^{-1}(\cdot)}{1+\l^2\chi_{\varrho}^2(d_{g_a}(a, \cdot))}+\frac{1}{2\l}\frac{\partial  H(a, \cdot)}{\partial a}+O\left(\frac{1}{\l^3}\right),
$$
where \;$O(1)$ is as in point 1).
\end{lem}
\vspace{4pt}

\noindent
\begin{lem}\label{eq:outbubbleest}
Assuming that $\epsilon$ is small and positive, $0<\varrho<\varrho_k$ and \;$0<\eta<\eta_k$ where $0<\varrho<\varrho_k$ and \;$0<\eta<\eta_k$ are as in Proposition \ref{eq:sharpest}, and $\l\geq \frac{1}{\epsilon}$, then we have:\\\\
1)) If $d_g(a, \partial M)\geq 4\ov C\varrho$, then
$$
\varphi_{a, \l}(\cdot)=G(a, \cdot)+\frac{1}{4\l^2}\D_{g_{a}}G(a, \cdot)+O\left(\frac{1}{\l^3}\right)\;\;\text{in} \;\;\;M\setminus B^{a}_{a}(\eta),
$$
and if $a\in \partial M$, then
$$
\varphi_{a, \l}(\cdot)=\frac{1}{2}G(a, \cdot)+\frac{1}{8\l^2}\D_{g_{a}}G(a, \cdot)+O\left(\frac{1}{\l^3}\right)\;\;\text{in} \;\;\;M\setminus B^{a}_{a}(\eta).
$$
where $O(1)$ means $O_{a, \l, \epsilon}(1)$ and for it meaning see Section \ref{eq:notpre}.\\\\
2) If $d_g(a, \partial M)\geq 4\ov C\varrho$, then
$$
\l\frac{\partial \varphi_{a, \l}(\cdot)}{\partial \l}=-\frac{1}{2\l^2}\D_{g_{a}}G_(a, \cdot)+O\left(\frac{1}{\l^3}\right)\;\;\text{in} \;\;\;M\setminus B^{a}_{a}(\eta),
$$
and if $a\in \partial M$, then
$$
\l\frac{\partial \varphi_{a, \l}(\cdot)}{\partial \l}=-\frac{1}{4\l^2}\D_{g_{a}}G_(a, \cdot)+O\left(\frac{1}{\l^3}\right)\;\;\text{in} \;\;\;M\setminus B^{a}_{a}(\eta),
$$
where  $O(1)$ is as in point 1).\\\\
3) If $d_g(a, \partial M)\geq 4\ov C\varrho$, then
$$
\frac{1}{\l}\frac{\partial \varphi_{a, \l}(\cdot)}{\partial a}=\frac{1}{\l}\frac{\partial G(a, \cdot)}{\partial a}+O\left(\frac{1}{\l^3}\right)\;\;\text{in} \;\;\;M\setminus B^{a}_{a}(\eta),
$$
and if $a\in \partial M$, then $$
\frac{1}{\l}\frac{\partial \varphi_{a, \l}(\cdot)}{\partial a}=\frac{1}{2\l}\frac{\partial G(a, \cdot)}{\partial a}+O\left(\frac{1}{\l^3}\right)\;\;\text{in} \;\;\;M\setminus B^{a}_{a}(\eta),
$$
where $O(1)$ is as in point 1).
\end{lem}
\vspace{4pt}

\noindent
Next, using the above two Lemmas, we obtain the following three Lemmas.
\begin{lem}\label{eq:selfintest}
Assuming that $\epsilon$ is small and positive, $0<\varrho<\varrho_k$ where $\varrho_k$ is as in Proposition \ref{eq:sharpest}, and $\l\geq \frac{1}{\epsilon}$, then we have:\\\\
1) If $d_g(a, \partial M)\geq 4\ov C\varrho$, then
$$
<P_g\varphi_{a, \l}, \varphi_{a, \l}>=32\pi^2\ln \l-\frac{40\pi^2}{3}+16\pi^2 H(a, a)+\frac{8\pi^2}{\l^2}\D_{g_a} H(a, a)+O\left(\frac{1}{\l^3}\right),
$$
and if $a\in \partial M$, then
$$
<P_g\varphi_{a, \l}, \varphi_{a, \l}>=16\pi^2\ln \l-\frac{20\pi^2}{3}+4\pi^2 H(a, a)+\frac{2\pi^2}{\l^2}\D_{g_a} H(a, a)+O\left(\frac{1}{\l^3}\right),
$$
where \;$O(1)$\; means \;$O_{a, \l, \epsilon}(1)$\; and for its meaning see Section \ref{eq:notpre}.\\\\
2) If $d_g(a, \partial M)\geq 4\ov C\varrho$, then
$$
<P_g\varphi_{a, \l}, \l\frac{\varphi_{a, \l}}{\partial \l} >=16\pi^2-\frac{8\pi^2}{\l^2}\D_{g_a}H(a, a)+O\left(\frac{1}{\l^3}\right),
$$
and if $a\in \partial M$, then
$$
<P_g\varphi_{a, \l}, \l\frac{\varphi_{a, \l}}{\partial \l} >=8\pi^2-\frac{2\pi^2}{\l^2}\D_{g_a}H(a, a)+O\left(\frac{1}{\l^3}\right),
$$
where $O(1)$ is as in point 1).\\\\
3) If $d_g(a, \partial M)\geq 4\ov C\varrho$, then
$$
<P_g\varphi_{a, \l}, \frac{1}{\l}\frac{\varphi_{a, \l}}{\partial a} >=\frac{16\pi^2}{\l}\frac{\partial H(a, a)}{\partial a}+O\left(\frac{1}{\l^3}\right),
$$
and if $a\in \partial M$, then
$$
<P_g\varphi_{a, \l}, \frac{1}{\l}\frac{\varphi_{a, \l}}{\partial a} > \, =  \, 3 \pi^2 \, + \, \frac{4\pi^2}{\l}\frac{\partial H(a, a)}{\partial a}+O\left(\frac{1}{\l^3}\right).
$$
where \;$O(1)$\; is as in point 1).
\end{lem}
\vspace{4pt}

\noindent
\begin{lem}\label{eq:interactest}
Assuming that $\epsilon$ is small and positive, $0<\varrho<\varrho_k$, $0<\eta<\eta_k$, $\L>\L_k$ where $0<\varrho<\varrho_k$, \;$0<\eta<\eta_k$, and $\L_k$ are as in Proposition \ref{eq:sharpest}, $a_i, a_j\in M$, \;$d_g(a_i, a_j)\geq 4\ov C\eta$, $\frac{1}{\L}\leq \frac{\l_i}{\l_j}\leq \L$, and $\l_i, \l_j\geq \frac{1}{\epsilon}$ where $\ov C$ is as in \eqref{eq:proua}, then we have:\\\\
1) \\
$(1)_{i, j}$\;\;\;If \;$d_g(a_i, \partial M)\geq 4\ov C\varrho$ \;and \;$d_g(a_j, \partial M)\geq 4\ov C\varrho$, then
$$
<P_g\varphi_{a_i, \l_i}, \varphi_{a_j, \l_j}>=16\pi^2G (a_j, a_i)+\frac{4\pi^2}{\l_i^2}\D_{g_{a_i}} G(a_i, a_j)++\frac{4\pi^2}{\l_j^2}\D_{g_{a_j}} G(a_j, a_i)+O\left(\frac{1}{\l^3_i}+\frac{1}{\l_j^3}\right),
$$
$(2)_{i, j}$\;\;\;if \;$d_g(a_i, \partial M)\geq 4\ov C\varrho$\; and \;$a_j\in\partial M$, then
$$
<P_g\varphi_{a_i, \l_i}, \varphi_{a_j, \l_j}>=8\pi^2G (a_j, a_i)+\frac{2\pi^2}{\l_i^2}\D_{g_{a_i}} G(a_i, a_j)++\frac{2\pi^2}{\l_j^2}\D_{g_{a_j}} G(a_j, a_i)+O\left(\frac{1}{\l^3_i}+\frac{1}{\l_j^3}\right),
$$
and\\\\
$(3)_{i, j}$\;\;\;if $a_i\in \partial M$ and $a_j\in \partial M$, then
$$
<P_g\varphi_{a_i, \l_i}, \varphi_{a_j, \l_j}>=4\pi^2G (a_j, a_i)+\frac{\pi^2}{\l_i^2}\D_{g_{a_i}} G(a_i, a_j)++\frac{\pi^2}{\l_j^2}\D_{g_{a_j}} G(a_j, a_i)+O\left(\frac{1}{\l^3_i}+\frac{1}{\l_j^3}\right),
$$
where \;$O(1)$\; means here \;$O_{A, \bar \l, \epsilon}(1)$\; with \;$A=(a_i, a_j)$\; and \;$\bar \l=(\l_i, \l_j)$\; and for the meaning of $O_{A, \bar \l, \epsilon}(1)$, see Section \ref{eq:notpre}.\\\\
2)\\
$(1)_{i, j}$\;\;\;If \;$d_g(a_i, \partial M)\geq 4\ov C\varrho$ \;and \;$d_g(a_j, \partial M)\geq 4\ov C\varrho$, then
$$
<P_g\varphi_{a_i, \l_i}, \l_j\frac{\partial \varphi_{a_j, \l_j}}{\partial \l_j}>=-\frac{8\pi^2}{\l_j^2}\D_{g_{a_j}}G(a_j, a_i)+O\left(\frac{1}{\l^3_j}\right),
$$
$(2)_{i, j}$\;\;\;if \;$d_g(a_i, \partial M)\geq 4\ov C\varrho$\; and \;$a_j\in\partial M$, then
$$
<P_g\varphi_{a_i, \l_i}, \l_j\frac{\partial \varphi_{a_j, \l_j}}{\partial \l_j}>=-\frac{4\pi^2}{\l_j^2}\D_{g_{a_j}}G(a_j, a_i)+O\left(\frac{1}{\l^3_j}\right),
$$
$(2)^{'}_{i, j}$\;\;\;if \;$a_i\in\partial M$ and \;$d_g(a_j, \partial M)\geq 4\ov C\varrho$\;, then
$$
<P_g\varphi_{a_i, \l_i}, \l_j\frac{\partial \varphi_{a_j, \l_j}}{\partial \l_j}>=-\frac{4\pi^2}{\l_j^2}\D_{g_{a_j}}G(a_j, a_i)+O\left(\frac{1}{\l^3_j}\right),
$$
$(3)_{i, j}$\;\;\;if $a_i\in \partial M$ and $a_j\in \partial M$, then
$$
<P_g\varphi_{a_i, \l_i}, \l_j\frac{\partial \varphi_{a_j, \l_j}}{\partial \l_j}>=-\frac{2\pi^2}{\l_j^2}\D_{g_{a_j}}G(a_j, a_i)+O\left(\frac{1}{\l^3_j}\right),
$$
where $O(1)$ is as in point 1).\\\\
3)\\
$(1)_{i, j}$\;\;\;If \;$d_g(a_i, \partial M)\geq 4\ov C\varrho$ \;and \;$d_g(a_j, \partial M)\geq 4\ov C\varrho$, then
$$
<P_g\varphi_{a_i, \l_i}, \frac{1}{\l_j}\frac{\partial \varphi_{a_j, \l_j}}{\partial a_j}>= \frac{16\pi^2}{\l_j}\frac{\partial G(a_j, a_i)}{\partial a_j}+O\left(\frac{1}{\l^3_j}\right),
$$
$(2)_{i, j}$\;\;\;if \;$d_g(a_i, \partial M)\geq 4\ov C\varrho$\; and \;$a_j\in\partial M$,
$$
<P_g\varphi_{a_i, \l_i}, \frac{1}{\l_j}\frac{\partial \varphi_{a_j, \l_j}}{\partial a_j}>= \frac{8\pi^2}{\l_j}\frac{\partial G(a_j, a_i)}{\partial a_j}+O\left(\frac{1}{\l^3_j}\right),
$$
$(2)^{'}_{i, j}$\;\;\;if \;$a_i\in\partial M$ and\;$d_g(a_j, \partial M)\geq 4\ov C\varrho$, then
$$
<P_g\varphi_{a_i, \l_i}, \frac{1}{\l_j}\frac{\partial \varphi_{a_j, \l_j}}{\partial a_j}>= \frac{8\pi^2}{\l_j}\frac{\partial G(a_j, a_i)}{\partial a_j}+O\left(\frac{1}{\l^3_j}\right),
$$

$(3)_{i, j}$\;\;\;if $a_i\in \partial M$ and $a_j\in \partial M$, then
$$
<P_g\varphi_{a_i, \l_i}, \frac{1}{\l_j}\frac{\partial \varphi_{a_j, \l_j}}{\partial a_j}>= \frac{4\pi^2}{\l_j}\frac{\partial G(a_j, a_i)}{\partial a_j}+O\left(\frac{1}{\l^3_j}\right).
$$

\end{lem}
\vspace{4pt}



\begin{lem}\label{l:A1}
1) If $\epsilon$ is small and positive, $a \in \partial M$, $q\in \N^*$, and $\l\geq \frac{1}{\epsilon}$ , then there holds
\begin{equation}
C^{-1}\l^{8q-4}\leq \int_Me^{4p\varphi_{a, \l}}dV_g\leq C\l^{8q-4},
\end{equation}
where $C$ is independent of $a$, $\l$, and $\epsilon$.\\\\
2) If $\epsilon$ is positive and small, $a_i, a_j \in \partial M$, $\l\geq \frac{1}{\epsilon}$ and $\l d_g(a_i,a_j) \geq 4\ov C R$, then we have
\begin{equation}
<P_g^{4, 3} \varphi_{a_i,\l},\varphi_{a_j,\l}> \, \leq \,4\pi^2 G(a_i,a_j) \, + \, O(1),
\end{equation}
where \;$O(1)$\; means here \;$O_{A,  \l, \epsilon}(1)$\; with \;$A=(a_i, a_j)$, and for the meaning of $O_{A, \l, \epsilon}(1)$, see section \ref{eq:notpre}.
3) If $\epsilon$ is positive and small, $a_i, a_j \in \partial M$, $\l_i, \l_j\geq \frac{1}{\epsilon}$, $\frac{1}{\L}\leq\frac{\L_i}{\l_j}\leq \L$ and $\l _id_g(a_i,a_j) \geq 4\ov C R$, then we have
\begin{equation}
<P_g^{4, 3} \varphi_{a_i,\l_i},\varphi_{a_j,\l_j}> \, \leq \, 4\pi^2  G(a_i,a_j) \, + \, O(1),
\end{equation}
where \;$O(1)$\; means here \;$O_{A, \bar \l, \epsilon}(1)$\; with \;$A=(a_i, a_j)$\; and \;$\bar \l=(\l_i, \l_j)$\; and for the meaning of $O_{A, \bar \l, \epsilon}(1)$, see Section \ref{eq:notpre}.
\end{lem}

\begin{lem}\label{l:A2}
Let  $q\in \N^*$,  $\hat R$ be a large positive constant, $\epsilon$ be a small positive number, $\alpha_i\geq 0$, $i=1, \cdots, q$, $\sum_{i=1}^q\alpha_i=k$, $\l\geq \frac{1}{\epsilon}$ and $ u = \sum_{i=1}^p \a_i \varphi_{a_i,\l}$ with $a_i\in \partial M$ for $i=1, \cdots, q$. Assuming that there exist two positive integer $i, j\in \{1, \cdots, p\}$ with  $i\neq j$ such that $\l d_g(a_i,a_j) \leq \frac{\hat R}{4\ov C}$, where $\ov C$ is as in \eqref{eq:proua}, then we have
\begin{equation}
II(u) \, \leq II(v) \, + \, O(\ln \hat R),
\end{equation}
with
$$
 v:= \sum_{k \leq p, k \neq i,j} \a_k \varphi_{a_k,\l} \, + (\a_i + \a_j) \varphi_{a_i,\l}.
$$
where here \;$O(1)$\;stand for \;$O_{{\bar\alpha}, A, \l, \epsilon}(1)$, with \;${\bar\alpha}=(\alpha_1, \cdots, \alpha_q)$\; and \;$A=(a_1, \cdots, a_q)$, and for the meaning of \;$O_{{\bar \alpha}, A, \l, \epsilon}(1)$, we refer the reader to Section \ref{eq:notpre}.
\end{lem}\begin{lem}\label{eq:lemma3}
1) If  $\epsilon$ is positive and small, $a_i, a_j\in \partial M$, $\l\geq \frac{1}{\epsilon}$ and $\l d_g(a_i, a_j)\geq 4\ov C R$, then
$$
\varphi_{a_j, \l}(\cdot)=\frac{1}{2}G(a_j,\cdot)+O(1)\;\;\;\text{in}\;\;B^{a_i}_{a_i}(\frac{R}{\l}),
$$
where here \;$O(1)$ means here \;$O_{ A, \l, \epsilon}(1)$, with  \;$A=(a_i, a_j)$, and for the meaning of \;$O_{ A, \l, \epsilon}(1)$, see section \ref{eq:notpre}.\\
2) If $\epsilon$ is positive and small, $a_i, a_j\in \partial M$, $\l_i, \l_j\geq \frac{1}{\epsilon}$, $\frac{1}{\L}\leq\frac{\L_i}{\l_j}\leq \L$, and $\l_i d_g(a_i, a_j)\geq 4\ov CR$, then
$$
\varphi_{a_j, \l_j}(\cdot)=\frac{1}{2}G(a_j,\cdot)+O(1)\;\;\;\text{in}\;\;B^{a_i}_{a_i}(\frac{R}{\l_i}),
$$
where here \;$O(1)$ means here \;$O_{ A, \bar\l, \epsilon}(1)$, with  \;$A=(a_i, a_j)$, $\bar \l=(\l_i, \l_j)$ and for the meaning of \;$O_{ A, \l, \epsilon}(1)$, see Section \ref{eq:notpre}.\\
\end{lem}
\vspace{4pt}

\begin{lem}\label{l:Keu}
Let $ \ov{u} = \sum_{i=1}^p \a_i \varphi_i \, +  \sum_{i= p+ 1}^{p+q}  \a_i \varphi_i \, + \sum_{r = 1}^{\b} \b_r v_r \in V(p,q,\e, \varrho, \eta)$. Then there holds

\begin{equation}\label{eq:Keu}
\int_M Ke^{4\ov{u}} dv_g \, = \, \G  \left [   1 \, + \frac{1}{2}  \sum_{i=p+1}^{p+q}  \frac{ \frac{\partial \mathcal{F}_i}{\partial n_{g_i}}(a_i) \, \G_i   }{ \G \l_i \mathcal{F}_i(a_i)}+ \, O(\frac{1}{\l^2})\right ]
\end{equation}
where
$$
\G_i:= \frac{\pi^2}{4}  \frac{\l_i^{8 \a_i -4} \mathcal{F}_i(a_i)  \mathcal{G}_i(a_i)}{ (2 \a_i -1)( 4 \a_i -1)},  \quad
\G :=   2 \sum_{i=1}^p \G_i \, + \, \sum_{i= p + 1}^{p +q} \G_i
$$

and $\mathcal{G}_i(a_i)$ is defined in \eqref{eq:Gi1} and \eqref{eq:Gi2}.
Moreover  setting
$$
\tau_i:= \, 1 - \, \frac{k \G_i}{ \int_M K e^{4 \ov{u}}}
$$
there holds
\begin{equation}\label{eq:tau}
\sum_{i=1}^p 2 \tau_i \, + \,  \sum_{i= p +1}^{p +q}   \tau_i \, = \, \frac{1}{2} \, \sum_{i=p+1}^{p+q} \frac{\frac{\partial \mathcal{F}_i}{\partial n_{g_i}}(a_i)   }{\l_i \mathcal{F}_i(a_i)}  \, + \,  O( \sum_k |\tau_k|^2 \, + \,  \frac{1}{\l^2})
\end{equation}
\end{lem}

\bigskip

\noindent
Next, using the relation between $(\R^4_+, g_{\R^4_+})$ and $(S^4_+, g_{S^4_+})$,  the spectral property of the Paneitz operator $P_{g_{S^4}}$, a standard doubling argument to deal with boundary points, and a standard  blow-up argument ( as in Brendle \, \cite{bren3}), we obtain the following last two lemmas of this subsection.
\begin{lem}\label{eq:positive}
Assuming that $0<\varrho<\varrho_k$, then there exists \;$\Gamma_0:=\Gamma_0(\varrho)$\; and \;$\tilde \L_0:=\tilde\L_0(\varrho)$ two large positive constant such that for every $a\in M$ such that either $d_g(a, \partial M)\geq4\ov C\varrho$ or $a\in \partial M$, $\l\geq \tilde \L_0$, and \;$w\in F_{a,  \l}:=\{w\in \mathcal{H}_{\frac{\partial}{\partial n}}: \ov w_{(Q, T)}=<\varphi_{a, \l}, w>_{P^{4, 3}}=0\}$, we have
\begin{equation}\label{eq:positivity}
 \int_Me^{4\hat\d_{a,\l}}w^2dV_{g_{a}}\leq \Gamma_0||w||^2.
\end{equation}
\end{lem}
\vspace{4pt}

\noindent
\begin{lem}\label{eq:positiveg}
Assuming that $0<\varrho<\varrho_k$, $0<\eta<\eta_k$, then there exists a small positive constant $c_0:=c_0(\varrho, \eta)$ and a large positive constant $\L_0:=\L_0(\varrho, \eta)$ such that for every  $(p, q)\in \N^2$ such that $2p+q=k$, for every $a_i\in M$ concentrations points $i=1, \cdots, p+q$ such that $d_g(a_i, a_j)\geq 4\ov C\eta$ for $i\neq j=1, \cdots, p+q$, $d_g(a_i, \partial M)\geq 4\ov C\varrho$, $i=1, \cdots, p$, where $\bar C$ is as in \eqref{eq:proua}, and  $a_i\in \partial M$, $i=p+1, p+q$, for every \;$\l_i>0$\; concentrations parameters satisfying $\l_i\geq \L_0$, with $i=1, \cdots, p+q$, and for every \;$w\in E_{A, \bar \l}^*=\cap_{i=1}^{p+q} E^*_{a_i, \l_i}$ with $A:=(a_1, \cdots, a_{p+q}$), $\bar \l:= (\l_1, \cdots, \l_{p+q})$ and $E^*_{a_i, \l_i}=\{w\in \mathcal{H}_{\frac{\partial}{\partial n}}:\;<\varphi_{a_i, \l_i}, w>_{P^{4, 3}}=<\frac{\partial\varphi_{a_i, \l_i}}{\partial \l_i}, w>_{P^{4, 3}}=<\frac{\partial\varphi_{a_i, \l_i}}{\partial a_i}, w>_{P^{4, 3}}=\ov{w}_{(Q, T)}=0\}$, there holds
\begin{equation}\label{eq:positivity}
||w||^2-24 \sum_{i=1}^{p+q}\int_Me^{4\hat\d_{a_i,\l_i}}w^2dV_{g_{a_i}}\geq c_0||w||^2.
\end{equation}
\end{lem}



\subsection{ Gradient  and energy estimates}

This section is devoted to the expansion of the functional $II$ and its gradient in the neighborhood of potential critical points at infinity.


\subsubsection{Expansion of the Euler-Lagrange functional near Infinity}\label{s:expe}
In this subsection, we derive an expansion of the Euler-Lagrange functional $II$ for the part at infinity which is characterized by a useful (topologically) piece of $B^p_q(M, \partial M)$ via the variational bubbles $\varphi_{a, \l}$ (where $B^p_q(M, \partial M)$ is as in \eqref{eq:barypq}), namely for elements at infinity with zero $w$-part in the representation \eqref{eq:para}. Indeed, we have:
\begin{lem}\label{eq:energyest}
Assuming that $(p, q)\in \N^2$  such that $2p+q=k$, $0<\varrho<\varrho_k$, $0<\eta<\eta_k$, and $0<\epsilon\leq \epsilon_k$, where $\varrho_k$ and $\eta_k$ are given by Proposition \ref{eq:sharpest}, and the $\epsilon_k$ is given by \eqref{eq:mini}, then for $a_i\in M$ concentration points,  $\alpha_i$ masses, $\l_i$ concentration parameters ($i=1,\cdots, p+q$), and $\beta_r$ negativity parameters ($r=1, \cdots, \bar k$) satisfying \eqref{eq:afpara}, we have
\begin{equation*}
\begin{split}
II(\sum_{i=1}^{p+q}\alpha_i\varphi_{a_i, \l_i}+\sum_{r=^1}^{\bar k} \beta_r(v_r-\ov{(v_r)}_{(Q, T)}))  =  -\frac{20}{3}k\pi^2-4k\pi^2\ln(\frac{k\pi^2}{6}) -8\pi^2\mathcal{F}_{p, q}(a_1, \dots, a_{p+q})   \\ 16 \pi^2 \left ( 2\sum_{i=1}^p(\alpha_i-1)^2 \ln \l_i  \, +\sum_{i=p+1}^{p+q}(\alpha_i-1)^2 \ln \l_i \right )  +  \sum_{r=1}^{\ov{k}} \mu_k \b^2_k \, - 4 \pi^2 \left [ \sum_{i=1}^p 2 \sigma_i^2 \, +\, \sum_{i= p+1}^{p+q}  \sigma_i^2     \right ] \\
 -  2 \pi^2 \sum_{i=p+1}^{p+q}\frac{1}{\l_i} \, \frac{1}{\mathcal{F}_i^A(a_i)}\frac{\partial \mathcal{F}_i^A}{\partial n_{g_{a_i}}}(a_i)   + O\left(   \sum_{r=1}^{\bar k}|\beta_r|^3 \, +  \, \sum_{k=1}^{p+q}|\alpha_k-1|^3+|\sigma_k|^3 + \frac{1}{\l_k^2}\right),
\end{split}
\end{equation*}



where $\mathcal{F}_{p, q}$ is as in \eqref{eq:limitfsintpq}, \;$O\left(1\right)$ means here \;$O_{\bar\alpha, A, \bar \l, \bar \beta, \epsilon}\left(1\right)$ \;with \;$\bar\alpha=(\alpha_1, \cdots, \alpha_{p+q})$, $A:=(a_1, \cdots, a_{p+q})$ $\bar \l:=(\l_1, \cdots, \l_{p+q})$, $\bar \beta:=(\beta_1, \cdots, \beta_{\bar k})$ and for $i=1, \cdots, p+q$, $$\tilde \sigma_i:=1-\frac{k\Gamma_i}{\Gamma},\;\;\;\; \Gamma:=2\sum_{i=1}^p\gamma_i+\sum_{i=p+1}^{p+q}\gamma_i, \;\;\;\Gamma_i:= \frac{\pi^2}{4} \frac{\l_i^{8\alpha_i-4}\mathcal{F}^{A}_i(a_i)\mathcal{G}_i(a_i)}{(2 \a_i -1)(4 \a_i -1)}, $$
 with for $i=1, \cdots, p$, $\mathcal{F}_i^A$ is as in \eqref{eq:partiallimitint},

\begin{equation}\label{eq:Gi1}
\begin{split}
 \mathcal{G}_i(a_i):=e^{4((\alpha_i-1)H(a_i, a_i)+\sum_{j=1, j\neq i}^p(\alpha_j-1)G(a_j, a_i))+\frac{1}{2}\sum_{j=p+1}^{p+q}(\alpha_j-1)G(a_j, a_i))}\\\times e^{\sum_{j=1, j\neq i}^p\frac{\alpha_j}{\l_j^2}\D_{g_{a_j}}G(a_j, a_i)}e^{\frac{1}{2}\sum_{j=p+1}^{p+q}\frac{\alpha_j}{\l_j^2}\D_{g_{a_j}}G(a_j, a_i)}e^{\frac{\alpha_i}{\l_i^2}\D_{g_{a_i}}H(a_i, a_i)}e^{4\sum_{r=1}^{\bar k}\beta_rv_r(a_i)},
\end{split}
\end{equation}
for $i=p+1, \cdots, p+q$, $\mathcal{F}_i^A$ is as in \eqref{eq:partiallimitbound},
\begin{equation}\label{eq:Gi2}
\begin{split}
 \mathcal{G}_i(a_i):=e^{4(\frac{1}{2}(\alpha_i-1)H(a_i, a_i)+\sum_{j=p+1, j\neq i}^{p+q}\frac{1}{2}(\alpha_j-1)G(a_j, a_i))+\sum_{j=1}^{p}(\alpha_j-1)G(a_j, a_i))}\\\times e^{\sum_{j=1}^p\frac{\alpha_j}{\l_j^2}\D_{g_{a_j}}G(a_j, a_i)}e^{\frac{1}{2}\sum_{j=p+1, j\neq i}^{p+q}\frac{\alpha_j}{\l_j^2}\D_{g_{a_j}}G(a_j, a_i)}e^{\frac{\alpha_i}{2\l_i^2}\D_{g_{a_i}}H(a_i, a_i)}e^{4\sum_{r=1}^{\bar k}\beta_rv_r(a_i)},
\end{split}
\end{equation}

\end{lem}



\subsubsection{Expansion of the gradient near infinity}
As mentioned above, in this subsection, we perform an expansion of $\n II$ on the same set as in the previous subsection. To do so, we start with the gradient of $II$ in the direction of $\bar \l$. Precisely, we have:
\begin{lem}\label{eq:gradientlambdaest}
Assuming that $(p, q)\in \N^2$  such that $2p+q=k$, $0<\varrho<\varrho_k$, $0<\eta<\eta_k$, and $0<\epsilon\leq \epsilon_k$, where $\varrho_k$ and $\eta_k$ are given by Proposition \ref{eq:sharpest}, and $\epsilon_k$ is given by \eqref{eq:mini}, then for $a_i\in M$ concentration points,  $\alpha_i$ masses, $\l_i$ concentration parameters ($i=1,\cdots, p+q$), and $\beta_r$ negativity parameters ($r=1, \cdots, \bar k$) satisfying \eqref{eq:afpara}, we have that for every $j=1,\cdots, p$, there holds
\begin{equation*}
\begin{split}
&<\n II(\sum_{i=1}^{p+q}\alpha_i\varphi_{a_i, \l_i}+\sum_{r=1}^{\bar k}\beta_r(v_r-\ov{(v_r)}_{(Q, T)}), \l_j\frac{\partial \varphi_{a_j, \l_j}}{\partial \l_j}>= 32\pi^2\alpha_j\tau_j-\frac{4\pi^2}{\l_j^2}\left(\frac{\D_{g_{a_j}} \mathcal{F}^{A}_j(a_j)}{\mathcal{F}^{A}_j(a_j)}-\frac{2}{3}R_g(a_j)\right)\\&
+O\left(\sum_{i=1}^{p+q}|\alpha_i-1|^2+ \sum_{i=1}^{p+q}  |\tau_i|^2 + \sum_{r=1}^{\bar k}|\beta_r|^3+\sum_{i=1}^{p+q}\frac{1}{\l_i^3}\right),
\end{split}
\end{equation*}
and for every $j=p+1,\cdots, p+q$, there holds
\begin{equation*}
\begin{split}
&<\n II(\sum_{i=1}^{p+q}\alpha_i\varphi_{a_i, \l_i}+\sum_{r=1}^{\bar k}\beta_r(v_r-\ov{(v_r)}_{(Q, T)}), \l_j\frac{\partial \varphi_{a_j, \l_j}}{\partial \l_j}>= 16\pi^2\alpha_j\tau_j-\frac{6\pi^2}{\l_j} \frac{1}{\mathcal{F}_j^A(a_j)}\frac{\partial \mathcal{F}_j^A}{\partial n_{g_{a_j}}}(a_j)\\&
+O\left(\sum_{i=1}^{p+q}|\alpha_i-1|^2+  \sum_{i=1}^{p+q}  |\tau_i|^2 + \sum_{r=1}^{\bar k}|\beta_r|^3+\sum_{i=1}^{p+q}\frac{1}{\l_i^2}\right),
\end{split}
\end{equation*}
 where $A:=(a_1, \cdots, a_{p+q})$, $O\left(1\right)$  and for $i=1, \cdots, p+q$, $\mathcal{F}_i^A$ is defined in \eqref{eq:partiallimitint}  ( resp. in \eqref{eq:partiallimitbound}) and $$\tau_i:=1-\frac{k\Gamma_i}{D},  \mbox{ where }  \Gamma_i := \frac{\pi^2}{4 \a_i(2\a_i -1)(4\a_i -1)}\; \;\;D:=\int_M K(x)e^{4(\sum_{i=1}^{p+q}\alpha_i\varphi_{a_i, \l_i}(x)+\sum_{r=1}^{\bar k}\beta_r v_r(x))}dV_g(x).$$
\end{lem}

\noindent
Using Lemma \ref{eq:gradientlambdaest}, we have the following corollary:
\begin{cor}\label{eq:cgradientlambdaest}
Assuming that $(p, q)\in \N^2$  such that $2p+q=k$ and $q \neq 0 $ and  $0<\varrho<\varrho_k$, $0<\eta<\eta_k$, and $0<\epsilon\leq \epsilon_k$, where $\varrho_k$ and $\eta_k$ are given by Proposition \ref{eq:sharpest}, and $\epsilon_k$ is given by \eqref{eq:mini}, then for $a_i\in M$ concentration points,  $\alpha_i$ masses, $\l_i$ concentration parameters ($i=1,\cdots, p+q$), and $\beta_r$ negativity parameters ($r=1, \cdots, \bar k$) satisfying \eqref{eq:afpara}, we have that
\begin{equation*}
 \begin{split}
&<\n II(\sum_{i=1}^{p+q}\alpha_i\varphi_{a_i, \l_i}+\sum_{r=1}^{\bar k}\beta_r(v_r-\ov{(v_r)}_{(Q, T)}), \sum_{j=1}^{p+q} \frac{\l_j}{\a_j} \frac{\partial \varphi_{a_j, \l_j}}{\partial \l_j}>= 2 \pi^2\sum_{i=p+1}^{p+q} \frac{1}{\l_i\mathcal{F}_i^A(a_i)}\frac{\partial \mathcal{F}_i^A}{\partial n_{g_{a_i}}}(a_i)\\&
+O\left(\sum_{i=1}^{p+q}|\alpha_i-1|^2+\sum_{r=1}^{\bar k}|\beta_r|^3+\sum_{i=p+1}^{p+q}\frac{|\tau_i|}{\l_i}\left|\frac{\partial \mathcal{F}_i^A}{\partial n_{g_{a_i}}}(a_i)\right|+\sum_{i=1}^{p+q}|\tau_i|^2 +\sum_{i=1}^{p+q}\frac{1}{\l_i^2}\right)
\end{split}
\end{equation*}
where $A:=(a_1, \cdots, a_{p+q})$, $O\left(1\right)$ as as in Lemma \ref{eq:energyest}, and for $i=1, \cdots, p+q$,  $\mathcal{F}_i^A$ is as in Lemma \ref{eq:energyest}.
\end{cor}
\vspace{4pt}




\noindent
Next, we give  the estimate of the gradient of the Euler-Lagrange functional $II$ in the direction of $A$. Precisely, we have that:
\begin{lem}\label{eq:gradientaest}
Assuming that $(p, q)\in \N^2$  such that $2p+q=k$, $0<\varrho<\varrho_k$, $0<\eta<\eta_k$, and $0<\epsilon\leq \epsilon_k$, where $\varrho_k$ and $\eta_k$ are given by Proposition \ref{eq:sharpest}, and $\epsilon_k$ is given by \eqref{eq:mini}, then for $a_i\in M$ concentration points,  $\alpha_i$ masses, $\l_i$ concentration parameters ($i=1,\cdots, p+q$), and $\beta_r$ negativity parameters ($r=1, \cdots, \bar k$) satisfying \eqref{eq:afpara}, we have that for every $j=1,\cdots, p$, there holds
\begin{equation*}
\begin{split}
<\n II(\sum_{i=1}^{p+q}\alpha_i\varphi_{a_i, \l_i}+\sum_{r=1}^{\bar k}\beta_r(v_r-\ov{(v_r)}_{(Q, T)}), \frac{1}{\l_j}\frac{\partial \varphi_{a_j, \l_j}}{\partial a_j}>&=-\frac{4\pi^2}{\l_j}\frac{\n_g\mathcal{F}_j^{A}(a_j)}{\mathcal{F}_j^{A}(a_j)}\\&+O\left(\sum_{i=1}^{p+q}|\alpha_i-1|^2+\sum_{i=1}^{p+q}\frac{1}{\l_i^2}+\sum_{r=1}^{\bar k}|\beta_r|^2+\sum_{i=1}^{p+q}\tau_i^2\right),
\end{split}
\end{equation*}
for every $j=p+1, \cdots, p+q$, there holds
\begin{equation*}
\begin{split}
<\n II(\sum_{i=1}^{p+q}\alpha_i\varphi_{a_i, \l_i}+\sum_{r=1}^{\bar k}\beta_r(v_r-\ov{(v_r)}_{(Q, T)}), \frac{1}{\l_j}\frac{\partial \varphi_{a_j, \l_j}}{\partial a_j}>&= \, 6 \pi^2 \tau_j \, -\frac{4\pi^2}{\l_j}\frac{\n_{\hat g}\mathcal{F}_j^{A}(a_j)}{\mathcal{F}_j^{A}(a_j)}  \,
 - \frac{4\pi^2}{\l_j}\frac{ \frac{ \partial F_i^A (a_j)}{\partial n_{g_j}}}{\mathcal{F}_j^{A}(a_j)}\, \\&+O\left(\sum_{i=1}^{p+q}|\alpha_i-1|^2+\sum_{i=1}^{p+q}\frac{1}{\l_i^2}+\sum_{r=1}^{\bar k}|\beta_r|^2+\sum_{i=1}^{p+q}\tau_i^2\right),
\end{split}
\end{equation*}
where $\hat g=:g_{|\partial M}$, $A:=(a_1, \cdots, a_{p+q})$, $O(1)$ is as in Lemma \ref{eq:energyest}, and for  $i=1, \cdots, p+q$, $\mathcal{F}_i^A$ is as in Lemma \ref{eq:energyest} and \;$\tau_i$ is as in Lemma \ref{eq:gradientlambdaest}.
\end{lem}












\bigskip
\bigskip
\bigskip
\bigskip
\bigskip


\begin{minipage}[l]{10.5cm}
    {\small Mohameden Ahmedou}\\
    {\footnotesize Mathematisches Institut der Justus-Liebig-Universit\"at Giessen}\\
    {\footnotesize Arndtsrasse 2, D-35392 Giessen}\\
    {\footnotesize Germany}\\
    {\tt\footnotesize Mohameden.Ahmedou@math.uni-giessen.de  }\\
    {}
\end{minipage}

$$
\begin{minipage}[l]{7.5cm}
    {\small  Sadok Kallel }\\
    {\footnotesize  American University of Sharjah (UAE)}\\
    {\footnotesize  and    Laboratoire Painlev\'e, USTL(France) }\\
    {\tt\footnotesize  sadok.kallel@math.univ-lille.fr,}\\
    {}
\end{minipage}
   \quad
\begin{minipage}[l]{7.5cm}
    {\small Cheikh Birahim Ndiaye}\\
    {\footnotesize  Mathematisches Institut der Justus-Liebig-Universit\"at Giessen }\\
    {\footnotesize   Arndtsrasse 2, D-35392 Giessen} \\
    {\footnotesize  Germany } \\
    {\footnotesize  and    T\"ubingen University, Auf der Morgenstelle 10, D-72076 T\"ubingen }\\
    {\tt\footnotesize ndiaye@math.uni-tuebingen.de}
\end{minipage}
$$

\end{document}